\documentclass[a4paper]{amsart}

\usepackage[english]{babel}
\usepackage{array}
\usepackage[utf8]{inputenc}
\usepackage{graphicx}
\usepackage{amsmath,latexsym,amsfonts,amssymb,amscd,euscript,epic,eepic,times,amsthm,enumerate}
\usepackage{xypic}
\usepackage{nicefrac}
\usepackage{faktor}
\usepackage{pdfsync}
\usepackage[all]{xy}

\usepackage{hyperref}
\usepackage{lscape}
\usepackage{array}

\newtheorem{theorem}{Theorem}[section]
\newtheorem{lemma}[theorem]{Lemma}
\newtheorem{proposition}[theorem]{Proposition}
\newtheorem{corollary}[theorem]{Corollary}
\newtheorem{definition}[theorem]{Definition}
\newtheorem{remark}[theorem]{Remark}

\newcommand{\bpr}{\begin{proposition}}
\newcommand{\epr}{\end{proposition}}
\newcommand{\btm}{\begin{theorem}}
\newcommand{\etm}{\end{theorem}}
\newcommand{\bco}{\begin{corollary}}
\newcommand{\eco}{\end{corollary}}
\newcommand{\blm}{\begin{lemma}}
\newcommand{\elm}{\end{lemma}}
\newcommand{\bdf}{\begin{definition}}
\newcommand{\edf}{\end{definition}}
\newcommand{\bpm}{\begin{pmatrix}}
\newcommand{\epm}{\end{pmatrix}}
\newcommand{\beq}{\begin{equation}}
\newcommand{\eeq}{\end{equation}}
\newcommand{\bit}{\begin{itemize}}
\newcommand{\eit}{\end{itemize}}
\newcommand{\brm}{\begin{remark}}
\newcommand{\erm}{\end{remark}}

\newcommand{\lie}{\mathfrak}
\newcommand{\nr}{\textnormal}
\newcommand{\wt}{\widetilde}

\newcommand{\ol}{\overline}

\newcommand{\CC}{\mathbb{C}}

\newcommand{\PP}{\mathbb{P}}

\newcommand{\RR}{\mathbb{R}}
\newcommand{\WW}{\mathbb{W}}
\newcommand{\ZZ}{\mathbb{Z}}

\newcommand{\Ee}{\mathcal{E}}
\newcommand{\Ff}{\mathcal{F}}

\newcommand{\Oo}{\mathcal{O}}

\newcommand{\Qq}{\mathcal{Q}}
\newcommand{\Rr}{\mathcal{R}}

\newcommand{\Vv}{\mathcal{V}}

\newcommand{\gG}{\lie{g}}

\newcommand{\lL}{\lie{l}}

\newcommand{\quotient}[2]{{\raisebox{.2em}{\thinspace $#1$}\left / \raisebox{-.15em}{ $#2$}\right.}}
\newcommand{\git}[2]{{\raisebox{.2em}{\thinspace $#1$}\left /\!\!/ \raisebox{-.15em}{$#2$}\right.}}

\newcommand\Quotient[2]{
        \mathchoice
            {
                \text{\raise1ex\hbox{\thinspace $#1$}\Big{/} \lower1ex\hbox{$#2$} \thinspace}%
            }
            {
                #1\,/\,#2
            }
            {
                #1\,/\,#2
            }
            {
                #1\,/\,#2
            }
    }

\newcommand\GIT[2]{
        \mathchoice
            {
                \text{\raise1ex\hbox{\thinspace $#1$}\Big{/}\!\!\!\!\Big{/} \lower1ex\hbox{$#2$} \thinspace}%
            }
            {
                #1\,/\,#2
            }
            {
                #1\,/\,#2
            }
            {
                #1\,/\,#2
            }
    }

\newcommand{\qua}{\thinspace}

\newcommand{\lra}{\longrightarrow}

\DeclareMathOperator{\id}{id}

\DeclareMathOperator{\ad}{ad}

\DeclareMathOperator{\Sym}{Sym}
\DeclareMathOperator{\sym}{{\mathfrak{S}}}

\DeclareMathOperator{\Pic}{Pic}

\DeclareMathOperator{\im}{im}

\DeclareMathOperator{\tr}{tr}
\DeclareMathOperator{\gr}{gr}
\DeclareMathOperator{\aj}{aj}

\DeclareMathOperator{\suma}{sum}
\DeclareMathOperator{\mult}{mul}
\DeclareMathOperator{\depth}{depth}
\DeclareMathOperator{\Sing}{Sing}
\DeclareMathOperator{\codim}{codim}

\DeclareMathOperator{\GL}{GL}
\DeclareMathOperator{\PGL}{PGL}
\DeclareMathOperator{\SL}{SL}
\DeclareMathOperator{\PSL}{PSL}

\DeclareMathOperator{\Ort}{O}
\DeclareMathOperator{\SO}{SO}
\DeclareMathOperator{\Sp}{Sp}
\DeclareMathOperator{\Spin}{Spin}

\DeclareMathOperator{\Hom}{Hom}
\DeclareMathOperator{\End}{End}

\newcommand{\MmM}{\mathcal{M}}
\newcommand{\Mm}{\mathcal{N}}

\title{Higgs bundles over elliptic curves}

\author{Emilio Franco}
\address{Emilio Franco \\ Institut für mathematik \\ Freie Universität von Berlin \\ 
Arnimallee, 3 \\ Berlin D-14195 (Germany)}
\email{emiliofranco@fu-berlin.de}

\author{\'Oscar Garc\'ia-Prada}
\address{\'Oscar Garc\'ia-Prada \\ Instituto de Ciencias Matem\'aticas \\ CSIC-UAM-UC3M-UCM \\ 
Calle Nicol\'as Cabrera, 15 \\ 28049 Madrid (Spain)}
\email{oscar.garcia-prada@icmat.es}

\author{P. E. Newstead}
\address{P. E. Newstead \\ Department of Mathematical Sciences \\ University of Liverpool \\ 
Peach Street \\ Liverpool L69 7ZL (United Kingdom)}
\email{newstead@liverpool.ac.uk}

\begin{document}

\date{\today}



\thanks{First author partially supported by Consejo Superior de Investigaciones Cient\'ificas through JAE-Predoc grant program and German Research Foundation through the project SFB 647. First and second authors partially supported by the Ministerio de Econom\'ia y Competitividad of Spain through Project MTM2010-17717 and Severo Ochoa Excellence Grant. The three authors thank the Isaac Newton Institute in Cambridge - which they visited while preparing the paper - for the excellent conditions provided.}

\begin{abstract}
In this paper we study $G$-Higgs bundles over an elliptic curve when the structure group $G$ is a classical complex reductive Lie group. Modifying the notion of family, we define a new moduli problem for the classification of semistable $G$-Higgs bundles of a given topological type over an elliptic curve and we give an explicit description of the associated moduli space as a finite quotient of a product of copies of the cotangent bundle of the elliptic curve. We construct a bijective morphism from this new moduli space to the usual moduli space of semistable $G$-Higgs bundles, proving that the former is the normalization of the latter. We also obtain an explicit description of the Hitchin fibration for our (new) moduli space of $G$-Higgs bundles and we study the generic and non-generic fibres.
\end{abstract}

\maketitle

\tableofcontents

\section{Introduction}
\label{sc introduction}

A systematic study of vector bundles over elliptic curves was initiated in 1957 by Atiyah \cite{atiyah}, where he describes the set of isomorphism classes of indecomposable vector bundles. After the development of GIT and the introduction by Mumford \cite{mumford} of the notions of stability for vector bundles, Atiyah's results were interpreted as the construction of an isomorphism $M(\GL(n,\CC))_d \cong \Sym^h X$, where $M(\GL(n,\CC))_d$ is the moduli space of semistable vector bundles of rank $n$ and degree $d$ over the elliptic curve $X$ and $h = \gcd(n,d)$. In \cite{ramanathan_stable}, Ramanathan extended the notion of stability to $G$-bundles where $G$ is an arbitrary complex reductive connected Lie group. Schweigert \cite{schweigert}, Friedman, Morgan and Witten \cite{friedman&morgan, friedman&morgan_2, friedman&morgan&witten} and for the topologically trivial case Laszlo \cite{laszlo} gave a description of the moduli space of semistable $G$-bundles with topological invariant $d$ over an elliptic curve $X$ in terms of a quotient 
\beq
\label{eq MG = Z/Gamma}
M(G)_d \cong \quotient{Z_{G,d}}{\Gamma_{G,d}},
\eeq
where $Z_{G,d}$ is the product of a certain number of copies of the curve and $\Gamma_{G,d}$ is a finite group. When $G$ is simple and simply connected with coroot lattice $\Lambda$ and Weyl group $W$, this quotient is $(X \otimes_{\ZZ} \Lambda)/W$. Friedman and Morgan \cite{friedman&morgan}, and Laszlo \cite{laszlo} when the topological type $d$ is trivial, constructed a bijective morphism from $Z_{G,d}/\Gamma_{G,d}$ to $M(G)_d$ which, since $M(G)_d$ is normal, is an isomorphism by Zariski's Main Theorem. Recall that $(X \otimes_{\ZZ} \Lambda)/W$ is isomorphic to a weighted projective space by a result of Looijenga \cite{looijenga} (see also the work of Bernstein-Shvartzman \cite{bernstein&shavartzman}).

\

In this paper we study $G$-Higgs bundles over an elliptic curve for classical complex Lie groups. If $G$ is a complex reductive Lie group, a $G$-Higgs bundle over a smooth projective curve is a pair $(P,\varphi)$ where $P$ is a principal $G$-bundle and $\varphi$, called the Higgs field, is a section of the adjoint bundle $\ad P$ tensored by $K$, the canonical line bundle of the curve. When the structure group $G$ is a classical reductive complex Lie group, there is a bijective correspondence between pairs $(P,\varphi)$ and triples $(E,\theta,\Phi)$ where $E$ is a vector bundle, $\theta$ is a reduction of structure group to $G$ of the $\GL(n,\CC)$-bundle associated to $E$ and $\Phi$ is a $K$-twisted endomorphism of $E$ compatible with the reduction of structure group $\theta$. We shall work with this latter description of Higgs bundles rather than the former.

Hitchin introduced $G$-Higgs bundles and their stability conditions in \cite{hitchin-self_duality_equations}. The existence of the moduli space of semistable Higgs bundles $\MmM(\GL(n,\CC))_d$ was proved by Hitchin in the case of rank $2$, and by Simpson \cite{simpson-hb&ls} and Nitsure \cite{nitsure} in arbitrary rank. In \cite{simpson1, simpson2} Simpson proved the existence of the moduli space $\MmM(G)_d$ of semistable $G$-Higgs bundles when $G$ is an arbitrary complex reductive Lie group. 

Let $\Gamma$ denote the universal central extension by $\ZZ$ of the fundamental group $\pi_1(X)$ of a compact Riemann surface,  and set $\Gamma_\RR = \RR \times_\ZZ \Gamma$. The moduli space of representations of $\Gamma_\RR$ in $G$ with topological type $d$ is the GIT quotient
\[
\Rr(G)_d = \git{\Hom^c(\Gamma_\RR,G)_d}{G},
\]
where $\Hom^c(\Gamma_\RR,G)$ is the space of central representations (i.e. those representations $\rho \in \Hom(\Gamma_\RR,G)$ satisfying $\rho(\RR) \subset Z_G(G)_0$).

As a consequence of a chain of theorems by Narasimhan and Seshadri \cite{narasimhan&seshadri}, Ramana\-than \cite{ramanathan_stable}, Donaldson \cite{donaldson}, Corlette \cite{corlette}, Hitchin \cite{hitchin-self_duality_equations} and Simpson \cite{simpson-hb&ls, simpson1, simpson2}, there exists a homeomorphism $\MmM(G)_d \simeq \Rr(G)_d$.

In \cite{simpson2} Simpson proved the Isosingularity Theorem which implies that $\MmM(G)_0$ is normal if and only if $\Rr(G)_0$ is normal. He proves that $\Rr(\GL(n,\CC))_0$ is normal for compact Riemann surfaces of genus $g \geq 2$, and therefore, in that case, $\MmM(G)_0$ is normal. His proof does not apply for the genus $1$ case but one can use results of Popov \cite{popov} and computations made by the computer program Macaulay (see \cite{hreinsdottir}) to prove that $\Rr(\GL(n,\CC))_0$ and $\MmM(\GL(n,\CC))_0$, hence also $\MmM(\SL(n,\CC))$ and $\MmM(\PGL(n,\CC))_0$, are normal for $n\le4$ (see Section \ref{sc normality}). For the rest of the cases, the normality of the moduli space $\MmM(G)_d$ of Higgs bundles on elliptic curves is a question that remains open.

\

A key result in our study of $G$-Higgs bundles for classical complex Lie groups over an elliptic curve $X$ is that a $G$-Higgs bundle is (semi)stable if and only if the underlying principal bundle is (semi)stable. This is a consequence of the fact that the canonical bundle of an elliptic curve is trivial, i.e. $K \cong \Oo$. Taking the underlying bundle of a semistable Higgs bundle we have a surjective morphism
\[
a_{G,d} : \MmM(G)_d \to M(G)_d. 
\]
Since the fibres of this surjective morphism are connected and so is $M(G)_d$, it follows that $\MmM(G)_d$ is connected.

We obtain an explicit description of semistable, stable and polystable $G$-Higgs bundles over an elliptic curve thanks to the previous result and the description of (semi)stable vector bundles and $G$-bundles given in \cite{atiyah} and \cite{friedman&morgan} respectively. The structure group of a polystable $G$-Higgs bundle can be reduced to a Levi subgroup $L$ of $G$ giving a stable $L$-Higgs bundle. In the elliptic case the conjugacy class of $L$ is the same for every polystable $G$-Higgs bundle with a given topological type. Let $Z_{G,d}$ and $\Gamma_{G,d}$ be as in (\ref{eq MG = Z/Gamma}). Using families of stable $L$-Higgs bundles we can construct families of polystable $G$-Higgs bundles $\Ee$ parametrized by $T^*Z_{G,d}$ such that every polystable $G$-Higgs bundle of topological type $d$ is isomorphic to $\Ee_z$ for some $z \in T^*Z_{G,d}$ and $\Ee_{z_1} \cong \Ee_{z_2}$ if and only if there exists $\gamma \in \Gamma_{G,d}$ giving $z_2 = \gamma \cdot z_1$. This family induces a bijective morphism
\beq
\label{eq Z/Gamma to MmMG}
\quotient{T^*Z_{G,d}}{\Gamma_{G,d}} \stackrel{1:1}{\lra} \MmM(G)_d.
\eeq
If $\MmM(G)_d$ were normal, this bijection would be an isomorphism by Zariski's Main Theorem. However, normality of $\MmM(G)_d$ for $g = 1$ is an open question except in the topologically trivial cases of $G = \GL(n,\CC)$, $G = \SL(n,\CC)$ and $G = \PGL(n,\CC)$ when $n \leq 4$ (see Section \ref{sc normality}). 

In view of this we construct a new moduli functor. The usual moduli functor associates to any scheme $T$ the set of families of $G$-Higgs bundles parametrized by $T$. We will consider a new moduli functor that associates a smaller set of families of $G$-Higgs bundles, the set of {\it locally graded} families (defined in Sections \ref{sc new moduli problem}, \ref{sc MmSLnC and MmPGLnC} and \ref{sc MmSp2mC}). For this new moduli functor the family $\Ee$ of polystable $G$-Higgs bundles constructed above has the local universal property. For the moduli space of Higgs bundles $\Mm(G)_d$ associated to this moduli functor, we have
\beq 
\label{eq MmG = T*Z/Gamma}
\Mm(G)_d \cong \quotient{T^*Z_{G,d}}{\Gamma_{G,d}}.
\eeq
From (\ref{eq Z/Gamma to MmMG}) we observe that there exists a bijective morphism $\Mm(G)_d \to \MmM(G)_d$; thus our new moduli space is not classifying extra structure. Furthermore, since $\Mm(G)_d$ is normal, it is the normalization of $\MmM(G)_d$.


The Hitchin map is defined in \cite{hitchin-duke} by evaluating a basis of the invariant polynomials $q_1, \dots, q_\ell$ on the Higgs field, 
\[
\begin{array}{ccc}
b_{G,d}: & \Mm(G)_d \lra & \bigoplus H^0(X,\Oo^{\otimes r_i}) 
\\
& (E,\Phi) \longmapsto & (q_1(\Phi), \dots, q_\ell(\Phi)).
\end{array}
\]
We observe that $b_{G,d}$ is not surjective in general. In order to preserve the surjectivity of the Hitchin map we redefine for each $d$ the Hitchin base $B(G,d)$ as the image of $b_{G,d}$. The explicit description of the moduli space $\Mm(G)_d$ allows us to study in detail the two fibrations
\[
\xymatrix{
 & \Mm(G)_d \ar[ld]_{a_{G,d}} \ar[rd]^{b_{G,d}} &
\\
M(G)_d & & B(G,d).
}
\]
In particular we describe all the fibres of the Hitchin fibration, not only the generic ones. Two Langlands dual groups have the same Hitchin base. For the two pairs of dual groups, $\SL(n,\CC)$ and $\PGL(n,\CC)$, and $\Sp(2m,\CC)$ and $\SO(2m+1,\CC)$, the Hitchin fibres over a non-generic point of the base are fibrations of projective spaces (in some cases quotients of projective spaces by finite groups) over isomorphic self-dual abelian varieties. 

By means of the quotients (\ref{eq MG = Z/Gamma}) and (\ref{eq MmG = T*Z/Gamma}) we define natural orbifold structures on $M(G)_d$ and $\Mm(G)_d$ and the projection $\Mm(G)_d \stackrel{a}{\lra} M(G)_d$ can be understood as the projection of the orbifold cotangent bundle.

\

The paper is structured as follows:

Section \ref{sc preliminaries} is a review on vector bundles and principal bundles for classical groups over an elliptic curve. It contains the description of stable and polystable bundles derived from \cite{atiyah} and \cite{friedman&morgan} and the subsequent description of the moduli spaces. This section is included not only to set up notation, but also to emphasize the isomorphism (\ref{eq MG = Z/Gamma}) which will be used in the description of the moduli spaces of $G$-Higgs bundles.

In Section \ref{sc definition of Higgs bundles} we give the definitions of $G$-Higgs bundles for classical groups and their stability notions. We discuss normality of the moduli space $\MmM(G)_d$ in Section \ref{sc normality} proving that $\MmM(\GL(n,\CC))_0$, $\MmM(\SL(n,\CC))$ and $\MmM(\PGL(n,\CC))_0$ are normal for $n \leq 4$ [Theorem \ref{tm normality of MmMGLnC_0 when n lower than 5}].

Section \ref{sc description of Mm} contains the explicit description of the moduli spaces of Higgs bundles. In Section \ref{sc stability of higgs bundles} we establish the equivalence between the stability of a Higgs bundle and the stability of its underlying bundle. This fact allows us to give a complete descrition of the polystable Higgs bundles.

We construct in Section \ref{sc MmMGLnC} a family $\Ee_{n,d}$ parametrizing all such bundles. The family $\Ee_{n,d}$ is parametrized by $T^*X \times \stackrel{h}{\dots} \times T^*X$ in such a way that two points parametrize isomorphic polystable Higgs bundles if and only if one is a permutation of the other. Using $\Ee_{n,d}$ we obtain a bijection between the symmetric product $\Sym^h T^*X$ of the cotangent bundle of the curve and the moduli space of Higgs bundles $\MmM(\GL(n,\CC))_d$ [Theorem \ref{tm MmGLnC is the normalization of MmMGLnC}]. We also study the smooth points of $\MmM(\GL(n,\CC))_d$ and its singular locus [Theorem \ref{tm singular set of MmMGLnC}].

We define locally graded families in Section \ref{sc new moduli problem} and we consider the modified moduli problem given by taking the image of the moduli functor to be the set of S-equivalence classes of locally graded families. We prove that $\Ee_{n,d}$ has the local universal property among locally graded families which implies that the moduli space associated to the new moduli functor $\Mm(\GL(n,\CC))_d$ is isomorphic to $\Sym^h T^*X$ [Theorem \ref{tm Moduli de GL-Higgs bundles}].

The work of Sections \ref{sc MmMGLnC} and \ref{sc new moduli problem} allows us to study in Section \ref{sc MmSLnC and MmPGLnC} the moduli spaces $\MmM(\SL(n,\CC))$ and $\MmM(\PGL(n,\CC))_{\wt{d}}$ for the usual moduli problem [Theorems \ref{tm MmMSLnC} and \ref{tm MmMPGLnC}] and $\Mm(\SL(n,\CC))$ and $\Mm(\PGL(n,\CC))_{\wt{d}}$ for the new moduli problem [Theorem \ref{tm MmSLnC and PGLnC}].

In Section \ref{sc MmSp2mC} we study the usual moduli spaces $\MmM(\Sp(2m,\CC))$, $\MmM(\Ort(n,\CC))_{k,a}$ and $\MmM(\SO(n,\CC))_{w_2}$ [Theorem \ref{tm MmG is the normalization of MmMG}]. Following an analogous procedure of that of Section \ref{sc new moduli problem} we obtain an explicit description of $\Mm(\Sp(2m,\CC))$, $\Mm(\Ort(n,\CC))_{k,a}$ and $\Mm(\SO(n,\CC))_{w_2}$ [Theorem \ref{tm MmG}].

We study the Hicthin map for these moduli spaces in Section \ref{sc hitchin map}, and we describe the generic and non-generic fibres explicitly.

Finally in the Appendix we define an action of the groups of torsion points of an abelian variety on a product of copies of the abelian variety. We study properties of this action and its quotient space. These results are used in Section \ref{sc Hitchin fibres of SLnC and PGLnC} to describe the Hitchin fibres for $\PGL(n,\CC)$ [Proposition \ref{pr PGLnC Hitchin fibre} and Remark \ref{rm SLnC and PGLnC Hitchin fibres}] (see also Remarks \ref{rm MPGLnC is also Z/sym_h} and \ref{rm MmPGLnC is also T*Z/sym_h}). 

We work in the category of algebraic schemes over $\CC$. All the bundles con\-si\-de\-red are algebraic bundles. The slope $\mu(E)$ of a vector bundle $E$ of rank $n$ and degree $d$ is defined by $\mu(E) := d/n$.

\

\noindent
{\it Acknowledgements.} This article is a modified version of part of the PhD thesis of the first author prepared under the supervision of the second and third authors at ICMAT (Madrid). The first author wishes to thank the second and third authors for their teaching, help and encouragement.

\section{Review on principal bundles over elliptic curves for classical groups}
\label{sc preliminaries}

\subsection{Vector bundles}
\label{sc vector bundles}

Let $X$ be a smooth projective curve of genus $g = 1$ and let $x_0$ be a distiguished point on it; we call the pair $(X,x_0)$ an {\it elliptic curve}. However, by abuse of notation, we usually refer to the elliptic curve simply as $X$.

The Abel-Jacobi map $\aj_h: \Sym^h X \to \Pic^h(X)$ sends the tuple $[x_1, \dots, x_h]_{\sym_h}$ to the line bundle $L(D)$, where $D$ is the divisor associated to the tuple of points. For $h > 2g - 2 = 0$ the map is surjective and the inverse image of $L \in \Pic^h(X)$ is given by the zeroes of the sections of $L$, i.e. it is the projective space
\beq
\label{eq fibre of aj}
\aj_h^{-1}(L) = \PP H^0(X,L) \cong \PP^{h-1}.
\eeq
For $h = 1$ this inverse image is a point and then $\aj_1 : X \stackrel{\cong}{\lra} \Pic^1(X)$ is an isomorphism. The distiguished point $x_0$ of the elliptic curve gives an isomorphism between $\Pic^{d}(X)$ and $\Pic^{d-h}(X)$,
\beq
\label{eq definition of t_x_0 h}
\begin{array}{cccc}
t^{x_0}_{h} : & \Pic^{d-h}(X)  & \lra &  \Pic^{d}(X) 
\\ 
& L &\longmapsto & L \otimes \Oo(x_0)^{h}.
\end{array}
\eeq

For every $d$ we define an isomorphism 
\beq
\label{eq definition of varsigma_x0 1 0}
\varsigma^{x_0}_{1,d} : X \lra \Pic^d(X),
\eeq
given by $\varsigma^{x_0}_{1,d} = t^{x_0}_{d-1} \circ \aj_1$. In particular $(\varsigma^{x_0}_{1,0})^{-1} : \Pic^0(X) \stackrel{\cong}{\lra} X$ defines an abelian group structure on $X$ with $x_0$ as the identity. The elliptic curve $(X,x_0)$ with this abelian group structure is an abelian variety and the diagram
\beq
\label{eq sum_X is the abel jacobi map}
\xymatrix{
[x_1, \dots, x_h]_{\sym_h} \ar@{|->}[d] & \Sym^h X \ar@{=}[r] \ar[d]_{\suma^h_X} & \Sym^h X \ar[d]^{\aj_h}  &  [x_1, \dots, x_h]_{\sym_h} \ar@{|->}[d]
\\
\sum_{i=1}^h x_i  &  X \ar[r]^{\cong \quad}_{\varsigma^{x_0}_{1,h}\quad } & \Pic^h(X)   &  \Oo(x_1) \otimes \dots  \otimes \Oo(x_h)
}
\eeq
commutes.

The vector bundle $E$ is {\it semistable} if every subbundle $F$ of $E$ satisfies
\[
\mu(F) \leq \mu(E).
\]
The vector bundle is {\it stable} if the above inequality is strict for every proper subbundle and it is {\it polystable} if it decomposes into a direct sum of stable vector bundles, all of the same slope.

Every semistable vector bundle $E$ possesses a {\it Jordan-H\"older filtration}
\[
0 = E_0 \subsetneq E_1 \subsetneq E_2 \subsetneq \dots \subsetneq E_m = E,
\]
where every quotient $E_i/E_{i-1}$ is stable of slope $\mu(E_i/E_{i-1}) = \mu(E)$. The {\it associated graded vector bundle} of $E$ is defined by
\[
\gr E := \bigoplus_i (E_i / E_{i-1}). 
\]
Although the Jordan-H\"older filtration of a given semistable vector bundle $E$ might not be unique, one can prove that the isomorphism class of $\gr E$ is unique. Two semistable vector bundles $E_1$ and $E_2$ are said to be {\it S-equivalent} if $\gr E_1 \cong \gr E_2$. 

A {\it family} of vector bundles over $X$ parametrized by a scheme $Y$ is a vector bundle $\Vv$ over $X \times Y$. We write $\Vv_y:=\Vv|X\times\{y\}$. Given a property of vector bundles which is satisfied by $\Vv_y$ for all $y\in Y$, we shall say that the family $\Vv$ satisfies the property {\it pointwise}. Two families of semistable vector bundles parametrized by the same variety $Y$ will be said to be {\it S-equivalent} if they are pointwise S-equivalent.

The moduli functor that associates to every scheme $Y$ the set of S-equivalence classes of families of semistable vector bundles of rank $n$ and degree $d$ parametrized by $Y$ possesses a coarse moduli space, which we denote by $M(\GL(n,\CC))_d$. (The notation is justified by the fact that there is a bijective correspondence between $\GL(n,\CC)$-bundles and vector bundles of rank $n$.) Every point of the moduli space represents a S-equivalence class of semistable vector bundles (or equivalently an isomorphism class of polystable vector bundles). The moduli space $M^{st}(\GL(n,\CC))_d$ of isomorphism classes of stable vector bundles is a smooth Zariski open subset of $M(\GL(n,\CC))_d$. Note also that, when $\gcd(n,d) = 1$, every semistable vector bundle is stable and then $M^{st}(\GL(n,\CC))_d = M(\GL(n,\CC))_d$.

These moduli spaces were described implicitly by Atiyah \cite{atiyah}, who did not have the notion of stability available. In 1991, Tu \cite{tu} interpreted Atiyah's results to give an explicit description of the moduli spaces (see also \cite{LePotier}). The following properties of vector bundles over elliptic curves are contained in \cite{atiyah} or \cite{tu} (with some changes of notation).

\begin{itemize}
\item If $\gcd(n,d)=1$, 
\begin{itemize}
\item the morphism given by the determinant
\[
\det:M(\GL(n,\CC))_d \stackrel{\cong}{\lra} \Pic^d(X)
\]
is an isomorphism;

\item a stable vector bundle $E$ of rank $n$ and degree $d$ satisfies $E\otimes L\cong E$ if and only if $L$ is a line bundle in $\Pic^0(X)[n]$ (i.e. $L$ is such that $L^{\otimes n}\cong \Oo$);

\item writing $\varsigma^{x_0}_{n,d} =  \det^{-1} \circ\qua \varsigma^{x_0}_{1,d}$ where $\varsigma^{x_0}_{1,d}$ is the map given in (\ref{eq definition of varsigma_x0 1 0}), we have
\beq
\label{eq definition of varsigma in the coprime case}
\varsigma^{x_0}_{n, d} : X \stackrel{\cong}{\lra} M(\GL(n,\CC))_d;
\eeq 

\item there exists a family $\Vv^{x_0}_{n,d}$ of stable vector bundles of rank $n$ and degree $d$ parametrized by $X$ such that  for every $x \in X$,
\beq
\label{eq varsigma Vv_x = x}
\varsigma^{x_0}_{n,d}(x) = [(\Vv^{x_0}_{n,d})_x]_S.
\eeq
Every family $\Ff \to X \times Y$ of semistable (and therefore stable) vector bundles with rank $n$ and degree $d$ defines naturally a morphism $\nu_{\Ff} : Y \to M(\GL(n,\CC))_d$. The composition with $(\varsigma^{x_0}_{n,d})^{-1}$ gives us a morphism $f : Y \to X$, which is canonically defined (up to the choice of $x_0 \in X$). Thanks to (\ref{eq varsigma Vv_x = x}) we know that $f^* \Ff \sim_S \Vv^{x_0}_{n,d}$, so $\Vv^{x_0}_{n,d}$ is a universal family in this sense.
\end{itemize}

\item There exists a unique indecomposable bundle $F_n$ of degree $0$ and rank $n$ such that $H^0(X,F_n) \neq 0$. Moreover $\dim H^0(X,F_n)=1$ and $F_n$ is a multiple extension of copies of $\Oo$. In particular $F_n$ is semistable. 

\item Every indecomposable bundle of degree $0$ and rank $n$ is of the form $F_n \otimes L$ for a unique line bundle $L$ of degree $0$.

\item If $\gcd(n,d)=h>1$,
\begin{itemize}

\item the fibre product over $X$ of $h$ copies of the family $\Vv^{x_0}_{n',d'}$ gives us a family $\Vv^{x_0}_{n,d}$ of polystable vector bundles parametrized by $Z_h = X \times \stackrel{h}{\dots} \times X$;

\item  every indecomposable bundle of rank $n$ and degree $d$ is of the form $E'\otimes F_h$ for a unique stable bundle $E'$ of rank $n'=\frac{n}h$ and degree $d'=\frac{d}h$;
\item every semistable bundle of rank $n$ and degree $d$ is of the form $\bigoplus_{j=1}^s(E_j'\otimes F_{h_j})$, where each $E_j'$ is stable of rank $n'$ and degree $d'$ and $\sum_{j=1}^s h_j=h$;
\item every polystable bundle of rank $n$ and degree $d$ is of the form $E_1'\oplus\ldots\oplus E_h'$, where each $E_i'$ is stable of rank $n'$ and degree $d'$;
\item as a consequence $M^{st}(\GL(n,\CC))_d$ is empty and the map to the moduli space induced by $\Vv^{x_0}_{n,d}$,
\beq
\label{eq cho 0 for GLnC}
\nu_{\Vv^{x_0}_{n,d}} : Z_h = X \times \stackrel{h}{\dots} \times X \lra M(\GL(n,\CC))_d,
\eeq
is surjective and factors through $\Sym^h X$ giving an isomorphism
\beq
\label{eq definition of varsigma in the noncoprime case}
\varsigma^{x_0}_{n,d} : \Sym^h X \stackrel{\cong}{\lra} M(\GL(n,\CC))_d.
\eeq 
\end{itemize}

\item If $E$ is stable, $\End E \cong \bigoplus_{L_i \in \Pic^0(X)[n]} L_i$.

\item $F_n\cong F_n^*$ and $F_n\otimes F_m$ is a direct sum of various $F_\ell$. In particular $\End F_n\cong F_1\oplus F_3\oplus\ldots\oplus F_{2n-1}$.
\end{itemize}

\subsection{Special linear and projective bundles}
\label{sc SLnC and PGLnC-bundles}

A {\it special linear or $\SL(n,\CC)$-bundle} over the elliptic curve $X$ is a pair $(E,\tau)$, where $E$ is a vector bundle and $\tau$ is a never vanishing section of $\det E$. An {\it isomorphism} between the $\SL(n,\CC)$-bundles $(E_1,\tau_1)$ and $(E_2,\tau_2)$ is an isomorphism of vector bundles $f : E_1 \to E_2$ such that $\tau_2 = \det f (\tau_1)$. It follows that $(E,\tau)$ is isomorphic to $(E,1)$. Therefore, a $\SL(n,\CC)$-bundle is completely determined by a vector bundle $E$ with trivial determinant. Note that, since $h = n$, there are no stable $\SL(n,\CC)$-bundles for $n \geq 2$.

A $\SL(n,\CC)$-bundle is {\it semistable or polystable} if it is, respectively, a semistable or polystable vector bundle. Two semistable $\SL(n,\CC)$-bundles are {\it S-equi\-va\-lent} if they are S-equivalent vector bundles. Again, we define S-equivalence for families pointwise. The moduli functor for $\SL(n,\CC)$-bundles under S-equivalence possesses a coarse moduli space $M(\SL(n,\CC))$.

By \cite{tu} and the commutativity of (\ref{eq sum_X is the abel jacobi map}), the diagram 
\[
\xymatrix{
\Sym^h X \ar[d]_{\suma^h_X} \ar[rr]^{\varsigma^{x_0}_{n,d} \qquad}_{\cong \qquad} & & M(\GL(n,\CC))_d \ar[d]^{\det} 
\\
X \ar[rr]^{\cong}_{\varsigma^{x_0}_{1,d}} & & \Pic^d(X).
}
\]
commutes. When $d = 0$, we have $h = \gcd(n,d) = n$, so 
\[
M(\SL(n,\CC)) \cong (\det)^{-1}(\Oo) \cong (\suma^n_X)^{-1}(x_0) = \aj_n^{-1}((\Oo)) \cong \PP^{n-1}.
\]
Take $A_n$ to be the subvariety of $Z_n = X \times \stackrel{n}{\dots} \times X$ given by
\beq
\label{eq definition of A_n}
A_n = \{ (x_1, \dots, x_n) | x_1 + \dots + x_n = x_0 \}.
\eeq
Let $u_n : A_n \to Z_{n-1}$ be the projection on the first $n-1$ factors; this morphism is an isomorphism and its inverse $u_n^{-1}$ sends $(x_1 \dots, x_{n-1})$ to $(x_1 \dots, x_{n-1}, -\sum x_i)$. Since the symmetric group $\sym_n$ preserves $A_n \subset Z_n$, we can use $u_n$ to define an action of $\sym_n$ on $Z_{n-1} = X \times \stackrel{n-1}{\dots} \times X$. This action gives an isomorphism between $Z_{n-1}/\sym_n$ and $A_n/\sym_n$: composing this with the restriction of $\varsigma^{x_0}_{n,0}$ gives
\beq
\label{eq MSLnC}
\hat{\varsigma}^{x_0}_n : \quotient{Z_{n-1}}{\sym_n} = \quotient{X \times \stackrel{n-1}{\dots} \times X}{\sym_n} \stackrel{\cong}{\lra}  M(\SL(n,\CC)).
\eeq

\

On a curve, every projective bundle or $\PGL(n,\CC)$-bundle is the projectivization of a vector bundle and we denote by $\PP(E)$ the projective bundle associated to $E$. It is well known that two vector bundles $E_1$ and $E_2$ give isomorphic projective bundles if and only if there exists a line bundle $L$ such that $E_2 \cong L \otimes E_1$. Note that the projectivization of a vector bundle of rank $n$ and degree $d$ is a $\PGL(n,\CC)$-bundle of degree $\wt{d} = (d \mod n)$. 

A $\PGL(n,\CC)$-bundle $\PP(E)$ is {\it semistable, stable or polystable} if $E$ is respectively a semistable, stable or polystable bundle. When $\PP(E)$ is semistable, we define its {\it associated graded object} as the projectivization $\PP(\gr E)$. Two se\-mista\-ble $\PGL(n,\CC)$-bundles are {\it S-equivalent} if they have isomorphic graded objects. A {\it family of projective bundles} over $X$ parametrized by $Y$ is a projective bundle over $X \times Y$. We define S-equivalence of families of semistable $\PGL(n,\CC)$-bundles pointwise. Let us consider the moduli functor that associates to every scheme $Y$ the set of S-equivalence classes of families of semistable $\PGL(n,\CC)$-bundles of degree $\wt{d} = (d \mod n)$ parametrized by $Y$. There exists a coarse moduli space $M(\PGL(n,\CC))_{\wt{d}}$ associated to this functor. Since they are no stable vector bundles if the rank $n$ and the degree $d$ are not coprime, the stable locus $M^{st}(\PGL(n,\CC))_{\wt{d}}$ is empty in that case, and $M^{st}(\PGL(n,\CC))_{\wt{d}} = M(\PGL(n,\CC))_{\wt{d}}$ if $n$ is coprime to $\wt{d}$.

\brm
\label{rm new family of projective bundles}
\nr{Since the dimension of $X \times Y$ is greater than $1$, not every projective bundle over $X \times Y$ comes from a vector bundle, i.e. not every family of $\PGL(n,\CC)$-bundles comes from a family of vector bundles. If we modify the notion of family of projective bundles, allowing only those that come from families of vector bundles, we obtain a different moduli functor. It can be proved that the two moduli problems have isomorphic coarse moduli spaces. Working with the se\-cond picture, one sees that $M(\PGL(n,\CC))$ is the quotient of $M(\GL(n,\CC))$ by the action of $\Pic(X)$ given by $(L,E) \mapsto L \otimes E$. One can always fix the determinant by tensoring by some element of $\Pic(X)$ and the only elements of $\Pic(X)$ that preserve the determinant are the $n$-th roots of the trivial bundle, so}
\[
M(\PGL(n,\CC))_{\wt{d}} \cong \quotient{\det^{-1}(\Oo(x_0)^{\otimes d})}{\Pic^0(X)[n]}.
\]
\erm

Take $h = \gcd(n,d)$ and set $n' = \frac{n}{h}$ and $d' = \frac{d}{h}$. Consider the following action of $X$ on $\Sym^h X$ with weight $n'$,
\beq
\label{eq definition of the action of X on Sym^h X}
\begin{array}{cccc}
X \times \Sym^h  & \lra &  \Sym^h X 
\\ 
\left( x, [x_1, \dots, x_h]_{\sym_h} \right) &\longmapsto & [x_1 + n'x, \dots, x_h + n' x]_{\sym_h}.
\end{array}
\eeq
By Atiyah's results, the diagram constructed using (\ref{eq definition of the action of X on Sym^h X})
\beq
\label{eq commutivity of the action of X on Sym^h and the action of Pic on MGLnC}
\xymatrix{
X \times \Sym^h X \ar[rr]  \ar[d]_{\varsigma^{x_0}_{1,0} \times \varsigma^{x_0}_{n,d}}^{\cong}  & & \Sym^h X \ar[d]^{\varsigma^{x_0}_{n,d}}_{\cong}
\\
\Pic^0(X) \times M(\GL(n,\CC))_d \ar[rr]^{\qquad  - \otimes - } &  & M(\GL(n,\CC))_d
}
\eeq
commutes. The action of $X[n]$ with weight $n'$ corresponds to the action of $X[h]$ with weight $1$. Recall that $(\varsigma^{x_0}_{n,d})^{-1}(\det^{-1}(\Oo(x_0)^{\otimes d}))$ is $A_h / \sym_h \cong Z_{h-1}/ \sym_h$ and clearly the action (with weight $1$) of $X[h]$ on $A_h$ corresponds naturally to the action (with weight $1$) of $X[h]$ on $Z_{h-1}$. For every $x \in X[h]$, one has that $x = -(h-1) x$, and then the action of $X[h]$ commutes with the action of the symmetric group $\sym_h$. Therefore, we have an isomorphism
\beq
\label{eq MPGLnC}
\check{\varsigma}^{x_0}_{n,\wt{d}} : \quotient{Z_{h-1}}{\sym_h \times X[h]} = \quotient{X \times \stackrel{h-1}{\dots} \times X}{\sym_h \times X[h]} \stackrel{\cong}{\lra} M(\PGL(n,\CC))_{\wt{d}},
\eeq
induced by the restriction of $\varsigma^{x_0}_{n,d}$ to $(\suma^h_{X})^{-1}(x_0)$.

\brm
\label{rm MPGLnC is also Z/sym_h}
\nr{Since $Z_{h-1}/(\sym_h \times X[h]) = (Z_{h-1}/X[h])/\sym_h$, Lemma \ref{lm the quotient of the weighted action is again X times dots times X} implies that}
\[
M(\PGL(n,\CC))_{\wt{d}} \cong \quotient{Z_{h-1}}{\sym_h} = \quotient{X \times \stackrel{h-1}{\dots} \times X}{\sym_h},
\]
\nr{although the action of $\sym_h$ is different from the action used in (\ref{eq MSLnC}). Note that the action of $X[h]$ on $Z_{h-1}$ is free by Lemma \ref{lm the the weighted action is free iff r=1}}. 
\erm 

\brm
\label{rm the action of X r on PP}
\nr{Recalling (\ref{eq sum_X is the abel jacobi map}) one has $(\suma^h_X)^{-1}(x) \cong \PP^{h-1}$ for every $x \in X$. The abelian group structure of $X$ gives an action of $X[h]$ on $\Sym^h X$ that preserves the fibres of $\suma^h_X$ and hence an action of $X[h]$ on $\PP^{h-1}$. Using this action, we have}
\[
M(\PGL(n,\CC))_{\wt{d}} \cong \quotient{\PP^{h-1}}{X[h]}.
\]
\erm

\subsection{Symplectic and orthogonal bundles}
\label{sc symplectic and orthogonal bundles}

A {\it symplectic bundle or $\Sp(2m,\CC)$-bundle} over the elliptic curve $X$ is a pair $(E,\Omega)$, where $E$ is a vector bundle of rank $2m$ over $X$ and $\Omega \in H^0(X,\Lambda^2 E^*)$ is a non-degenerate symplectic form on $E$. The $\Sp(2m,\CC)$-bundles $(E,\Omega)$ and $(E',\Omega')$ are {\it isomorphic} if there exists an isomorphism $f : E' \to E$ such that $\Omega' = f^t \Omega f$.

Similarly, an {\it orthogonal bundle or $\Ort(n,\CC)$-bundle} over $X$ is a pair $(E,Q)$, where $E$ is a vector bundle of rank $n$ and $Q \in H^0(X,\Sym^2 E^*)$ is a non-degenerate symmetric form on $E$. Again, two $\Ort(n,\CC)$-bundles $(E,Q)$ and $(E',Q')$ are {\it isomorphic} if there exists an isomorphism $f : E' \to E$ such that $Q' = f^t Q f$.

A {\it special orthogonal bundle or $\SO(n,\CC)$-bundle} is a triple $(E,Q,\tau)$ such that $(E,Q)$ is a $\Ort(n,\CC)$-bun\-dle and $\tau$ is a trivialization of $\det E$ (a never vanishing section of $\det E$) compatible with $Q$, that is $\tau^2 = (\det Q)^{-1}$. An {\it isomorphism} between the $\SO(n,\CC)$-bundles $(E_1,Q_1,\tau_1)$ and $(E_2,Q_2,\tau_2)$ is an isomorphism of the underlying $\Ort(n,\CC)$-bun\-dles that sends $\tau_1$ to $\tau_2$. Note that the existence of a trivialization of $\det E$ implies that $\det E \cong \Oo$. A direct sum of various $\SO(n_i,\CC)$-bundles is the $\SO(n,\CC)$-bundle given by the $\Ort(n,\CC)$-bundle which is the direct sum of the underlying $\Ort(n_i,\CC)$-bundles plus the trivialization of the determinant induced by those of the $\SO(n_i,\CC)$-bundles.

Symplectic and orthogonal bundles are particular cases of pairs $(E,\Theta)$, where $E$ is a vector bundle and $\Theta: E \to E^*$ is an isomorphism that satisfies $\Theta = b \Theta^t$. If $b = 1$, we have an $\Ort(n,\CC)$-bundle, and if $b = -1$ it is a $\Sp(2m,\CC)$-bundle. 

Given the isomorphism $\Theta: E \to E^*$, for every subbundle $F$ of $E$ we can define its orthogonal complement with respect to $\Theta$,
\[
F^{\perp_{\Theta}} = \{ v \in E | \Theta(v)(u) = 0 \nr{ for every } u \in F \}.
\]
A subbundle is {\it isotropic} with respect to $\Theta$ if $F \subseteq F^{\perp_\Theta}$. It is {\it coisotropic} if $F^{\perp_{\Theta}} \subseteq F$. Since $E \cong E^*$, the exact sequence 
\[
0 \lra F^{\perp_\Theta} \lra E \lra F^* \lra 0
\]
implies that 
\beq
\label{eq deg F^perp = deg F}
\deg(F^{\perp_{\Theta}}) = \deg(F).
\eeq 

Let $(E,\Theta)$ be a $\Sp(2m,\CC)$-bundle or an $\Ort(n,\CC)$-bundle. Note that $\mu(E) = 0$ since $E \cong E^*$. We say that $(E, \Theta)$ is {\it semistable} if and only if, for any isotropic subbundle $F$ of $E$,
\[
\mu(F)  \leq \mu(E) = 0,
\]
and it is {\it stable} if the above inequality is strict for any proper isotropic subbundle. Recall that every parabolic subgroup of $\Ort(n,\CC)$ or $\Sp(2m,\CC)$ may be described as the subgroup that preserves a partial flag of isotropic subspaces in the standard representation (see for instance \cite[Section 23.3]{fulton&harris}). By \cite[Proposition 3.12]{ramanathan1}, if $(E,\Theta)$ is a semistable $\Sp(2m,\CC)$-bundle or $\Ort(n,\CC)$-bundle, there exists a reduction of structure group to a parabolic subgroup giving a {\it Jordan-Hölder filtration}, 
\[
0 = E_{0} \subsetneq E_{1} \subsetneq \dots \subsetneq E_{k-1} \subsetneq  E_k \subseteq E_k^{\perp_\Theta} \subsetneq E_{k-1}^{\perp_\Theta} \subsetneq \dots \subsetneq  E_1^{\perp_\Theta} \subsetneq E_0^{\perp_\Theta} = E,
\] 
where for every $i \leq k$, $E_i/E_{i-1}$ and $E^{\perp_\Theta}_{i-1}/E^{\perp_\Theta}_i$ are stable vector bundles and $\Theta$ induces an isomorphism $\theta_i : E_i/E_{i-1} \stackrel{\cong}{\lra} (E^{\perp_\Theta}_{i-1}/E^{\perp_\Theta}_i)^*$. If $E_k^{\perp_\Theta}/E_k$ is non-zero, $\Theta$ induces a non-degenerate quadratic form $\wt{\Theta}$ on it in such a way that $(E_k^{\perp_\Theta}/E_k,\wt{\Theta})$ is a stable $\Sp(2m',\CC)$ or $\Ort(n',\CC)$-bundle.

For every semistable $\Sp(2m,\CC)$ or $\Ort(n,\CC)$-bundle $(E,\Theta)$ we define its {\it asso\-cia\-ted graded object} 
\[
\gr(E,\Theta) := (E_k^{\perp_\Theta}/E_k,\wt{\Theta}) \oplus \bigoplus_{i = 1}^{k} \left( (E_i/E_{i-1}) \oplus (E^{\perp_\Theta}_{i-1}/E^{\perp_\Theta}_{i}) , \bpm 0 & b\theta_i^t \\   \theta_i & 0 \epm \right),
\]
where $b=-1$ for $\Sp(2m,\CC)$-bundles and $b = 1$ for $\Ort(n,\CC)$-bundles. As happens in the case of vector bundles, the Jordan-H\"older filtration may not be unique but $\gr(E,\Theta)$ is unique up to isomorphism. We say that a semistable $\Sp(2m,\CC)$ or $\Ort(n,\CC)$-bundle $(E,\Theta)$ is {\it polystable} if $(E,\Theta) \cong \gr(E,\Theta)$. The notion of S-equivalence is clear.

The stability notions of principal bundles were introduced by Ramanathan \cite{ramanathan_stable} in terms of the degrees of line bundles constructed with antidominant cha\-rac\-ters applied to the reduction of the structure group of the bundle to parabolic subgroups. By \cite[Remark 3.1]{ramanathan_stable} and \cite[Remark 4.3]{ramanan} the notions of stability, semistability and polystability of $\Ort(n,\CC)$ and $\SO(n,\CC)$-bundles are equivalent to those described above. This statement can be extended to the case of symplectic bundles.

Since $\SO(2,\CC) \cong \CC^*$ is abelian, every $\SO(2,\CC)$-bundle is stable. Whenever $n>2$ we have that $Z(\SO(n,\CC))$ is a subgroup of $Z(\Ort(n,\CC))$ and then by \cite[Proposition 7.1 and Corollary to Theorem 7.1]{ramanathan_stable}, a $\SO(n,\CC)$-bundle is semistable, stable or polystable if it is semistable, stable or polystable as an $\Ort(n,\CC)$-bundle. The graded object of a $\SO(n,\CC)$-bundle is given by the graded object of the underlying $\Ort(n,\CC)$-bundle:
\[
\gr(E,Q,\tau) = (\gr(E,Q),\tau). 
\]

S-equivalence of families of symplectic, orthogonal and special orthogonal bundles is defined pointwise, as in the case of families of vector bundles. The moduli functor associating to any scheme $Y$ the set of S-equivalence classes of families of semistable $\Sp(2m,\CC)$, $\Ort(n,\CC)$ or $\SO(n,\CC)$-bundles parametrized by $Y$ posseses a coarse moduli space, which we denote by $M(\Sp(2m,\CC))$, $M(\Ort(n,\CC))$ or $M(\SO(n,\CC))$. Although Ra\-ma\-na\-than works in \cite{ramanathan1} and \cite{ramanathan2} with principal bundles over smooth projective curves of genus $g \geq 2$, his construction can be extended to the $g=1$ case and therefore these moduli spaces exist and are normal projective varieties.

Let us recall the following well known result about symplectic and othogonal bundles over smooth projective curves of arbitrary genus. 

\bpr 
\label{pr E Theta semistable implies E semistable} 
Let $(E,\Theta)$ be a semistable $\Sp(2m,\CC)$ or $\Ort(n,\CC)$-bundle. Then $E$ is semistable. Let $n>2$ and let $(E,Q,\tau)$ be a semistable $\SO(n,\CC)$-bundle. Then $E$ is semistable. 
\epr

\proof
The proof of the semistability of the underlying vector bundle of a semistable orthogonal bundle is given in \cite[Proposition 4.2]{ramanan} (see also Proposition \ref{pr E Theta Phi semistable implies E Phi semistable} below). The same proof applies to semistable symplectic bundles. Since for $n > 2$ a $\SO(n,\CC)$-bundle is semistable if and only if its underlying orthogonal bundle is semistable, the result can be extended to special orthogonal bundles.
\qed

We denote the elements of $\Pic^0(X)[2]$ by $\Oo (=J_0)$, $J_1$, $J_2$ and $J_3$. 

\btm
\label{tm MSp2mC and MstOnC}
Suppose $n>4$. We have
\begin{align*}
& M^{st}(\Ort(1,\CC)) = \bigsqcup_{a=0}^{3} \{(E^{st}_{1,a}, Q^{st}_1)\},
\\
& M^{st}(\Ort(2,\CC)) = \bigsqcup_{a'=0}^{5} \{(E^{st}_{2,a}, Q^{st}_2)\},
\\
& M^{st}(\Ort(3,\CC)) = \bigsqcup_{a=0}^{3} \{(E^{st}_{3,a}, Q^{st}_3)\},
\\
& M^{st}(\Ort(4,\CC)) = \{(E^{st}_{4,0}, Q^{st}_{4})\},
\\
& M^{st}(\Ort(n,\CC)) = \emptyset,
\end{align*}
where $(E^{st}_{k,a},Q^{st}_k)$ are defined in (\ref{eq definition of E^st_1}), (\ref{eq definition of E^st_2}), (\ref{eq definition of E^st_3}) and (\ref{eq definition of E^st_4}). For every $m>0$,
\[
M^{st}(\Sp(2m,\CC)) = \emptyset.
\]
\etm

\proof
By \cite[Proposition 4.5]{ramanan}, the orthogonal bundle $(E,Q)$ is stable if and only if $(E,Q)$ is an orthogonal direct sum of subbundles $(E_i,Q_i)$ which are mutually nonisomorphic with each $E_i$ stable. The only stable vector bundles with degree $0$ are the line bundles, so $E_i$ are line bundles. Every $\Ort(1,\CC)$-bundle is isomorphic to $(J_a,1)$ so the only possible stable $\Ort(n,\CC)$-bundles are
\beq
\label{eq definition of E^st_1}
(E^{st}_{1,a}, Q^{st}_1) = (J_a,1),  
\eeq
\beq
\label{eq definition of E^st_2}
(E^{st}_{2,a'}, Q^{st}_2) = (J_{b_i},1)  \oplus (J_{b_j},1),  
\eeq
where $a'$ indexes the pairs $\{ i,j\}$ such that $b_i \neq b_j$,
\beq
\label{eq definition of E^st_3}
(E^{st}_{3,a}, Q^{st}_3) = (J_{b_1},1)  \oplus(J_{b_2},1)  \oplus(J_{b_3},1), 
\eeq
where $b_i \neq a$ and $b_i \neq b_j$ if $i \neq j$, and
\beq
\label{eq definition of E^st_4}
(E^{st}_{4,0}, Q^{st}_4) = (J_0,1)  \oplus (J_1,1) \oplus (J_2,1)  \oplus (J_3,1).  
\eeq

The proof of \cite[Proposition 4.5]{ramanan} can be extended to symplectic bundles. Then, since there are no stable vector bundles of even rank and degree $0$, there are no stable $\Sp(2m,\CC)$-bundles.
\qed

\bco
\label{co stable SOnC-bundles}
Let $n \neq 2$. Any stable $\SO(n,\CC)$-bundle is isomorphic to one of
\begin{enumerate}
\item the $\SO(1,\CC)$-bundle $(E^{st}_{1,0},Q^{st}_1,1) = (\Oo,1,1)$,

\item the $\SO(3,\CC)$-bundle $(E^{st}_{3,0},Q^{st}_3,1)$
\item the $\SO(4,\CC)$-bundle $(E^{st}_{4,0},Q^{st}_4,1)$
\end{enumerate}
\eco 

The following is immediate after the definition of the associated graded object, Theorem \ref{tm MSp2mC and MstOnC} and Corollary \ref{co stable SOnC-bundles}. 

\bco 
\label{co E Theta polystable implies E polystable} 
Let $(E,\Theta)$ be a polystable $\Sp(2m,\CC)$, $\Ort(n,\CC)$ or $\SO(n,\CC)$-bundle. Then $E$ is polystable.
\eco

Every $\SO(2,\CC)$-bundle is stable and isomorphic to
\[
(E,Q,\tau) \cong \left( L \oplus L^*, \bpm  & 1 \\ 1 &  \epm, \sqrt{-1} \right), 
\]
where $L$ is a line bundle. The {\it degree} of a $\SO(2,\CC)$-bundle is the degree of $L$.

We can now give a description of the polystable $\Sp(2m,\CC)$, $\Ort(n,\CC)$ and $\SO(n,\CC)$-bundles.

\bpr
\label{pr decomposition of polystable Sp2mC OnC or SOnC-bundles}
A $\Sp(2m,\CC)$-bundle over an elliptic curve is polystable if and only if it is isomorphic to a direct sum of polystable $\Sp(2,\CC)$-bundles. 

An $\Ort(n,\CC)$-bundle is polystable if and only if it is isomorphic to a direct sum of polystable $\Ort(2,\CC)$-bundles or it is isomorphic to a direct sum of polystable $\Ort(2,\CC)$-bundles and a stable $\Ort(m,\CC)$-bundle (where $m = 1,3$ or $4$).

Let $n>2$. Let $(E,Q,\tau)$ be a polystable $\SO(n,\CC)$-bundle with Stiefel-Whitney class $w_2$,
\begin{enumerate}
\item \label{num SO2m with w_2 = 0} if $n=2n'$ and $w_2 = 0$, then $(E,Q,\tau)$ is a direct sum of $m_{n,w_{2}} = n'$ stable $\SO(2,\CC)$-bundles of trivial degree,
\item \label{num SO2m+1 with w_2 = 0} if $n=2n'+1$ and $w_2 = 0$, then $(E,Q,\tau)$ is a direct sum of $(E^{st}_{1,0},Q^{st}_1,1)$ and $m_{n,w_{2}} = n'$ stable $\SO(2,\CC)$-bundles of trivial degree,
\item \label{num SO2m+1 with w_2 = 1} if $n=2n'+1$ and $w_2 = 1$, then $(E,Q,\tau)$ is a direct sum of $(E^{st}_{3,0},Q^{st}_3,1)$ and $m_{n,w_{2}} = n'-1$ stable $\SO(2,\CC)$-bundles of trivial degree,
\item \label{num SO2m with w_2 = 1}if $n=2n'$ and $w_2 = 1$, then $(E,Q,\tau)$ is a direct sum of $(E^{st}_{4,0},Q^{st}_4,1)$ and $m_{n,w_{2}} = n'-2$ stable $\SO(2,\CC)$-bundles of trivial degree.
\end{enumerate}
\epr

\proof
By definition $(E,\Theta)$ is polystable if and only if it is isomorphic to $\gr(E,\Theta)$ and then decomposes as follows,
\[
(E,\Theta) \cong (E_k^{\perp_\Theta}/E_k,\wt{\Theta}) \oplus \bigoplus_{i = 1}^{k} \left( (E_i/E_{i-1}) \oplus (E^{\perp_\Theta}_{i-1}/E^{\perp_\Theta}_{i}) , \bpm 0 & b\theta_i^t \\   \theta_i & 0 \epm \right),
\]
where $(E_k^{\perp_\Theta}/E_k,\wt{\Theta})$ is stable (if it is not zero) and $(E_i/E_{i-1})$ and $(E^{\perp_\Theta}_{i-1}/E^{\perp_\Theta}_{i})$ are stable vector bundles of degree $0$ and therefore have rank $1$. So the factors
\[
\left( (E_i/E_{i-1}) \oplus (E^{\perp_\Theta}_{i-1}/E^{\perp_\Theta}_{i}) , \bpm 0 & b\theta_i^t \\   \theta_i & 0 \epm \right)
\]
are polystable $\Sp(2,\CC)$ or $\Ort(2,\CC)$-bundles. This, together with Theorem \ref{tm MSp2mC and MstOnC}, proves the statement for $\Sp(2m,\CC)$ and $\Ort(n,\CC)$-bundles.

Every $\SO(2,\CC)$-bundle is stable. Recall that for $n>2$, a $\SO(n,\CC)$-bundle is polystable if and only if the underlying $\Ort(n,\CC)$-bundle is polystable. 

Let us take $n>2$. By the description of polystable $\Ort(n,\CC)$-bundles we have given above, a $\SO(n,\CC)$-bundle is polystable if and only if it is a direct sum of stable $\SO(2,\CC)$-bundles of trivial degree and perhaps a $\SO(m,\CC)$-bundle with stable underlying $\Ort(m,\CC)$-bundle. From Theorem \ref{tm MSp2mC and MstOnC} we see that the only possible $\SO(m,\CC)$-bundles with stable underlying $\Ort(m,\CC)$-bundles are the stable $\SO(m,\CC)$-bundles given in Corollary \ref{co stable SOnC-bundles}.

By \cite[Proposition 7.7 and Theorem 7.8]{friedman&morgan}, if $(E,Q,\tau)$ is a semistable $\SO(2n',\CC)$-bundle that lifts to a $\Spin(2n',\CC)$-bundle, then the underlying vector bundle of $\gr(E,Q,\tau)$ is isomorphic to $\bigoplus_{i= 1}^m (L_i \oplus L_i^*)$. On the other hand, if $(E,Q,\tau)$ does not lift to the group $\Spin(2n',\CC)$, its underlying vector bundle is isomorphic to $\Oo \oplus J_1 \oplus J_2 \oplus J_3 \oplus \bigoplus_{i= 1}^{m-2} (L_i \oplus L_i^*)$. This implies that the $\SO(2n',\CC)$-bundles of type \ref{num SO2m with w_2 = 0} lift to $\Spin(2n',\CC)$ and therefore they have trivial Stiefel-Whitney class, while the $\SO(2n',\CC)$-bundles of type \ref{num SO2m with w_2 = 1} do not lift to $\Spin(2n',\CC)$ and they have non-trivial Stiefel-Whitney class.

The odd case is analogous.
\qed

\btm
\label{tm MstSOnC}
We have
\begin{align*}
M^{st}(\SO(1,\CC))_{0} & = \{  [(E^{st}_{1,0},Q^{st}_1,1)]_{\cong}  \},
\\
M^{st}(\SO(1,\CC))_{1} & = \emptyset,
\\
M^{st}(\SO(2,\CC))_d & = M(\SO(2,\CC))_d \cong X,
\\
M^{st}(\SO(3,\CC))_{0} & = \emptyset,
\\
M^{st}(\SO(3,\CC))_{1} & = \{ [(E^{st}_{3,0},Q^{st}_3,1)]_{\cong}  \},
\\
M^{st}(\SO(4,\CC))_0 & = \emptyset,
\\
M^{st}(\SO(4,\CC))_{1} & = \{   [(E^{st}_{4,0},Q^{st}_4,1)]_{\cong}   \}.
\end{align*}

Suppose $n>4$ and let $w_2$ be either $0$ or $1$, then
\[
M^{st}(\SO(n,\CC))_{w_2} = \emptyset.
\]
\etm
 
\proof
This is a straightforward consequence of Corollary \ref{co stable SOnC-bundles} and Proposition \ref{pr decomposition of polystable Sp2mC OnC or SOnC-bundles}. The description of $M(\SO(2,\CC))$ follows from the isomorphism of Lie groups $\SO(2,\CC) \cong \CC^*$.
\qed

\blm
\label{lm E1 Theta1 cong E2 Theta2 iff E1 cong E2}
Let $(E_1,\Theta_1)$ and $(E_2,\Theta_2)$ be two polystable $\Sp(2m,\CC)$ or $\Ort(n,\CC)$-bundles. If $E_1 \cong E_2$, then $(E_1,\Theta_1) \cong (E_2,\Theta_2)$.
\elm

\proof
If $(E_j,\Theta_j)$ are polystable $\Ort(n,\CC)$-bundles, then by Proposition \ref{pr decomposition of polystable Sp2mC OnC or SOnC-bundles} they have, up to isomorphism, the form
\[
(E_j,\Theta_j) \cong \bigoplus_i (J_{a_{j,i}}, 1) \oplus \bigoplus_k \left( L_{j,k} \oplus L_{j,k}^*, \bpm & 1 \\ 1 & \epm \right).
\]
If $E_1 \cong E_2$, we have, after possible reordering, that $J_{a_{1,i}} = J_{a_{2,i}}$ and $L_{1,k} = L_{2,k}$ or $L_{1,k} = L^*_{2,k}$. Since 
\[
\left(L_{j,k} \oplus L_{j,k}^*, \bpm & 1 \\ 1 & \epm \right) \cong \left( L^*_{j,k} \oplus L_{j,k}, \bpm & 1 \\ 1 & \epm \right),
\]
it follows that $(E_1,\Theta_1)$ and $(E_2,\Theta_2)$ are isomorphic $\Ort(n,\CC)$-bundles. The statement for $\Sp(2m,\CC)$-bundles follows from the discussion above and the fact that 
\[
\left( L_{j,k} \oplus L_{j,k}^*, \bpm & -1 \\ 1 & \epm \right) \cong \left( L^*_{j,k} \oplus L_{j,k}, \bpm & -1 \\ 1 & \epm \right),
\]
where the isomorphism is given by
\[
\bpm \sqrt{-1} & \\ & -\sqrt{-1} \epm.
\]
\qed

\blm
\label{lm E1 Theta1 tau1 cong E2 Theta2 tau2 iff E1 Theta1 cong E2 Theta2}
Let $(E_1, Q_1, \tau_1)$ and $(E_2, Q_2, \tau_2)$ be two polystable $\SO(n,\CC)$-bundles with invariants $(n,w_2)$ equal to $(2n', 1)$, $(2n'+1,0)$ or $(2n'+1,1)$. If $(E_1, Q_1) \cong (E_2, Q_2)$, then $(E_1, Q_1, \tau_1) \cong (E_2, Q_2, \tau_2)$.
\elm

\proof
By Proposition \ref{pr decomposition of polystable Sp2mC OnC or SOnC-bundles}, a $\SO(n,\CC)$-bundle with invariants $(n,w_2)$ equal to $(2n', 1)$, $(2n'+1,0)$ or $(2n'+1,1)$ has the form
\[
(E_j, Q_j, \tau_j) \cong  (E^{st}_{k,0}, Q^{st}_k, 1) \oplus \bigoplus_i  (E_{i,j}, Q_{i,j}, \tau_{i,j}), 
\]
where $k = 1$, $3$ or $4$ and the $(E_{i,j}, Q_{i,j}, \tau_{i,j})$ are stable $\SO(2,\CC)$-bundles of trivial degree.

Recall that an isomorphism of $\SO(n,\CC)$-bundles is an isomorphism of $\Ort(n,\CC)$-bundles that preserves the trivialization $\tau$. If $(E_1,Q_1)\cong (E_2,Q_2)$, then the $\Ort(2,\CC)$-bundles $(E_{1,j},Q_{1,j})$ and $(E_{2,j},Q_{2,j})$ are isomorphic, possibly after reordering of the factors. For this order, either 
\[
(E_{1,j},Q_{1,j},\tau_{1,j}) \cong (E_{2,j},Q_{2,j},\tau_{2,j}),
\]
or 
\[
(E_{1,j},Q_{1,j},\tau_{1,j}) \cong (E_{2,j},Q_{2,j},-\tau_{2,j}).
\]
On the other hand, for $k=1$, $3$ and $4$, we have
\[
(E^{st}_{k,0}, Q^{st}_k, 1) \cong (E^{st}_{k,0}, Q^{st}_k, -1)
\]
since both $\SO(k,\CC)$-bundles are stable and by Theorem \ref{tm MstSOnC} there a unique stable $\SO(k,\CC)$-bundle up to isomorphism.

As a consequence, we can construct a morphism that inverts the trivialization $\tau_j$ combined with the morphism that inverts the trivialization of the stable factor. This leaves the trivialization $\tau$ of the total $\SO(n,\CC)$-bundle unchanged.
\qed

\blm
\label{lm polystable SOnC-bundles}
Let $(E_1,Q_1,\tau_1)$ and $(E_2,Q_2,\tau_2)$ be two polystable $\SO(2m,\CC)$-bundles with $w_2 = 0$ of the form 
\[
(E_i,Q_i,\tau_i) \cong \bigoplus_{j=1}^m (E_{i,j},Q_{i,j},\tau_{i,j})
\]
where the $(E_{i,j},Q_{i,j},\tau_{i,j})$ are $\SO(2,\CC)$-bundles of degree $0$. The $\SO(n,\CC)$-bun\-dles are isomorphic if and only if one can form $m$ pairs of the form 
\[
((E_{1,i_{\ell}}, Q_{1,i_{\ell}}, \tau_{1,i_{\ell}}), (E_{2,j_{\ell}}, Q_{2,j_{\ell}}, \tau_{2,j_{\ell}}))
\]
such that for every $\ell$, we have $(E_{1,i_{\ell}},Q_{1,i_{\ell}}) \cong (E_{2,j_{\ell}},Q_{2,j_{\ell}})$ and the number of pairs such that 
\[
(E_{1,i_{\ell}},Q_{1,i_{\ell}}, \tau_{1,i_{\ell}}) \ncong (E_{2,j_{\ell}},Q_{2,j_{\ell}},\tau_{2,j_{\ell}}) 
\]
is even.
\elm

\proof
If $(E_{1},Q_{1})$ and $(E_{2},Q_{2})$ are isomorphic $\Ort(2m,\CC)$-bundles, one can form $m$ pairs of isomorphic $\Ort(2,\CC)$-bundles $(E_{1,i_{\ell}}, Q_{1,i_{\ell}})$ and $(E_{2,j_{\ell}}, Q_{2,j_{\ell}})$. The associated $\SO(2,\CC)$-bundles, $(E_{1,i_{\ell}}, Q_{1,i_{\ell}}, \tau_{1,i_{\ell}})$ and $(E_{2,j_{\ell}}, Q_{2,j_{\ell}}, \tau_{2,j_{\ell}})$, satisfy
\[
(E_{1,i_{\ell}},Q_{1,_{\ell}},\tau_{1,_{\ell}}) \cong (E_{2,j_{\ell}},Q_{2,j_{\ell}},\tau_{2,j_{\ell}})
\]
or 
\[
(E_{1,i_{\ell}},Q_{1,_{\ell}},\tau_{1,_{\ell}}) \cong (E_{2,j_{\ell}},Q_{2,j_{\ell}},-\tau_{2,j_{\ell}}).
\]
In the second situation the $\SO(2,\CC)$-bundles are not isomorphic unless $L_j \cong L_j^*$. If we have an isomorphism of $\Ort(2m,\CC)$-bundles that inverts an even number of $\tau_j$ then the product of all of them remains unchanged and then the $\SO(2m,\CC)$-bundles are isomorphic.

If our isomorphism of $\Ort(2m,\CC)$-bundles inverts an odd number of $\tau_j$, then the product of all of them changes its sign and then the $\SO(2m,\CC)$-bundles cannot be isomorphic.
\qed

Recall the universal family of line bundles $\Vv^{x_0}_{1,0}$ of degree $0$. Let us note that $\Lambda^2 (\Vv_{1,0}^{x_0} \oplus (\Vv^{x_0}_{1,0})^*) \cong \Vv_{1,0}^{x_0} \otimes (\Vv_{1,0}^{x_0})^* \cong \Oo_{X \times X}$ and take the non-vanishing section of $\Lambda^2 (\Vv_{1,0}^{x_0} \oplus (\Vv_{1,0}^{x_0})^*)$, 
\[
\varOmega = \bpm 0 & -1 \\ 1 & 0 \epm.
\]
Then
\[
\wt{\Vv}_2^{x_0} = \left( \Vv_{1,0}^{x_0} \oplus (\Vv^{x_0}_{1,0})^*,\varOmega \right)
\]
is a family of polystable $\Sp(2,\CC)$-bundles parametrized by $X$. If, instead of $\varOmega$, we take the non-vanishing section of $\Sym^2 (\Vv_{1,0}^{x_0} \oplus (\Vv^{x_0}_{1,0})^*)$, 
\[
\Qq = \bpm 0 & 1 \\ 1 & 0 \epm,
\]
we obtain a family of polystable $\Ort(2,\CC)$-bundles parametrized by $X$,
\[
\mathring{\Vv}_2^{x_0} = \left( \Vv_{1,0}^{x_0} \oplus (\Vv^{x_0}_{1,0})^*,\Qq \right).
\]
We see that $(\det \Qq)^{-1}$ is the section $-1$ of $\Oo_{X \times X}$. Then the section $\tau$ of $\det (\Vv_{1,0}^{x_0} \otimes (\Vv_{1,0}^{x_0})^*)$ can be taken to be the imaginary number $\sqrt{-1}$ or $-\sqrt{-1}$. We fix $\tau = \sqrt{-1}$ and we construct the following family of $\SO(2,\CC)$-bundles of degree $0$
\[
\overline{\Vv}_{2}^{x_0} = \left( \Vv_{1,0}^{x_0} \oplus (\Vv_{1,0}^{x_0})^*,\Qq, \sqrt{-1} \right).
\]

\brm
\label{rm Vv2_1 cong Vv2_2 if they are related by pm}
\nr{The restriction of $\wt{\Vv}_2^{x_0}$ and $\mathring{\Vv}_2^{x_0}$ to two different points of $X$, $x_1$ and $x_2$, give S-equivalent (isomorphic) bundles $(\wt{\Vv^{x_0}}_2)_{x_1} \sim_S (\wt{\Vv}_2^{x_0})_{x_2}$ and $(\mathring{\Vv}_2^{x_0})_{x_1} \sim_S (\mathring{\Vv}_2^{x_0})_{x_2}$ if and only if $x_1 = - x_2$}.
\erm

Now we define $\wt{\Vv}_{2m}^{x_0}$ to be the family of polystable $\Sp(2m,\CC)$-bundles induced by taking the fibre product of $m$ copies of $\wt{\Vv}_2^{x_0}$, which is parametrized by $Z_m = X \times \stackrel{m}{\dots} \times X$. By Proposition \ref{pr decomposition of polystable Sp2mC OnC or SOnC-bundles} this family includes representatives of all S-equivalence classes of semistable $\Sp(2m,\CC)$-bundles. It follows that $M(\Sp(2m,\CC))$ is connected.

Recall the stable $\Ort(k,\CC)$-bundle $(E^{st}_{k,a},Q^{st}_{k})$ appearing in Theorem \ref{tm MSp2mC and MstOnC}. Define $\mathring{\Vv}^{x_0}_{n,k,a}$ to be the family of polystable $\Ort(n,\CC)$-bundles given by the direct sum of $(E^{st}_{k,a},Q^{st}_{k})$ and $m_k = \frac{n-k}{2}$ copies of $\mathring{\Vv}_2^{x_0}$. This family is parametrized by $Z_{m_k} = X \times \stackrel{m_k}{\dots} \times X$.

Analogously, for all values of $(n,w_2)$, we define families of polystable $\SO(n,\CC)$-bundles $\overline{\Vv}_{n,w_2}^{x_0}$ as follows:
\bit
\item $\ol{\Vv}_{2m,0}^{x_0}$ is given by $m$ copies of $\ol{\Vv}_2^{x_0}$ and therefore it is parametrized by $Z_m = X \times \stackrel{m}{\dots} \times X$,
\item $\ol{\Vv}_{2m+1,0}^{x_0}$ is given by the direct sum of $(E^{st}_{1,0},Q^{st}_1,1)$ and $m$ copies of $\ol{\Vv}_2^{x_0}$ and therefore it is parametrized by $Z_m = X \times \stackrel{m}{\dots} \times X$,
\item $\ol{\Vv}_{2m+1,1}^{x_0}$ is given by the direct sum of $(E^{st}_{3,0},Q^{st}_3,1)$ and $m-1$ copies of $\ol{\Vv}_2^{x_0}$ and therefore it is parametrized by $Z_{m-1} = X \times \stackrel{m-1}{\dots} \times X$,
\item $\ol{\Vv}_{2m,1}^{x_0}$ is given by the direct sum of $(E^{st}_{4,0},Q^{st}_4,1)$ and $m-2$ copies of $\ol{\Vv}_2^{x_0}$ and therefore it is parametrized by $Z_{m-2} = X \times \stackrel{m-2}{\dots} \times X$.
\eit 

It follows from Proposition \ref{pr decomposition of polystable Sp2mC OnC or SOnC-bundles} that $M(\SO(n,\CC))$ has two connected components, denoted by $M(\SO(n,\CC))_{w_2}$, where $w_2 = 0,1$.

The symmetric group $\sym_m$ acts naturally on $(\ZZ_2 \times \stackrel{m}{\dots} \times \ZZ_2)$ permuting the factors. Using this action we define $\Gamma_m$ as the semidirect product 
\beq
\label{eq definition of Gamma_m}
\Gamma_m := (\ZZ_2 \times \stackrel{m}{\dots} \times \ZZ_2) \rtimes \sym_m
\eeq
determined by the commutation relations $\sigma \overline{c} = (\sigma \cdot \overline{c}) \sigma$, for any $\sigma \in \sym_m$ and any $\overline{c} \in (\ZZ_2 \times \stackrel{m}{\dots} \times \ZZ_2)$.

The permutation action of the symmetric group and the action of $\ZZ_2$ on $X$ given by $-1 \cdot x = -x$ induce an action of $\Gamma_m$ on $Z_m = X \times \stackrel{m}{\dots} \times X$. The quotient of this space by $\Gamma_m$ under this action is
\[
\quotient{Z_m}{\Gamma_m} = \quotient{X \times \stackrel{m}{\dots} \times X}{\Gamma_m} = \Sym^m(X/\ZZ_2) \cong \Sym^m (\PP^1) \cong \PP^m. 
\]

Let $D_m$ be the subgroup of $\ZZ_2 \times \stackrel{m}{\dots} \times \ZZ_2$ given by the tuples $\overline{c} = (c_i, \dots, c_m)$ such that only an even number of $c_i$ are equal to $-1$. We recall that $\Gamma_m$ is the semidirect product of $(\ZZ_2 \times \stackrel{m}{\dots} \times \ZZ_2)$ and $\sym_{m}$. We define $\Delta_m \subset \Gamma_m$ as the subgroup 
\beq
\label{eq definition of Delta_m}
\Delta_m: = \{ \sigma \overline{c} \in \Gamma_m \nr{ such that } \overline{c} \in D_m \}.
\eeq
The action of $\Gamma_m$ on $X \times \stackrel{m}{\dots} \times X$ induces an action of $\Delta_m$.

\brm
\label{rm wtVv_z1 S-equivalent to wtVv_z2 if z1 and z2 are related by Gamma}
\nr{Consider the families $\wt{\Vv}_{2m}^{x_0}$, $\mathring{\Vv}_{n,k,a}^{x_0}$ and $\ol{\Vv}_{n,w_2}^{x_0}$ when $(n,w_2)$ is $(2m,1)$, $(2m+1,0)$ or $(2m+1,1)$. These families are parametrized by $Z_{m'} = X \times \stackrel{m'}{\dots} \times X$ for some $m'$. By Remark \ref{rm Vv2_1 cong Vv2_2 if they are related by pm} and Lemma \ref{lm E1 Theta1 tau1 cong E2 Theta2 tau2 iff E1 Theta1 cong E2 Theta2} we have that any two points $z_1, z_2 \in Z_{m'}$ parametrize isomorphic (S-equivalent) bundles if and only if $z_1 = \gamma \cdot z_2$ for some $\gamma \in \Gamma_m$.}
\erm

\brm
\label{rm ovlineVv_z1 sim ovlineVv_z2 if they are related by the action of Delta}
\nr{Let $z_1,z_2 \in Z_m = X \times \stackrel{m}{\dots} \times X$. It follows by Lemma \ref{lm polystable SOnC-bundles} that $(\overline{\Vv}_{2m,0}^{x_0})_{z_{1}}$ is isomorphic (S-equivalent) to $(\overline{\Vv}_{2m,0}^{x_0})_{z_{2}}$ if and only if there exists $\gamma \in \Delta_m$ such that $\gamma \cdot z_1 = z_2$.}
\erm

\bpr
\label{pr connected components of MOnC}
The connected components of $M(\Ort(2m+1,\CC))$ are indexed by $k=1,3$ and $a = 0, \dots, n_k-1$ where $n_1 = 4$, and $n_3 = 4$.

The connected components of $M(\Ort(2,\CC))$ are indexed by $k=0,2$ and $a = 0, \dots, n_k-1$ where $n_0 = 1$ and $n_2 = 6$.

If $m > 1$, the connected components of $M(\Ort(2m,\CC))$ are indexed by $k=0,2,4$ and $a = 0, \dots, n_k-1$ where $n_0 = 1$, $n_2 = 6$ and $n_4 = 1$.
\epr

\proof
Since the families $\mathring{\Vv}^{x_0}_{n,k,a}$ are parametrized by the connected variety $Z_{(n-k)/2}$, all the S-equivalence classes of semistable $\Ort(n,\CC)$-bundles parametrized by $\mathring{\Vv}^{x_0}_{n,k,a}$ lie in the same connected component of $M(\Ort(n,\CC))$. 

On the other hand if $(E^{st}_{k,a},Q^{st}_{k}) \ncong (E^{st}_{k',a'},Q^{st}_{k'})$, we observe that there is no family of semistable $\Ort(n,\CC)$-bundles connecting an S-equivalence class parametrized by $\mathring{\Vv}^{x_0}_{n,k,a}$ and an S-equivalence class parametrized by $\mathring{\Vv}^{x_0}_{n,k',a'}$.
\qed

\btm
\label{tm MG}
Denote the connected components of $M(\Ort(n,\CC))$ by $M(\Ort(n,\CC))_{k,a}$ where $k,a$ are as in Proposition \ref{pr connected components of MOnC}. There are natural isomorphisms
\[
\wt{\varsigma}^{x_0}_m : \quotient{Z_m}{\Gamma_m} = \Sym^m(X/\ZZ_2) \stackrel{\cong}{\lra} M(\Sp(2m,\CC)),
\]
\[
\mathring{\varsigma}^{x_0}_{n,k,a} : \quotient{Z_{(n-k)/2}}{\Gamma_{(n-k)/2}} = \Sym^{(n-k)/2}(X/\ZZ_2) \stackrel{\cong}{\lra} M(\Ort(n,\CC))_{k,a},
\]
\[
\ol{\varsigma}^{x_0}_{2m,0} : \quotient{Z_{m}}{\Delta_{m}} \stackrel{\cong}{\lra} M(\SO(2m,\CC))_{0}.
\]
For $(n,w_2) = (2m+1,0)$, $(2m+1,1)$ and $(2m,1)$ we set respectively $m' = m$, $m-1$ and $m-2$. There are natural isomorphisms
\[
\ol{\varsigma}^{x_0}_{n,w_2} : \quotient{Z_{m'}}{\Gamma_{m'}} = \Sym^{m'}(X/\ZZ_2) \stackrel{\cong}{\lra} M(\SO(n,\CC))_{w_2},
\]
\etm

\proof
The family $\wt{\Vv}^{x_0}_{m}$ induces a morphism
\[
\nu_{\wt{\Vv}^{x_0}_{m}} : Z_m = X \times \stackrel{m}{\dots} \times X  \lra M(\Sp(2m,\CC)),
\]
and by Remark \ref{rm wtVv_z1 S-equivalent to wtVv_z2 if z1 and z2 are related by Gamma} it factors through $\wt{\varsigma}^{x_0}_m$. By Zariski's Main Theorem, this map is an isomorphism since it is bijective and $M(\Sp(2m,\CC))$ is normal.

Using $\mathring{\Vv}^{x_0}_{n,k,a}$ and $\ol{\Vv}^{x_0}_{n,w_2}$, we can apply the same construction to describe $M(\Ort(n,\CC))_{k,a}$ and $M(\SO(n,\CC))_{w_2}$.
\qed

\section{Higgs bundles over elliptic curves}
\label{sc definition of Higgs bundles}

\subsection{Higgs bundles}
\label{sc vbhb}

A {\it Higgs bundle} over an elliptic curve $X$ is a pair $(E,\Phi)$, where $E$ is a vector bundle on $X$ and $\Phi$ is an endomorphism of $E$ called the {\it Higgs field}. 
Two Higgs bundles, $(E_1,\Phi_1)$ and $(E_2,\Phi_2)$, are {\it isomorphic} if there exists an isomorphism of vector bundles $f: E_1 \to E_2$ such that $\Phi_2 = f \circ \Phi_1 \circ f^{-1}$.

Given the Higgs bundle $(E,\Phi)$, we say that a subbundle $F \subset E$ is {\it $\Phi$-invariant} if $\Phi(F)$ is contained in $F$. A Higgs bundle $(E,\Phi)$ is {\it semistable} if the slope of any $\Phi$-invariant subbundle $F$ satisfies
\[
\mu(F) \leq \mu(E).
\]
The Higgs bundle is {\it stable} if the above inequality is strict for every proper $\Phi$-invariant subbundle and {\it polystable} if it is semistable and isomorphic to a direct sum of stable Higgs bundles.

If $(E,\Phi)$ is semistable, then it has a {\it Jordan-H\"older filtration} of $\Phi$-invariant subbundles
\[
0 = E_0 \subsetneq E_1 \subsetneq E_2 \subsetneq \dots \subsetneq E_m = E
\]
where the restriction of the Higgs field to every quotient $E_i/E_{i-1}$ induces a stable Higgs bundle $(E_i/E_{i-1}, \wt{\Phi_i})$ with slope $\mu(E_i/E_{i-1}) = \mu(E)$. For every semistable Higgs bundle $(E,\Phi)$ we define its {\it associated graded object}
\[
\gr(E,\Phi) := \bigoplus_i (E_i/E_{i-1}, \wt{\Phi_i}).
\]
Although the Jordan-H\"older filtration may not be uniquely determined by $(E,\Phi)$, the isomorphism class of $\gr(E,\Phi)$ is. Two semistable Higgs bundles $(E,\Phi)$ and $(E',\Phi')$ are said to be {\it S-equivalent} if $\gr(E,\Phi) \cong \gr(E',\Phi')$. Denote by $[(E,\Phi)]_S$ the S-equivalence class of $(E,\Phi)$. 

A {\it family of Higgs bundles} parametrized by $Y$ is a pair $\Ee = (\Vv, \varPhi)$, where $\Vv$ is a family of vector bundles parametrized by $Y$ and $\varPhi$ is a section of $\End \Vv$. For every $y \in Y$, we will write $\Ee_y$ for the Higgs bundle over $X$ obtained by restricting $\Ee$ to $X \times \{ y \}$. Two families of semistable Higgs bundles are {\it S-equivalent} if they are pointwise S-equivalent. 

Consider the moduli functor that associates to every scheme $Y$ the set of S-equivalence classes of families of semistable Higgs bundles parametrized by $Y$. By \cite{nitsure} and \cite{simpson-hb&ls} there exists a coarse moduli space $\MmM(\GL(n,\CC))_d$ of S-equivalence classes of semistable Higgs bundles of rank $n$ and degree $d$ associated to this moduli functor. The points of $\MmM(\GL(n,\CC))_d$ correspond to S-equivalence classes of semistable Higgs bundles and can be identified also with isomorphism classes of polystable Higgs bundles since in every S-equivalence class there is always a polystable Higgs bundle which is unique up to isomorphism. The locus of stable Higgs bundles $\MmM^{st}(\GL(n,\CC))_d$ is an open subvariety of $\MmM(\GL(n,\CC))_d$.

\subsection{Special linear and projective Higgs bundles}
\label{sc SLnC and PGLnC}

A {\it special linear or $\SL(n,\CC)$-Higgs bundle} over the elliptic curve $X$ is a triple $(E,\Phi,\tau)$ where $(E,\tau)$ is a $\SL(n,\CC)$-bundle and $\Phi$ is an endomorphism of $E$ such that $\tr \Phi = 0$. As we saw in Section \ref{sc SLnC and PGLnC-bundles}, we can forget about $\tau$ and then a {\it $\SL(n,\CC)$-Higgs bundle} over $X$ is a Higgs bundle $(E,\Phi)$ with trivial determinant and traceless Higgs field, i.e. $\det E \cong \Oo$, $\tr \Phi = 0$. Two $\SL(n,\CC)$-Higgs bundles are {\it isomorphic} if they are isomorphic Higgs bundles.

A $\SL(n,\CC)$-Higgs bundle is {\it semistable, stable or polystable} if it is, respectively, a semistable, stable or polystable Higgs bundle. Two semistable $\SL(n,\CC)$-Higgs bundles are {\it S-equi\-va\-lent} if they are S-equivalent Higgs bundles, and two families of $\SL(n,\CC)$-Higgs bundles are {\it S-equivalent} if they are pointwise S-equi\-va\-lent. Using the moduli problem for Higgs bundles, we define the moduli functor for $\SL(n,\CC)$-Higgs bundles by restricting to trivial determinant and traceless Higgs fields. We denote by $\MmM(\SL(n,\CC))$ the coarse moduli space associated to this moduli problem and by $\MmM^{st}(\SL(n,\CC))$ the stable locus, which is a Zariski open subset of $\MmM(\SL(n,\CC))$. 

Take the morphism $(\det,\tr)$ from $\MmM(\GL(n,\CC))_0$ to $\Pic^0(X) \times H^0(X,\Oo)$ and note that 
\beq
\label{eq MmSLnC is the preimage of det tr}
\MmM(\SL(n,\CC)) \cong (\det, \tr)^{-1}(\Oo,0) \subset \MmM(\GL(n,\CC))_0. 
\eeq

A {\it $\PGL(n,\CC)$-Higgs bundle} over the elliptic curve $X$ is a pair $(\PP(E),\Phi)$ where $\PP(E)$ is a projective bundle given by the vector bundle $E$ and $\Phi$ is an element of $H^0(X,\End E)$ with $\tr \Phi = 0$. There exists a natural isomorphism between $\End E$ and $\End (L \otimes E)$ and then an isomorphism between the $\PGL(n,\CC)$-Higgs bundles $(\PP(E_1),\Phi_1)$ and $(\PP(E_2),\Phi_2)$ corresponds to an isomorphism $f: E_1 \stackrel{\cong}{\lra} E_2 \otimes L$, for some $L \in \Pic(X)$, such that $f \circ \Phi_1 \circ f^{-1} = \Phi_2 \otimes \id_L$.

A $\PGL(n,\CC)$-Higgs bundle $(\PP(E),\Phi)$ is {\it semistable, stable or polystable} if $(E,\Phi)$ is respectively a semistable, stable or polystable Higgs bundle. If $(\PP(E),\Phi)$ is semistable and $\gr(E,\Phi) = (\gr E,\gr \Phi)$, we define {\it the associated graded object} of $(\PP(E),\Phi)$ as the pair $(\PP(\gr E), \gr \Phi)$. Two semistable $\PGL(n,\CC)$-Higgs bundles are S-equivalent if they have isomorphic graded objects. 


Following Remark \ref{rm new family of projective bundles}, we require that a family of $\PGL(n,\CC)$-Higgs bundles $(\PP(\Ee),\varPhi)$ comes from a family of Higgs bundles $(\Ee,\varPhi)$. Define S-equivalence of families of semista\-ble $\PGL(n,\CC)$-Higgs bundles pointwise and consider the moduli functor that associates to every scheme $Y$ the set of S-equivalence classes of families of semistable $\PGL(n,\CC)$-Higgs bundles of topological type $\wt{d}$ pa\-ra\-me\-tri\-zed by $Y$. There exists a coarse moduli space $\MmM(\PGL(n,\CC))_{\wt{d}}$ associated to this functor and, if $d$ is  a representative of $\wt{d}$, it can be proved that $\MmM(\PGL(n,\CC))_{\wt{d}}$ is the quotient of $(\tr)^{-1}(0) \subset \MmM(\GL(n,\CC))_d$ by $\Pic(X)^{0}$ and therefore
\beq
\label{eq MmPGLnC is MmMGLnC quotiented by Pic}
\MmM(\PGL(n,\CC))_{\wt{d}} \cong \quotient{(\det,\tr)^{-1}(\Oo(x_0)^{\otimes d},0)}{\Pic^0(X)[n]}.
\eeq

\subsection{Symplectic and orthogonal Higgs bundles}
\label{sc sp2mC o OnC-HB}

A {\it $\Sp(2m,\CC)$-Higgs bundle} over the elliptic curve $X$ is a triple $(E,\Omega,\Phi)$, where $(E,\Omega)$ is a $\Sp(2m,\CC)$-bundle and $\Phi \in H^0(X,\End E)$ is an endomorphism of $E$ which anticommutes with $\Omega$, i.e. for every $x \in X$ and every $u,v \in E_x$, 
\[
\Omega(u,\Phi(v)) = -\Omega(\Phi(u),v).
\]
Two $\Sp(2m,\CC)$-Higgs bundles, $(E,\Omega,\Phi)$ and $(E',\Omega',\Phi')$ are {\it isomorphic} if there exists an isomorphism of $\Sp(2m,\CC)$-bundles $f : (E',\Omega') \to (E,\Omega)$ such that $\Phi' = f^{-1} \Phi f$.

An {\it $\Ort(n,\CC)$-Higgs bundle} (resp. {\it a $\SO(n,\CC)$-Higgs bundle}) over $X$ is a triple $(E,Q,\Phi)$ (resp. a quadruple $(E,Q,\Phi,\tau)$), where $(E,Q)$ is an $\Ort(n,\CC)$-bundle (resp. $(E,Q,\tau)$ is a $\SO(n,\CC)$-bundle) and $\Phi \in H^0(X,\End E)$ is an endomorphism of $E$ which anticommutes with $Q$, i.e. for every $x \in X$ and every $u,v \in E_x$, 
\[
Q(u,\Phi(v)) = -Q(\Phi(u),v).
\]
Two $\Ort(n,\CC)$-Higgs bundles (resp. $\SO(n,\CC)$-Higgs bundles) $(E,Q,\Phi)$ and $(E',Q',\Phi')$ (resp. $(E,Q,\Phi,\tau)$ and $(E',Q',\Phi',\tau')$) are {\it isomorphic} if there exists an isomorphism of $\Ort(n,\CC)$-bundles $f : (E',Q') \to (E,Q)$ (resp. an isomorphism of $\SO(n,\CC)$-bundles $f : (E',Q',\tau') \to (E,Q,\tau)$) such that $\Phi' = f^{-1} \Phi f$.

The following notions of stability, semistability and polystability of $\Sp(2m,\CC)$, $\Ort(n,\CC)$ and $\SO(n,\CC)$-Higgs bundles are the notions of stability worked out in \cite{oscar&ignasi&gothen} (see also \cite{marta&oscar} for the stability of $\SO(n,\CC)$-Higgs bundles).


Let $(E,\Theta,\Phi)$ be a $\Sp(2m,\CC)$-Higgs bundle or an $\Ort(n,\CC)$-Higgs bundle. Note that $\mu(E) = 0$ since $E \cong E^*$. We say that $(E, \Theta, \Phi)$ is {\it semistable} if and only if, for any $\Phi$-invariant isotropic subbundle $F$ of $E$,
\[
\mu(F)  \leq \mu(E) = 0,
\]
and it is {\it stable} if the above inequality is strict for any proper $\Phi$-invariant isotropic subbundle. If $(E,\Theta,\Phi)$ is a semistable $\Sp(2m,\CC)$-Higgs bundle or $\Ort(n,\CC)$-Higgs bundle we have a {\it Jordan-H\"older filtration} (see for instance \cite{oscar&ignasi&gothen}) of $\Phi$-invariant subbundles, consisting of
\[
0 = E_{0} \subsetneq E_{1} \subsetneq \dots \subsetneq E_{k-1} \subsetneq  E_k \subseteq E_k^{\perp_\Theta} \subsetneq E_{k-1}^{\perp_\Theta} \subsetneq \dots \subsetneq  E_1^{\perp_\Theta} \subsetneq E_0^{\perp_\Theta} = E,
\] 
such that, if we denote by $\ol{\Phi}_i$ and $\ol{\Phi}'_i$ the Higgs fields on $E_i/E_{i-1}$ and $E^{\perp_\Theta}_{i-1}/E^{\perp_\Theta}_i$ induced by $\Phi$, we have for every $i \leq k$ that $(E_i/E_{i-1},\ol{\Phi}_i)$ and $(E^{\perp_\Theta}_{i-1}/E^{\perp_\Theta}_i, \ol{\Phi}'_i)$ are stable Higgs bundles and $\Theta$ induces an isomorphism $\theta_i : (E_i/E_{i-1},\ol{\Phi}_i) \stackrel{\cong}{\lra} ((E^{\perp_\Theta}_{i-1}/E^{\perp_\Theta}_i)^*,-(\ol{\Phi}'_i)^t)$. If $E_k^{\perp_\Theta}/E_k$ is non-zero and $\wt{\Phi}$ is the Higgs field on it induced by $\Phi$, $\Theta$ induces a non-degenerate symplectic or symmetric form $\wt{\Theta}$ anticommuting with $\wt{\Phi}$, in such a way that $(E_k^{\perp_\Theta}/E_k,\wt{\Theta},\wt{\Phi})$ is a stable $\Sp(2m',\CC)$ or $\Ort(n',\CC)$-Higgs bundle.

We define the {\it asso\-cia\-ted graded object} of a semistable $\Sp(2m,\CC)$ or $\Ort(n,\CC)$-Higgs bundle $(E,\Theta,\Phi)$ in terms of a Jordan-Hölder filtration, 
\begin{align*}
\gr(E,\Theta,\Phi) := & (E_k^{\perp_\Theta}/E_k,\wt{\Theta},\wt{\Phi}) \quad \oplus 
\\
&  \oplus \bigoplus_{i = 1}^{k} \left( (E_i/E_{i-1}) \oplus (E^{\perp_\Theta}_{i-1}/E^{\perp_\Theta}_{i}) , \bpm 0 & b\theta_i^t \\   \theta_i & 0 \epm, \bpm \ol{\Phi}_i & \\ & \ol{\Phi}'_i \epm \right),
\end{align*}
where $b=-1$ for $\Sp(2m,\CC)$-Higgs bundles and $b = 1$ for $\Ort(n,\CC)$-Higgs bundles. As in the previous cases, the Jordan-H\"older filtration may not be unique, but $\gr(E,\Theta,\Phi)$ is unique up to isomorphism.

For $n>2$, a $\SO(n,\CC)$-Higgs bundle is {\it semistable, stable or polystable} if it is semistable, stable or polystable as an $\Ort(n,\CC)$-Higgs bundle. The graded object of a $\SO(n,\CC)$-Higgs bundle is given by the graded object of the underlying $\Ort(n,\CC)$-Higgs bundle.

We define stability and S-equivalence of families pointwise, as we did for families of Higgs bundles. S-equivalence between families of stable objects implies isomorphism pointwise. The moduli functors are defined by associating to every scheme $Y$ the set of S-equivalence classes of families parametrized by $Y$. In \cite{simpson2} it is proved that there exist moduli spaces associated to these moduli functors; we denote them respectively by $\MmM(\Sp(2m,\CC))$, $\MmM(\Ort(n,\CC))_{k,a}$ and $\MmM(\SO(n,\CC))_{\omega_2}$. We denote the stable loci by $\MmM^{st}(\Sp(2m,\CC))$, $\MmM^{st}(\Ort(n,\CC))$ and $\MmM^{st}(\SO(n,\CC))$.

\subsection{Normality of the moduli spaces}
\label{sc normality}

Normality is an important local property of some algebraic varieties whose study is simplified by the following result (see for example \cite[Chap II.8]{hartshorne}).

\bpr {\bf (Serre's Criterion of normality)}
\label{tm Serre's Criterion}
The algebraic variety $Y$ is normal if and only if the following properties are satisfied: 
\begin{enumerate}

\item[(R1)] \label{num R1} the codimension of the singular locus $\Sing(Y)$ is strictly greater than $1$;

\item[(S2)] \label{num S2} for every $y \in Y$ we have 
\[
\depth(\Oo_{Y,y}) \geq \min \{ 2, \dim(\Oo_{Y,y}) \}.
\]
\end{enumerate} 
\epr   

Recall that $\Gamma_\RR = \RR \times_\ZZ \Gamma$, where $\Gamma$ denotes the universal central extension by $\ZZ$ of the fundamental group $\pi_1(X)$. To study the normality of the moduli space $\MmM(G)_d$ of $G$-Higgs bundles it is enough, by the Isosingularity Theorem of \cite{simpson2}, to study the normality of the moduli space $\Rr(G)_d$ of central representations of $\Gamma_\RR$ into $G$. Recall that the latter is a GIT quotient of the space $\Hom^c(\Gamma_\RR,G)_d$ of central representations of $\Gamma_\RR$ in $G$ with topological invariant $d \in \pi_1(G)$ by the adjoint action of $G$. 

\bpr {\bf (\cite[Section 0.2]{mumford&fogarty&kirwan})}
\label{pr normality is preserved under GIT}
Suppose that $Y$ is a categorical quotient of $Z$ by the complex reductive Lie group $G$. If $Z$ is normal and locally integral, then $Y$ is normal and locally integral.
\epr     

We focus our study on the normality of $\Hom^c(\Gamma_\RR,G)_d$. If $\rho$ is a central representation of topological type $0$, one has $\rho(\RR) = 1$. Therefore the space $\Hom^c(\Gamma_\RR,G)_0$ is identified with the space $\Hom(\pi_1(X), G)_0$ of representations of the fundamental group of the curve with topologically trivial invariant. Since $\pi_1(X) \cong \ZZ \times \ZZ$, we have that $\Hom(\pi_1(X), G)$ is the {\it commuting variety of the group} $G$,
\[
C(G) := \{ (g,g') \in G \times G \qua | \qua [g,g'] = \id \}.
\]
Also the {\it commuting variety of a Lie algebra} $\gG$, defined by
\[
C(\gG) := \{ (A,A') \in \gG \times \gG \qua | \qua [A,A'] = 0 \},
\]    
will be important for us. Commuting varieties have been extensively studied. For example, one can check that they are reduced algebraic varieties and, by the following result of \cite{richardson}, $C(\gG)$ is irreducible and therefore integral. When $G$ is semisimple and simply connected, so is $C(G)$. 

\bpr {\bf (\cite[Corollary 2.5 and Theorem C]{richardson})}
\label{pr C G irreducible}
Let $G$ be a semisimple simply connected complex Lie group and let $\gG$ be a reductive Lie algebra. Then $C(G)$ and $C(\gG)$ are irreducible algebraic varieties.
\epr 

The property of $C(G)$ and $C(\gG)$ being normal has not been determined in general but there is a long-standing conjecture stating that the commuting variety $C(\gG)$ is always normal (see \cite{popov} and \cite{premet}). The following result states that $C(\gG)$ satisfies Serre's condition (R1).

\bpr {\bf (\cite[Corollary 1.9]{popov})}
\label{pr the singular locus of C_n has codimension greater than 2}
Let $\gG$ be a noncommutative complex reductive Lie algebra. The singular locus $\Sing(C(\gG))$ has codimension greater than or equal to $2$.
\epr 

If $C(\gG)$ were Cohen-Macaulay then it would automatically satisfy Serre's condition (S2) and would therefore be a normal variety. Up to now, this has only been determined for $\gG\lL(n,\CC)$ and lower rank using computations performed by the computer program Macaulay.

\bpr {\bf (\cite{bayer&stillman&stillman} for $n = 3$ and \cite{hreinsdottir} for $n = 4$)}
\label{pr Cohen-Macaulay in lower rank}
Let $n \leq 4$. Then the commuting variety $C(\gG\lL(n,\CC))$ is Cohen-Macaulay.
\epr 

Note that $C(\GL(n,\CC))$ is an open subset of $C(\gG\lL(n,\CC))$ given by the non-vanishing of the determinant. 

\bco
\label{co C GLnC is normal for n lower than 5}
Let $n \leq 4$. Then the commuting varieties $C(\GL(n,\CC))$ and $C(\gG\lL(n,\CC))$ are normal.
\eco 




\btm
\label{tm normality of MmMGLnC_0 when n lower than 5}
Let $n \leq 4$. Then the moduli spaces $\MmM(\GL(n,\CC))_0$, $\MmM(\SL(n,\CC))$ and $\MmM(\PGL(n,\CC))_0$ are normal.
\etm 

\proof
By Proposition \ref{pr C G irreducible} and Corollary \ref{co C GLnC is normal for n lower than 5}, the variety $\Hom(\pi_1(X),\GL(n,\CC))$ is normal and integral. By Proposition \ref{pr normality is preserved under GIT} the moduli space of representations $\Rr(\GL(n,\CC))_0$ is normal. This implies that $\MmM(\GL(n,\CC))_0$ is normal by the Isosingularity Theorem \cite[Theorem 10.6]{simpson2}.   

For $\MmM(\SL(n,\CC))$, in view of (\ref{eq MmSLnC is the preimage of det tr}) and the fact that 
$\Pic^0(X) \times H^0(X,\Oo)$ is smooth, it follows that $\MmM(\SL(n,\CC))$ is a complete intersection in some open set in $\MmM(\GL(n,\CC))_0$, so (S2) holds. It is easy to check that (R1) also holds (or see Theorem \ref{tm MmMSLnC}).

Finally, for $\MmM(\PSL(n,\CC))_0$, use the result for $\MmM(\SL(n,\CC))$ and \eqref{eq MmPGLnC is MmMGLnC quotiented by Pic}.
\qed

\section{Description of the moduli spaces}
\label{sc description of Mm}

\subsection{Stability in terms of the underlying bundle}
\label{sc stability of higgs bundles}

The triviality of the canonical line bundle simplifies the study of the semistability of Higgs bundles over elliptic curves.

\bpr
\label{pr E Phi semistable iff E semistable} 
The Higgs bundle $(E,\Phi)$ is semistable if and only if $E$ is semistable.
\epr

\proof 
If the vector bundle $E$ is semistable, every subbundle $F$ satisfies $\mu(F) \leq \mu(E)$, so the Higgs bundle $(E,\Phi)$ is semistable too. 

Suppose $E$ is not semistable and take its Harder-Nara\-sim\-han filtration
\[
0 = F_0 \subset F_1 \subset \dots \subset F_{s-1} \subset F_s = E,
\]
where the $E_i= F_i/F_{i-1}$ are semistable with $\mu(E_i) > \mu(E_j)$ if $i < j$. In particular we have $H^0(X,\Hom(E_i,E_j)) = 0$ if $i < j$.

The subbundle $F_1 = E_1$ has $\mu(F_1) > \mu(E)$. Suppose $\Phi(F_1)$ is non-zero and take $F_\ell$ such that $\Phi(F_1) \subseteq F_\ell$ but $\Phi(F_1)\nsubseteq F_{\ell-1}$. Thus $\Phi$ induces a non-zero morphism from $E_1 = F_1$ to $E_\ell = F_\ell / F_{\ell -1}$, but there are no non-zero morphisms unless $\ell = 1$. Then $F_1$ is $\Phi$-invariant and $(E,\Phi)$ is not semistable.  
\qed

\bco
\label{co if gdc=1 E Phi stable then E stable}
If $\gcd(n,d)=1$, then $(E,\Phi)$ is stable if and only if $E$ is stable.
\eco

We need to extend Corollary \ref{co if gdc=1 E Phi stable then E stable} to the non-coprime case.

\bpr
\label{pr E Phi stable iff E stable} 
$(E,\Phi)$ is stable if and only if $E$ is stable.
\epr 

\proof
Take first $E$ strictly semistable and indecomposable, so $E\cong E'\otimes F_h$, where $h=\gcd(n,d)>1$ and $E'$ is stable of rank $n'=\frac{n}h$ and degree $d'=\frac{d}h$. The endomorphism bundle satisfies
\begin{align*}
H^0(X,\End E) \cong & \quad H^0(X, \End E'\otimes \End F_h) 
\\
\cong & \bigoplus_{L_i \in \Pic^0(X)[n']} H^0(X,L_i \otimes \End F_h).
\end{align*}
We have $\End F_h \cong F_1 \oplus \dots F_{2h-1}$, and $H^0(X,L_i\otimes F_j) = 0$ for every $F_j$ and every non-trivial $L_i \in \Pic^0(X)[n]$, so $H^0(X,L_i \otimes \End F_h) = 0$ if $L_i \ncong \Oo$. This implies that $H^0(X,\End E) \cong H^0(X,\End F_h)$ and every endomorphism of $E$ has the form $\Phi=\id_{E'}\otimes\phi$ for some $\phi\in \End F_h$. 

The vector bundle $F_h$ has a unique subbundle isomorphic to $\Oo$. Let $\phi \in \End F_h$. If $\phi$ is non-zero, $\phi(\Oo)$ is either zero or a subbundle of $F_h$ with a non-zero section, so $\phi(\Oo) \subseteq \Oo$. If $\phi = 0$, $\Oo$ is again $\phi$-invariant. The subbundle $E'\otimes\Oo$ contradicts the stability of $(E,\Phi)$.

Now consider $E$ strictly semistable and decomposable. One can write 
\[
E \cong \bigoplus^s_{j=1} E_j \cong \bigoplus^s_{j=1} (E'_j \otimes (F_{h_{j,1}}\oplus \dots \oplus F_{h_{j,k_j}}))
\]
where the $E'_j$ are stable vector bundles of rank $n'$ and degree $d'$ such that $E_j'\not\cong E_k'$ if $j \neq k$. For $j \neq k$,
\[
\Hom(E_j'\otimes F_{h_j},E_k'\otimes F_{h_k})=0,
\]
so $\Phi(E_j) \subset E_j$ and then the Higgs bundle $(E,\Phi)$ decomposes in a direct sum of $(E_j,\Phi_j)$ where $\Phi_j$ is the restriction to $E_j$. This contradicts stability of $(E, \Phi)$ when $s \neq 1$.

Now consider $E \cong E' \otimes (F_{h_1} \oplus \dots \oplus F_{h_k})$ with $E'$ stable, so 
\[
\End E\cong \End \left(\bigoplus_{j=1}^k F_{h_j} \right),\ \ \Phi=\id_{E'} \otimes\phi,\ \ \phi\in H^0\left(X,\End \left(\bigoplus_{j=1}^k F_{h_j}\right)\right).
\]
Now $\bigoplus_{j=1}^k F_{h_j}$ has a unique subbundle $\Oo^k$ and $H^0(\Oo^k)\stackrel{\cong}{\hookrightarrow}H^0(\bigoplus_{j=1}^k F_{h_j})$. So $\Oo^k$ is $\phi$-invariant and $E'\otimes \Oo^k$ is $\Phi$-invariant. This implies that $(E,\Phi)$ is strictly semistable unless all $h_j=1$ and $E\cong E'\otimes \Oo^h$. In this case $\phi\in \End \Oo^k=\{k\times k\mbox{ matrices}\}$. Choose an eigenvector $v$ for $\phi$. Then $E'\otimes v$ contradicts stability of $(E,\Phi)$ and $(E,\Phi)$ is strictly semistable.
\qed

\bco
\label{co E Phi polystable only if E polystable}
Let $(E,\Phi)$ be polystable of rank $n$ and degree $d$, and $h = \gcd(n,d)$. Then 
\[
(E,\Phi) = \bigoplus_{i = 1}^h (E_i,\Phi_i)
\]
where $E_i$ is a stable bundle of rank $n'= \frac{n}{h}$ and degree $d' = \frac{d}{h}$ and
\[
\Phi_i \in H^0(X,\End E_i) \cong H^0(X,\Oo) \otimes \id_{E_i}.
\]
Furthermore $(E,\Phi)$ is polystable only if $E$ is polystable.
\eco

\brm
\nr{Although the polystability of the underlying vector bundle is a necessary condition for the polystability of the Higgs bundle, it is not sufficent.} 

\nr{To illustrate this fact, take $(E,\Phi)$ such that $E \cong \Oo \oplus \Oo$ and $\Phi = A \otimes 1_X$, where $A$ is non-diagonalizable. It follows that $(E,\Phi)$ is an indecomposable Higgs bundle. Since $E$ is a strictly polystable vector bundle, $(E,\Phi)$ is not stable. Also, $(E,\Phi)$ is indecomposable so it is not possible to express $(E,\Phi)$ as a direct sum of stable Higgs bundles.} 
\erm

The following result is well known for $\Sp(2m,\CC)$, $\Ort(n,\CC)$ or $\SO(n,\CC)$ over smooth projective curves of arbitrary genus. We give here a proof specific for $g = 1$.

\bpr
\label{pr E Theta Phi semistable implies E Phi semistable} 
If $(E,\Theta,\Phi)$ is a semistable $\Sp(2m,\CC)$ or $\Ort(n,\CC)$-Higgs bundle, then $(E,\Phi)$ is semistable. If $n>2$ and $(E,Q,\Phi,\tau)$ is a semistable $\SO(n,\CC)$-Higgs bundle, then $(E,\Phi)$ is semistable. 
\epr

\proof
Suppose that $(E,\Phi)$ is not semistable and take the first term of its Harder-Narasim\-han filtration
\[
0 = F_0 \subset F_1 \subset \dots \subset F_n = E.
\]
The subbundle $F_1$ is $\Phi$-invariant, $\mu(F_1) > \mu(E) = 0$ and $(F_1,\Phi_{F_1})$ is semistable; by Proposition \ref{pr E Phi semistable iff E semistable} so is $F_1$. By \cite[Lemma 10.1]{atiyah&bott} $F_1^* \otimes F_1^*$ is semistable of negative degree. If the bundle $F_1$ is not isotropic then $\Theta_{F_1,F_1} \in H^0(X,F_1^* \otimes F_1^*)$ must be non-zero and there exists a line subbundle of $F_1^* \otimes F_1^*$ of degree $\ge0$, contradicting the semistability of $F_1^* \otimes F_1^*$. So $F_1$ is isotropic and contradicts the semistability of $(E,\Theta,\Phi)$.

The statement for $\SO(n,\CC)$-Higgs bundles follows from the fact that (when $n>2$) a $\SO(n,\CC)$-Higgs bundle $(E,Q,\Phi,\tau)$ is semistable if and only if the underlying $\Ort(n,\CC)$-Higgs bundle $(E,Q,\Phi)$ is semistable.
\qed

Combining Proposition \ref{pr E Theta Phi semistable implies E Phi semistable} with Proposition \ref{pr E Phi semistable iff E semistable} gives us

\bco
\label{co E Theta semistable if E Theta Phi semistable}
If $(E,\Theta,\Phi)$ is a semistable $\Ort(n,\CC)$ or $\Sp(2m,\CC)$-Higgs bundle, then $(E,\Theta)$ is semistable. If $(E,Q,\Phi,\tau)$ is a semistable $\SO(n,\CC)$-Higgs bundle, then $(E,Q,\tau)$ is a semistable $\SO(n,\CC)$-bundle. 
\eco

\bpr
\label{pr E Theta Phi stable if E Phi is stable}
Let $(E,\Theta,\Phi)$ be a stable $\Ort(n,\CC)$ or $\Sp(2m,\CC)$-Higgs bundle. Then $\Phi = 0$ and $(E,\Theta)$ is a stable $\Ort(n,\CC)$ or $\Sp(2m,\CC)$-bundle.

Let $(E,Q,\Phi,\tau)$ be a stable $\SO(n,\CC)$-Higgs bundle. Then $\Phi = 0$ and $(E,Q,\tau)$ is a stable $\SO(n,\CC)$-bundle.
\epr 

\proof
By Proposition \ref{pr E Theta Phi semistable implies E Phi semistable} $(E,\Phi)$ is semistable. If $(E,\Phi)$ is stable, then $E$ is stable by Proposition \ref{pr E Phi stable iff E stable}; since $\deg E = 0$, it follows that $E$ has rank $1$. Since $\Phi$ anticommutes with $\Theta$, this implies that $\Phi = 0$.

Assume now that $(E,\Phi)$ is not stable and take $(F_1,\Phi_1)$ to be the first term of a Jordan-H\"older filtration; since $\Phi$ anticommutes with $\Theta$, one infers that $F_1^{\perp_\Theta}$ is $\Phi$-invariant. Following \cite{ramanan}, we define the $\Phi$-invariant subbundles $W$ and $V$ generated by $F_1 \cap F_1^{\perp_\Theta}$ and $F_1 + F_1^{\perp_\Theta}$ respectively. Note that $V = W^{\perp_\Theta}$ and then the exact sequence
\[
0 \lra W \lra F_1 \oplus F_1^{\perp_{\Theta}} \lra V \lra 0
\]
implies that $\deg(W) + \deg(W^{\perp_\Theta}) = \deg(F_1) + \deg(F_1^{\perp_\Theta})$. Recalling (\ref{eq deg F^perp = deg F}), one has that $\deg(W) = \deg(F_1) = 0$. This implies that $W = 0$ by the stability of $(E,\Theta,\Phi)$ and the fact that $W$ is isotropic and $\Phi$-invariant. Then $V = E = F_1 \oplus F_1^{\perp_{\Theta}}$ and
\[
(E,\Phi) = (F_1,\Phi_1) \oplus (F_1^{\perp_\Theta}, \Phi_2),
\]
where $\Phi_2$ is the restriction of $\Phi$ to $F_1^{\perp_\Theta}$. Note that $(F_1,\Phi_1)$ is a stable Higgs bundle and $(F_1^{\perp_\Theta},\Phi_2)$ is semistable. By induction, we decompose $(F_1^{\perp_\Theta},\Phi_2)$ into stable factors pro\-ving that $(E,\Phi)$ is polystable. Since $\deg(E)=0$, by Corollary \ref{co E Phi polystable only if E polystable} $(E,\Phi)$ decomposes as a direct sum of Higgs line bundles of degree $0$,  
\[
(E,\Phi) \cong \bigoplus (L_i,\phi_i).
\]
No $L_i$ is isotropic, otherwise it would contradict the stability of $(E,\Theta,\Phi)$ since $L_i$ is $\Phi$-invariant. Thus $\Theta$ restricts to $L_i$ and since $\Theta$ anticommutes with the Higgs field we have $\phi_i = -\phi_i$. Therefore $\Phi = 0$. 
\qed

\bco
The only stable $\Ort(k,\CC)$-Higgs bundles are $(E^{st}_{k,a},Q^{st}_{k},0)$, where the underlying $\Ort(k,\CC)$-bundles $(E^{st}_{k,a},Q^{st}_{k})$ are the stable $\Ort(k,\CC)$-bundles appearing in Theorem \ref{tm MSp2mC and MstOnC}.
\eco

\bco
Let $(E,\Theta,\Phi)$ be a polystable $\Ort(n,\CC)$ or $\Sp(2m,\CC)$-Higgs bundle. Then $(E,\Theta)$ is a polystable $\Ort(n,\CC)$ or $\Sp(2m,\CC)$-bundle.

Let $(E,Q,\Phi,\tau)$ be a polystable $\SO(n,\CC)$-Higgs bundle. Then $(E,Q,\tau)$ is a polystable $\SO(n,\CC)$-bundle.
\eco 

\proof
This follows from Proposition \ref{pr E Theta Phi stable if E Phi is stable} and the definition of polystability.
\qed

\bco
\label{pr decomposition of polystable Sp2mC OnC or SOnC-Higgs bundles}
A $\Sp(2m,\CC)$-Higgs bundle over an elliptic curve is polystable if and only if it is isomorphic to a direct sum of polystable $\Sp(2,\CC)$-Higgs bundles. 

An $\Ort(n,\CC)$-Higgs bundle is polystable if and only if it is isomorphic to a direct sum of polystable $\Ort(2,\CC)$-Higgs bundles or it is isomorphic to a direct sum of polystable $\Ort(2,\CC)$-Higgs bundles and a stable $\Ort(m,\CC)$-Higgs bundle (where $m = 1,3$ or $4$).

Let $n>2$. Let $(E,Q,\Phi,\tau)$ be a polystable $\SO(n,\CC)$-Higgs bundle with Stiefel-Whitney class $w_2$,
\begin{enumerate}
\item \label{num SO2m-HB with w_2 = 0} if $n=2n'$ and $w_2 = 0$, it is a direct sum of $m_{n,w_{2}} = n'$ stable $\SO(2,\CC)$-Higgs bundles of trivial degree,
\item \label{num SO2m+1-HB with w_2 = 0} if $n=2n'+1$ and $w_2 = 0$, it is a direct sum of $(E^{st}_{1,0},Q^{st}_1,0, 1)$ and $m_{n,w_{2}} = n'$ stable $\SO(2,\CC)$-Higgs bundles of trivial degree,
\item \label{num SO2m+1-HB with w_2 = 1} if $n=2n'+1$ and $w_2 = 1$, it is a direct sum of $(E^{st}_{3,0}, Q^{st}_3, 0, 1)$ and $m_{n,w_{2}} = n'-1$ stable $\SO(2,\CC)$-Higgs bundles of trivial degree,
\item \label{num SO2m-HB with w_2 = 1}if $n=2n'$ and $w_2 = 1$, it is a direct sum of $(E^{st}_{4,0},Q^{st}_4,0,1)$ and $m_{n,w_{2}} = n'-2$ stable $\SO(2,\CC)$-Higgs bundles of trivial degree.
\end{enumerate}
\eco

\proof
This follows from Propositions \ref{pr E Theta Phi stable if E Phi is stable} and \ref{pr decomposition of polystable Sp2mC OnC or SOnC-bundles}. 
\qed

\bpr
Let $G$ be $\GL(n,\CC)$, $\SL(n,\CC)$, $\PGL(n,\CC)$, $\Sp(2m,\CC)$, $\Ort(n,\CC)$ or $\SO(n,\CC)$. Then we have a morphism between moduli spaces
\[
a_G : \MmM(G)_d \lra M(G)_d
\]
defined by the underlying principal bundle of a Higgs bundle for these structure groups. The fibres of this map are connected.
\epr 

\proof
The existence of this map follows from Proposition \ref{pr E Phi semistable iff E semistable} and Corollary \ref{co E Theta semistable if E Theta Phi semistable}.

For any Higgs bundle $(E,\Phi)$ in the fibre of $E$ there exists a path connecting $(E,\Phi)$ to $(E,0)$. So the fibre is connected. 
\qed

\bco
\label{co conn comp of MmMG are classified by the conn comp of MG}
For any topological type $d$, $\MmM(G)_d$ is connected.
\eco 

\proof
This follows since $M(G)_d$ is connected by (\ref{eq definition of varsigma in the noncoprime case}), (\ref{eq MSLnC}), (\ref{eq MPGLnC}) and Theorem \ref{tm MG}.
\qed

\subsection{The moduli space of Higgs bundles}
\label{sc MmMGLnC}

The abelian group structure defined on $X$ induces naturally an abelian group structure on $T^*X$. Recall here that the canonical bundle is trivial, so
\[
T^*X \cong X \times \CC. 
\]
The moduli space $\MmM(\GL(1,\CC))_0$ is naturally identified with $T^*\Pic^0(X)$, so the isomorphism $\varsigma^{x_0}_{1,d} : X \to \Pic^d(X)$ of (\ref{eq definition of varsigma_x0 1 0}) gives an isomorphism 
\[
\eta^{x_0}_{1,d} : T^*X \stackrel{\cong}{\lra} \MmM(\GL(1,\CC))_d 
\]
for every $d$.

\btm
\label{tm no stable for gcd greter than 1}
Let $n$ and $d$ be two integers and write $h = \gcd(n,d)$. If $h > 1$, we have
\[
\MmM^{st}(\GL(n,\CC))_d = \emptyset.
\]
\etm

\proof
This follows from Proposition \ref{pr E Phi stable iff E stable} and Atiyah’s results.
\qed

\bpr
\label{pr for gcd = 1 det tr is an isomorphism}
If $\gcd(n,d) = 1$, then $\MmM^{st}(\GL(n,\CC))_d = \MmM(\GL(n,\CC))_d$ and the morphism
\[
(\det,\tr) : \MmM(\GL(n,\CC))_d \stackrel{\cong}{\lra} \MmM(\GL(1,\CC))_d
\]
is an isomorphism.
\epr

\proof
If $\gcd(n,d) = 1$, then there are no strictly semistable vector bundles nor strictly semistable Higgs bundles, so $M^{st}(\GL(n,\CC))_d = M(\GL(n,\CC))_d$ and $\MmM^{st}(\GL(n,\CC))_d = \MmM(\GL(n,\CC))_d$. Moreover, the determinant gives an isomorphism
\[
\det : M(\GL(n,\CC))_d \stackrel{\cong}{\lra} M(\GL(1,\CC))_d.
\]
Taking differentials we obtain the isomorphism
\[
(\det,\tr) : T^*M(\GL(n,\CC))_d \stackrel{\cong}{\lra} T^*M(\GL(1,\CC))_d.
\]
Now, for any $(n,d)$, $T^*M^{st}(\GL(n,\CC))_d$ is an open subscheme of $\MmM^{st}(\GL(n,\CC))_d$. Due to Proposition \ref{pr E Phi stable iff E stable}, every stable Higgs bundle has a stable underlying vector bundle, so $T^*M^{st}(\GL(n,\CC))_d = \MmM^{st}(\GL(n,\CC))_d$. Since $T^*M(\GL(n,\CC))_d = \MmM(\GL(n,\CC))_d$, the result follows.
\qed

If $\gcd(n,d) = 1$, we set $\eta^{x_0}_{n,d} = (\det, \tr)^{-1} \circ \qua \eta^{x_0}_{1,d}$.

\btm
\label{tm moduli space of HB if gdc=1}
If $\gcd(n,d)=1$, then
\beq
\label{eq definition of xi in the coprime case}
\eta^{x_0}_{n,d} : T^*X \stackrel{\cong}{\lra} \MmM(\GL(n,\CC))_d
\eeq
is an isomorphism.
\etm

\proof 
This follows from Proposition \ref{pr for gcd = 1 det tr is an isomorphism} and the definition of the isomorphism $\eta^{x_0}_{1,d}$.
\qed

\bpr
\label{pr Ee_nd with h=1 is universal}
Let $\gcd(n,d) = 1$. There exists a universal family $\Ee^{x_0}_{n,d} = (\Vv^{x_0}_{n,d},\varPhi_{n,d})$ of stable Higgs bundles of rank $n$ and degree $d$ parametrized by $T^*X$.
\epr

\proof
Consider the family of Higgs bundles over $X \times \CC \cong T^*X$ such that the restriction to $X \times {(x,t)}$ is given by the pair $((\Vv^{x_0}_{n,d})_x, \lambda \otimes \id_{(\Vv^{x_0}_{n,d})_x})$ where $\lambda = \frac{1}{n} t$. We can check that for any $(x,t) \in T^*X$ one has
\[
\eta^{x_0}_{n,d}((x,t)) = [(\Ee^{x_0}_{n,d})_{(x,t)}]_S.
\]
Any family $\Ff \to X \times Y$ induces canonically a morphism $\nu_{\Ff} : Y \to \MmM^{st}(\GL(n,\CC))_d$. Clearly, the composition $f = (\eta^{x_0}_{n,d})^{-1} \circ \nu_{\Ff}$ is a morphism $f : Y \to T^*X$ such that $\Ff$ is S-equivalent to $f^* \Ee^{x_0}_{n,d}$.
\qed

If $\gcd(n,d) = h > 1$, $n'=\frac{n}{h}$ and $d'=\frac{d}{h}$, we take the fibre product of $h$ copies of $\Ee^{x_0}_{n',d'}$ to define the family $\Ee^{x_0}_{n,d}$ of polystable Higgs bundles parametrized by $T^*Z_h = T^*X \times \stackrel{h}{\dots} \times T^*X$.

\brm
\label{rm Ee_z1 S-equivalent to Ee_z2 if z1 is a permutation of z2}
\nr{The action of $\sym_h$ on $Z_h$ induces an action of $\sym_h$ on $T^*Z_h$. If $w_1$ and $w_2$ are two points of $T^*Z_h$, the Higgs bundles $(\Ee^{x_0}_{n,d})_{w_1}$ and $(\Ee^{x_0}_{n,d})_{w_2}$ are S-equivalent if and only if for some $\gamma \in \sym_h$ one has $w_2 = \gamma \cdot w_1$ (i.e. $w_1$ is a permutation of $w_2$).} 
\erm

\btm
\label{tm MmGLnC is the normalization of MmMGLnC}
Let $h = \gcd(n,d)$. Consider the moduli space $\MmM(\GL(n,\CC))_d$ associated to the usual moduli problem.
\begin{enumerate}[(i)]
\item There exists a bijective morphism 
\[
\eta^{x_0}_{n,d} : \Sym^h T^*X \lra \MmM(\GL(n,\CC))_d.
\]
\item $\Sym^h T^*X$ is the normalization of $\MmM(\GL(n,\CC)_d$.
\item $\MmM(\GL(n,\CC))_d$ is irreducible.
\item If $n\le4$, then $\eta^{x_0}_{n,0}$ is an isomorphism.
\end{enumerate}
\etm

\proof
{\it (i)} The family $\Ee^{x_0}_{n,d}$ induces a morphism
\[
\nu_{\Ee^{x_0}_{n,d}} : T^*Z_h = T^*X \times \stackrel{h}{\dots} \times T^*X  \lra \MmM(\GL(n,\CC))_d,
\]
and by Remark \ref{rm Ee_z1 S-equivalent to Ee_z2 if z1 is a permutation of z2} it factors through $\Sym^h T^*X$ giving the required bijective morphism.

Since $\Sym^h T^*X$ is normal, {\it (ii)} follows from Zariski's Main Theorem. The continuity of $\eta^{x_0}_{n,d}$ and the irreducibility of $\Sym^h T^*X$ imply {\it (iii)}. {\it (iv)} follows from Theorem \ref{tm normality of MmMGLnC_0 when n lower than 5}.
\qed

Although we cannot prove normality except when $n\le4$, $d=0$, or $\gcd(n,d)=1$, we can identify the singular set of $\MmM(\GL(n,\CC)_d$ in all cases.

\btm
\label{tm singular set of MmMGLnC} The singular set $\Sing(\MmM(\GL(n,\CC)_d)$ coincides with the set $S$ of points represented by polystable Higgs bundles for which at least two of the direct summands are isomorphic. In particular, if $h=\gcd(n,d)\ge2$, this set has codimension $2$.
\etm

\proof
Note first that $S$ is the image under $\nu_{\Ee^{x_0}_{n,d}}$ of the union $\Delta$ of the diagonals in $T^*Z_h = T^*X \times \stackrel{h}{\dots} \times T^*X$. This union coincides with the union of the fixed point sets of the elements of $\sym _h$ acting on $T^*Z_h$. Since  $\codim\Delta=2$, the image of $\Delta$ in $T^*Z_h/\sym_h$ is $\Sing(T^*Z_h/\sym_h)$. Since $S=\eta^{x_0}_{n,d}(\Sing(T^*Z_h/\sym_h))$ and $\eta^{x_0}_{n,d}$ is bijective, it follows that $S$ is contained in $\Sing(\MmM(\GL(n,\CC)_d))$.

To see that the points outside $S$ are smooth, we consider the deformation complex of a Higgs bundle $(E,\Phi)$. This gives rise to an exact sequence
\[
H^0(X,\End E)\stackrel{e_0(\Phi)}{\lra}H^0(X,\End E)\lra T\lra H^1(X,\End E)\stackrel{e_1(\Phi)}{\lra}H^1(X,\End E),
\]
where $T$ is the infinitesimal deformation space of $(E,\Phi)$, the maps $e_i(\Phi)$ are given by $e_i(\Phi)(\Psi)=[\Psi,\Phi]$ and $e_1(\Phi)$ is the Serre dual of $e_0(\Phi)$. (For this, see \cite{nitsure} and note that in our case the canonical bundle is trivial.) If $(E,\Phi)$ is polystable with non-isomorphic summands $(E_i,\Phi_i)$ ($1\le i\le h$), we have $\Phi_i=\lambda_i\id_{E_i}$; moreover, if $E_i\cong E_j$ with $i\ne j$, then $\lambda_i\ne\lambda_j$. It follows that
\[
\Psi \in\ker e_0(\Phi) \Longleftrightarrow \Psi = \bigoplus \mu_i \id_{E_i}
\]
for some $\mu_i\in\CC$. Hence $\codim(\im e_0(\Phi)) = h$. By duality $\dim(\ker e_1(\Phi))=h$, so 
\[
\dim(T)=2h.
\]
On the other hand, we have a family of polystable Higgs bundles parametrized by $T^*Z_h$, which has dimension $2h$. It follows that the local deformation space of $(E,\Phi)$ is a complete family and is smooth. So the point of the moduli space represented by $(E,\Phi)$ is smooth.
\qed

\brm
\nr{Note that this gives a direct proof that $\Sing(\MmM(\GL(n,\CC)_d))$ has codimension $2$ without using Proposition \ref{pr the singular locus of C_n has codimension greater than 2}. Moreover, we have shown that a small deformation of a polystable Higgs bundle with non-isomorphic direct summands is polystable. In fact, one can show directly, without using deformation theory, that, if $(E,\Phi)$ is a semistable Higgs bundle such that $\gr(E,\Phi)$ has non-isomorphic direct summands, then $(E,\Phi)$ is polystable. To see this, it is sufficient to show that, if $E_1$ and $E_2$ are stable bundles of the same slope and $(E_1,\Phi_1)\not\cong(E_2,\Phi_2)$, then any extension}
\[
0\lra(E_1,\Phi_1)\lra(E,\Phi)\lra(E_2,\Phi_2)\lra0
\]
\nr{splits. If $E_1\not\cong E_2$, this is clear since then $H^1(E_2^*\otimes E_1)=0$ and $H^0(E_i^*\otimes E_j)=0$ when $i\ne j$. If $E_1=E_2$ but $E\not\cong E_1\oplus E_2$, then $E\cong E_1\otimes F_2$ and any endomorphism of $E$ has the form $\id_{E_1}\otimes\phi$ for some endomorphism $\phi$ of $F_2$. Now $\phi=\lambda\id_{F_2}+\nu$ where $\nu^2=0$, which implies that $\Phi_1=\Phi_2=\lambda\id_{E_1}$, contradicting the hypothesis that $(E_1,\Phi_1)\not\cong(E_2,\Phi_2)$. Finally, if $E\cong E_1\oplus E_1$ but $\Phi_1\ne\Phi_2$, then $\Phi$ can be written as a diagonalisable $2\times 2$ matrix, giving the required splitting of $(E,\Phi)$.}
\erm

\subsection{A new moduli problem} 
\label{sc new moduli problem}

In view of the fact that the normality of the  moduli space is in general an open question, we shall now define a moduli functor whose associated moduli space is the normalization of $\MmM(\GL(n,\CC))_d$. To do so, we modify the definition of family. 

We say that a family of semistable Higgs bundles $\Ff \to X \times Y$ is {\it locally graded} if for every point $y$ of $Y$ there exists an open subset $U \subset Y$ containing $y$ and a set of families $\Ee_1, \dots, \Ee_\ell$ where each $\Ee_i$ is a family of stable Higgs bundles parametrized by $U$ such that for every point $y'$ of $U$ we have 
\[
\gr \Ff |_{X \times \{ y' \}} \cong \bigoplus_{i = 1}^{\ell} \Ee_i |_{X \times \{ y' \}}, 
\]
and therefore $\Ff|_{X \times U} \sim_S \bigoplus_i \Ee_i$. The new moduli functor associates to a scheme $Y$ the set of S-equivalence classes of locally graded families of semistable Higgs bundles of rank $n$ and degree $d$ parametrized by $Y$.

A family $\Ee$ parametrized by $Z$ is said to have the {\it local universal property} if, for any family $\Ff$ parametrized by $Y$ and any point $y \in Y$, there exists a neighbourhood $U$ containing $y$ and a (not necessarily unique) morphism $f : U \to Z$ such that $\Ff|_{U} \sim f^*\Ee$. Families with the local universal property are very useful for describing moduli spaces as we can see in the following result.

\bpr {\bf (Proposition 2.13 of \cite{newstead})}
\label{pr Mm described by a family with the local universal property}
Let us suppose that there exists a family $\Ee$ parametrized by $Z$ with the local universal property. Suppose that there exists a group $\Gamma$ acting on $Z$ such that $\Ee_{z_1} \sim \Ee_{z_2}$ if and only if $z_1$ and $z_2$ lie in the same orbit of this action. Then a categorical quotient of $Z$ by $\Gamma$ is a coarse moduli space if and only if it is an orbit space.
\epr 

The definition of locally graded families is justified thanks to the next result.

\bpr
\label{pr Ee_nd has the local universal property}
The family $\Ee^{x_0}_{n,d}$ has the local universal property among locally graded families of semistable Higgs bundles.
\epr

\proof
Take $\Ff$ to be any locally graded family of semistable Higgs bundles of rank $n$ and degree $d$ parametrized by $Y$. Set $h = \gcd(n,d)$, $n' = n/h$ and $d' = d/h$. Since $\Ff$ is locally graded, for every $y \in Y$ there exists an open neighborhood $U$ and a set of families $\Ee_1, \dots \Ee_h$ of stable Higgs bundles of rank $n'$ and degree $d'$ parametrized by $U$ and such that
\[
\Ff |_{X \times U} \sim_S \bigoplus_{i = 1}^{h} \Ee_i.
\]
Since $\Ee^{x_0}_{n',d'}$ is a universal family, for every $\Ee_i$ there exists $f_i : U \to T^*X$ such that $\Ee_i \sim_S f_i^* \Ee^{x_0}_{n',d'}$. Setting $f = (f_1, \dots, f_h)$, we have
\[
\Ff|_{X \times U} \sim_S f^* \Ee^{x_0}_{n,d}.
\]
\qed

\btm
\label{tm Moduli de GL-Higgs bundles}
There exists a coarse moduli space $\Mm(\GL(n,\CC))_d$ of S-equivalence classes of semistable Higgs bundles of rank $n$ and degree $d$. Furthermore, there is a morphism $T^*X \to \MmM(\GL(n,\CC))_d$ which induces an isomorphism
\[
\xi^{x_0}_{n, d} :  \Sym^h T^*X \stackrel{\cong}{\lra} \Mm(\GL(n,\CC))_d.
\]
Moreover $\Mm(\GL(n,\CC))_d$ is the normalization of $\MmM(\GL(n,\CC))_d$.
\etm

\proof
Since $\Sym^h T^*X = T^*Z_h/\sym_h$ is an orbit space, by Propositions \ref{pr Mm described by a family with the local universal property} and \ref{pr Ee_nd has the local universal property} and Remark \ref{rm Ee_z1 S-equivalent to Ee_z2 if z1 is a permutation of z2}, the coarse moduli space $\Mm(\GL(n,\CC))_d$ exists and is isomorphic to $\Sym^h T^*X$. The isomorphism is given by
\[
\begin{array}{cccc}
\xi^{x_0}_{n,d} : & \Sym^h T^*X  & \lra &  \Mm(\GL(n,\CC))_d 
\\ 
& [(x_1,t_1), \dots, (x_h,t_h)]_{\sym_h} &\longmapsto & \left[ \Ee^{x_0}_{n',d'}|_{X \times (x_1,t_1)} \oplus \dots \oplus \Ee^{x_0}_{n',d'}|_{X \times (x_h,t_h)} \right]_S.
\end{array}
\]

The last statement follows from Theorem \ref{tm MmGLnC is the normalization of MmMGLnC}.
\qed

\subsection{Moduli spaces of special linear and projective Higgs bundles}
\label{sc MmSLnC and MmPGLnC}

\btm
\[
\MmM^{st}(\SL(1,\CC)) = \{ pt \},
\]
and for $n > 1$
\[
\MmM^{st}(\SL(n,\CC)) = \emptyset. 
\]
\etm

\proof
When the determinant is trivial the vector bundle has degree $0$. The only stable Higgs bundles with degree $0$ are the Higgs line bundles and therefore $(\Oo,0)$ is the only stable $\SL(n,\CC)$-Higgs bundle.
\qed

\btm
\label{tm Mm^st PGLnC}
\[
\MmM^{st}(\PGL(n,\CC))_{\wt{d}} = \emptyset
\]
unless $n \in \ZZ^+$ and $\wt{d} \in \ZZ_n$ are coprime. In that case
\[
\MmM(\PGL(n,\CC))_{\wt{d}} \qua  = \MmM^{st}(\PGL(n,\CC))_{\wt{d}} \qua \cong \{ pt \}.
\]
\etm

\proof
We recall from Theorems \ref{tm no stable for gcd greter than 1} and \ref{tm moduli space of HB if gdc=1} that there exist stable Higgs bundles of rank $n$ and degree $d$ if and only if $\gcd(n,d)=1$. Recall that a stable $\PGL(n,\CC)$-Higgs bundle of degree $\wt{d}$ can be represented by a pair $(\PP(E), \Phi)$, where $E$ is a stable vector bundle of degree $d$ and $\tr \Phi = 0$. When $\gcd(n,d) =  1$, $E$ is determined up to tensoring by a line bundle $L$, so $\PP(E)$ is uniquely determined. Since $E$ is stable, the endomorphism $\Phi$ is a scalar multiple of the identity; so $\Phi = 0$. Thus $(\PP(E),\Phi)$ is uniquely determined.
 
\qed

We now study the preimage of $(\det,\tr)$. We have a morphism
\[
\begin{array}{cccc}
\suma^h_{T^*X} : & \Sym^h T^*X  & \lra &  T^*X 
\\ 
& [(x_1,t_1), \dots, (x_h,t_h)]_{\sym_h} &\longmapsto & \sum_{i=1}^h (x_i,t_i).
\end{array}
\]
Take the bijective morphism given in Theorem \ref{tm MmGLnC is the normalization of MmMGLnC} to construct the diagram
\beq
\label{eq suma_T*X commutes with det tr for MmM}
\xymatrix{
\Sym^h T^*X \ar[d]_{\suma^h_{T^*X}}  \ar[rr]^{\eta^{x_0}_{n,d}}_{1:1} & & \MmM(\GL(n,\CC))_d \ar[d]^{(\det,\tr)}
\\
T^*X  \ar[rr]^{\cong \qquad}_{\eta^{x_0}_{1,d}\qquad } & & \Pic(X)^d \times H^0(X,\Oo).
}
\eeq
One can easily check that this diagram commutes.

\btm
\label{tm MmMSLnC}
Consider the moduli space $\MmM(\SL(n,\CC))$ associated to the usual moduli problem.
\begin{enumerate}[(i)]
\item There exists a bijective morphism
\[
\hat{\eta}^{x_0}_n : \quotient{(T^*X \times \stackrel{n-1}{\dots} \times T^*X)}{\sym_n}\lra\MmM(\SL(n,\CC)),
\]
where the action of $\sym_n$ on $T^*Z_{n-1} = T^*X \times \stackrel{n-1}{\dots} \times T^*X$ is induced by the action on $Z_{n-1} = X \times \stackrel{n-1}{\dots} \times X$ used in (\ref{eq MSLnC}).
\item $\MmM(\SL(n,\CC))$ is irreducible.
\item $\Sing(\MmM(\SL(n,\CC)))$ coincides with the set of points represented by polystable Higgs bundles for which at least two of the direct summands are isomorphic.
\item  If $n\ge2$, $\Sing(\MmM(\SL(n,\CC)))$ has codimension $2$.
\item If $n\le4$, the morphism of (i) is an isomorphism.
\end{enumerate}
\etm

\proof 
{\it (i)} By \eqref{eq MmSLnC is the preimage of det tr} and \eqref{eq suma_T*X commutes with det tr for MmM}, we have a bijective morphism
\[
\hat{\eta}_n : (\suma^n_{T^*X})^{-1}((x_0,0))\lra\MmM(\SL(n,\CC)).
\]
The action of $\sym_n$ on $T^*Z_{n-1}$ induced by the action of $\sym_n$ on $Z_{n-1}$ used in (\ref{eq MSLnC}) gives us an isomorphism between $T^*Z_{n-1}/\sym_n$ and $(\suma^n_{T^*X})^{-1}((x_0,0))$.

{\it (ii)} follows from the continuity of $\hat{\eta}^{x_0}_n$ and the irreducibility of the source.

{\it (iii)} and {\it (iv)} are proved in the same way as Theorem \ref{tm singular set of MmMGLnC}. The only changes are to replace $\End E$ by $\End_0 E$ and note that $h=n$. This gives $\dim(T)=2(n-1)=\dim(T^*Z_{n-1})$.

{\it (v)} follows by Zariski's Main Theorem from Theorem \ref{tm normality of MmMGLnC_0 when n lower than 5}.
\qed

Let $h = \gcd(n,d)$, $n' = \frac{n}{h}$ and $d' = \frac{d}{h}$. Recall the action of $X$ on $\Sym^{h}X$ defined  in (\ref{eq definition of the action of X on Sym^h X}); we extend it to an action of $X$ on $\Sym^h T^*X$. The commutativity of (\ref{eq commutivity of the action of X on Sym^h and the action of Pic on MGLnC}) implies the commutativity of 
\beq
\label{eq f_n'h commutes with tensorization for MmM}
\xymatrix{
X \times \Sym^h T^*X \ar[d]_{\varsigma^{x_0}_{1,0} \times \eta^{x_0}_{n,d}}^{1:1} \ar[rr] & &   \Sym^h T^*X \ar[d]^{\eta^{x_0}_{n,d}}_{\cong}
\\
\Pic^0(X) \times \MmM(\GL(n,\CC))_d \ar[rr]^{\qquad  - \otimes - } &  & \MmM(\GL(n,\CC))_d,
}
\eeq
where the right arrow is the bijective morphism of Theorem \ref{tm MmGLnC is the normalization of MmMGLnC}.
Recall that the action of $X[n]$ with weight $n'$ corresponds to the action of $X[h]$ with weight $1$ and that the action of $X[h]$ commutes with the action of the symmetric group $\sym_h$. 

Note that $h = \gcd(n,d)$ is equal to $n$ when $d = 0$.

\btm 
\label{tm MmMPGLnC}
Consider the moduli space $\MmM(\PGL(n,\CC))_{\wt{d}}$ associated to the usual mo\-du\-li problem.
\begin{enumerate}[(i)]
\item There exists a bijective morphism
\[
\check{\eta}^{x_0}_{n,\wt{d}} : \quotient{(T^*X \times \stackrel{h-1}{\dots} \times T^*X)}{\sym_h \times X[n]}\lra\MmM(\PGL(n,\CC)_{\tilde d}),
\]
where the action of $\sym_h\times X[h]$ on $T^*Z_{h-1} = T^*X \times \stackrel{h-1}{\dots} \times T^*X$ is induced by the action on $Z_{h-1} = X \times \stackrel{h-1}{\dots} \times X$ used in (\ref{eq MPGLnC}).
\item $\MmM(\PGL(n,\CC)_{\tilde d}$ is irreducible.
\item $\Sing(\MmM(\PGL(n,\CC))_{\tilde{d}})$ coincides with the set of points represented by polystable Higgs bundles for which at least two of the direct summands are isomorphic.
\item  If $h \geq 2$, $\Sing(\MmM(\PGL(n,\CC))_{\wt{d}})$ has codimension $2$.
\item If $n \leq 4$ and $\tilde{d}=0$, the morphism of (i) is an isomorphism.
\end{enumerate}
\etm

\proof 
The proof is similar to that of Theorem \ref{tm MmMSLnC}.

{\it (i)} By \eqref{eq MmPGLnC is MmMGLnC quotiented by Pic}, \eqref{eq suma_T*X commutes with det tr for MmM} and \eqref{eq f_n'h commutes with tensorization for MmM}, we have a bijective morphism
\[
((\suma^h_{T^*X})^{-1}((x_0,0)))/X[h]\lra\MmM(\PGL(n,\CC))_{\tilde{d}}.
\]
The action of $\sym_h \times X[h]$ on $T^*Z_{h-1}$ induced by the action of $\sym_h$ on $Z_{h-1}$ used in (\ref{eq MPGLnC}) gives an isomorphism between $T^*Z_{h-1}/(\sym_h \times X[h])$ and $((\suma^h_{T^*X})^{-1}((x_0,0)))/X[h]$.

{\it (ii)} follows from the continuity of $\check{\eta}^{x_0}_{n,\wt{d}}$ and the irreducibility of the source.

To prove {\it (iii)} we again modify the proof of Theorem \ref{tm singular set of MmMGLnC} by replacing $\End E$ by $\End_0 E$, giving $\dim(T)  =2(h-1)=\dim(\MmM(\PGL(n,\CC))_{\tilde{d}})$.

{\it (iv)} follows from Theorem \ref{tm MmMSLnC} when we note that, by Atiyah's results, $X[h]$ acts freely on $((\suma^h_{T^*X})^{-1}((x_0,0)))$.

{\it (v)} follows by Zariski's Main Theorem from Theorem \ref{tm normality of MmMGLnC_0 when n lower than 5}.
\qed

In order to obtain normal moduli spaces in the cases not covered by Theorems \ref{tm MmMSLnC} (iv) and \ref{tm MmMPGLnC} (iv), we again need to restrict ourselves to locally graded families.

A {\it locally graded family of semistable $\SL(n,\CC)$-Higgs bundles} is a family of $\SL(n,\CC)$-Higgs bundles which is locally graded as a family of semistable Higgs bundles. The existence of a moduli space $\Mm(\SL(n,\CC))$ for the corresponding moduli functor follows from the existence of $\Mm(\GL(n,\CC))_0$ and 
\[
\Mm(\SL(n,\CC)) \cong (\det, \tr)^{-1} ((\Oo,0)) \subset \Mm(\GL(n,\CC))_0.
\]

A {\it locally graded family of semistable $\PGL(n,\CC)$-Higgs bundles} is a fa\-mi\-ly $(\PP(\Vv),\varPhi)$ of $\PGL(n,\CC)$-Higgs bundles which is the projectivization of a locally graded family of semistable Higgs bundles $(\Vv,\varPhi)$. 
It follows that the existence of $\Mm(\GL(n,\CC))_d$ implies the existence of $\Mm(\PGL(n,\CC))_{\wt{d}}$ and furthermore
\[
\Mm(\PGL(n,\CC))_{\wt{d}} \cong \quotient{(\det,\tr)^{-1}((\Oo(x_0)^{\otimes d},0))}{\Pic^0(X)[n]}.
\]

As in the case of (\ref{eq suma_T*X commutes with det tr for MmM}) and (\ref{eq f_n'h commutes with tensorization for MmM}), one can easily check that the diagrams 
\beq
\label{eq suma_T*X commutes with det tr for Mm}
\xymatrix{
\Sym^h T^*X \ar[d]_{\suma^h_{T^*X}}  \ar[rr]^{\xi^{x_0}_{n,d}}_{\cong} & & \Mm(\GL(n,\CC))_d \ar[d]^{(\det,\tr)}
\\
T^*X  \ar[rr]^{\cong \qquad}_{\xi^{x_0}_{1,d}\qquad } & & \Pic(X)^d \times H^0(X,\Oo)
}
\eeq
and 
\beq
\label{eq f_n'h commutes with tensorization for Mm}
\xymatrix{
X \times \Sym^h T^*X \ar[d]_{\varsigma^{x_0}_{1,0} \times \xi^{x_0}_{n,d}}^{\cong} \ar[rr] & &   \Sym^h T^*X \ar[d]^{\xi^{x_0}_{n,d}}_{\cong}
\\
\Pic^0(X) \times \Mm(\GL(n,\CC))_d \ar[rr]^{\qquad  - \otimes - } &  & \Mm(\GL(n,\CC))_d
}
\eeq
commute.

Note that in the commuting diagrams (\ref{eq suma_T*X commutes with det tr for Mm}) and (\ref{eq f_n'h commutes with tensorization for Mm}) we have an isomorphism $\xi^{x_0}_{n,d}$ for any value of $n$ and $d$, while in the case of (\ref{eq suma_T*X commutes with det tr for MmM}) and (\ref{eq f_n'h commutes with tensorization for MmM}) the bijection $\eta^{x_0}_{n,d}$ is proved to be an isomorphism only in lower rank.

\btm
\label{tm MmSLnC and PGLnC}
Consider the action of $\sym_n$ on $T^*Z_{n-1} = T^*X \times \stackrel{n-1}{\dots} \times T^*X$ induced by the action on $Z_{n-1} = X \times \stackrel{n-1}{\dots} \times X$ used in (\ref{eq MSLnC}). We have isomorphisms
\beq
\label{eq MmSLnC}
\xymatrix{
\hat{\xi}^{x_0}_{n} : \quotient{(T^*X \times \stackrel{n-1}{\dots} \times T^*X)}{\sym_n} \ar[r]^{\qquad \quad \cong} & \Mm(\SL(n,\CC))
}
\eeq
and 
\beq
\label{eq MmPGLnC}
\xymatrix{
\check{\xi}^{x_0}_{n,\wt{d}} : \quotient{(T^*X \times \stackrel{h-1}{\dots} \times T^*X)}{\sym_h \times X[h]} \ar[r]^{\qquad \qquad \cong} & \Mm(\PGL(n,\CC))_{\wt{d}}.
}
\eeq 

Moreover, there are natural bijective morphisms
\[
\Mm(\SL(n,\CC)) \lra \MmM(\SL(n,\CC))
\]
and 
\[
\Mm(\PGL(n,\CC))_{\wt{d}} \lra \MmM(\PGL(n,\CC))_{\wt{d}}.
\]
Hence $\Mm(\SL(n,\CC))$ is the normalization of $\MmM(\SL(n,\CC))$ and $\Mm(\PGL(n,\CC))_{\wt{d}}$ the normalization of $\MmM(\PGL(n,\CC))_{\wt{d}}$.
\etm

\proof 
The proofs of (\ref{eq MmSLnC}) and (\ref{eq MmPGLnC}) follow the proofs of Theorems \ref{tm MmMSLnC} and \ref{tm MmMPGLnC}.

The remaining statements follow from the universal property of $\MmM(\SL(n,\CC))$ and $\MmM(\PGL(n,\CC))_{\wt{d}}$, the normality of the quotients of $T^*X \times  \stackrel{n-1}{\dots} \times T^*X$ by the finite groups $\sym_n$ and $\sym_n \times X[n]$ and Zariski's Main Theorem.
\qed

\brm
\label{rm MmPGLnC is also T*Z/sym_h}
\nr{By Lemma \ref{lm the quotient of the weighted action is again X times dots times X} we have that $Z_{h-1}/X[h] \cong Z_{h-1}$. Since $T^*X \cong X\times \CC$ we have that $T^*Z_{h-1}/X[h] \cong T^*Z_{h-1}$. Thus}
\[
\Mm(\PGL(n,\CC))_{\wt{d}} \cong \quotient{T^*Z_{h-1}}{\sym_h} = \quotient{T^*X \times \stackrel{h-1}{\dots} \times T^*X}{\sym_h},
\]
\nr{where the action of $\sym_h$ of $T^*Z_{h-1}$ is induced by the action of $\sym_h$ on $Z_{h-1}$ used in Remark \ref{rm MPGLnC is also Z/sym_h}  and is different from the action used in (\ref{eq MmSLnC})}. 
\erm

\subsection{Moduli spaces of symplectic and orthogonal Higgs bundles}
\label{sc MmSp2mC}

\btm
\begin{align*}
& \MmM^{st}(\Sp(2m,\CC)) \cong  M^{st}(\Sp(2m,\CC)) = \emptyset,
\\
& \MmM^{st}(\Ort(n,\CC)) \cong M^{st}(\Ort(n,\CC)),
\\
& \MmM^{st}(\SO(n,\CC)) \cong M^{st}(\SO(n,\CC)).
\end{align*}
\etm

\proof
This follows from Proposition \ref{pr E Theta Phi stable if E Phi is stable}.
\qed

Recall the universal family of Higgs line bundles $\Ee^{x_0}_{1,0} = (\Vv^{x_0}_{1,0},\varPhi_{1,0})$ with degree $0$ and the family of $\Sp(2,\CC)$-bundles $\wt{\Vv}^{x_0}_2 = (\Vv^{x_0}_{1,0} \oplus \Vv^*_{1,0},\varOmega)$. We can check that $\varOmega$ anticommutes with $\varPhi_{1,0} \oplus (-\varPhi_{1,0})$ and then
\[
\wt{\Ee}_2 = (\Vv^{x_0}_{1,0} \oplus \Vv^*_{1,0},\varOmega,\varPhi_{1,0} \oplus (-\varPhi_{1,0}))
\]
is a family of polystable $\Sp(2,\CC)$-Higgs bundles parametrized by $T^*X$. Similarly, we can construct a family of polystable $\Ort(2,\CC)$-Higgs bundles parametrized by $T^*X$,
\[
\mathring{\Ee}_2 = (\Vv^{x_0}_{1,0} \oplus (\Vv^{x_0}_{1,0})^*,\Qq,\varPhi_{1,0} \oplus (-\varPhi_{1,0})),
\]
and a family of stable $\SO(2,\CC)$-Higgs bundles of degree $0$,
\[
\overline{\Ee}^{x_0}_{2} = \left( \Vv^{x_0}_{1,0} \oplus (\Vv^{x_0}_{1,0})^*,\Qq,\varPhi_{1,0} \oplus (-\varPhi_{1,0}), \sqrt{-1} \right).
\]

\brm
\label{rm Ee2_1 cong Ee2_2 if they are related by Gamma}
\nr{The restrictions of $\wt{\Ee}^{x_0}_2$ and $\mathring{\Ee}^{x_0}_2$ to two different points of $T^*X$, $z_1 = (x_1,t_1)$ and $z_2 = (x_2,t_2)$, give S-equivalent (isomorphic) Higgs bundles $(\wt{\Ee}_2)_{z_1} \sim_S (\wt{\Ee}_2)_{z_2}$ and $(\mathring{\Ee}^{x_0}_2)_{z_1} \sim_S (\mathring{\Ee}^{x_0}_2)_{z_2}$ if and only if $z_1 = - z_2$. Two points $z_1 \neq z_2$ of $T^*X$ give non-S-equivalent $\SO(2,\CC)$-Higgs bundles, $(\ol{\Ee}^{x_0}_2)_{x_1} \not\sim_S (\ol{\Ee}^{x_0}_2)_{x_1}$}.
\erm

Now we define $\wt{\Ee}^{x_0}_{2m}$ to be the family of polystable $\Sp(2m,\CC)$-Higgs bundles induced by taking the fibre product of $m$ copies of $\wt{\Ee}_2$. Note that $\wt{\Ee}^{x_0}_{2m}$ is parametrized by $T^*Z_m = T^*X \times \stackrel{m}{\dots} \times T^*X$. 

We define $\mathring{\Ee}^{x_0}_{n,k,a}$ to be the family of polystable $\Ort(n,\CC)$-Higgs bundles given by the direct sum of $(E^{st}_{k,a},Q^{st}_{k},0)$ and $m_k = \frac{n-k}{2}$ copies of $\mathring{\Ee}_2$. This family is parametrized by $T^*Z_{m_k} = T^*X \times \stackrel{m_k}{\dots} \times T^*X$.

Analogously, for all values of $(n,w_2)$, we define families of polystable $\SO(n,\CC)$-Higgs bundles $\overline{\Ee}^{x_0}_{n,w_2}$ as follows:
\bit
\item $\ol{\Ee}^{x_0}_{2m,0}$ is given by $m$ copies of $\ol{\Ee}_2$ and therefore it is parametrized by $T^*Z_m = T^*X \times \stackrel{m}{\dots} \times T^*X$,
\item $\ol{\Ee}^{x_0}_{2m+1,0}$ is given by the direct sum of $(E^{st}_{1,0},Q^{st}_1,0,1)$ and $m$ copies of $\ol{\Ee}_2$ and therefore it is parametrized by $T^*Z_m = T^*X \times \stackrel{m}{\dots} \times T^*X$,
\item $\ol{\Ee}^{x_0}_{2m+1,1}$ is given by the direct sum of $(E^{st}_{3,0},Q^{st}_3,0,1)$ and $m-1$ copies of $\ol{\Ee}_2$ and therefore it is parametrized by $T^*Z_{m-1} = T^*X \times \stackrel{m-1}{\dots} \times T^*X$,
\item $\ol{\Ee}^{x_0}_{2m,1}$ is given by the direct sum of $(E^{st}_{4,0},Q^{st}_4,0,1)$ and $m-2$ copies of $\ol{\Ee}_2$ and therefore it is parametrized by $T^*Z_{m-2} = T^*X \times \stackrel{m-2}{\dots} \times T^*X$.
\eit

\brm
\label{rm wtEe_z1 S-equivalent to wtEe_z2 if z1 and z2 are related by Gamma}
\nr{Consider the action of $\Gamma_{m'}$ on $T^*Z_{m'}$ induced by the action of $\Gamma_{m'}$ on $Z_{m'}$ used in Remark \ref{rm wtVv_z1 S-equivalent to wtVv_z2 if z1 and z2 are related by Gamma} (i.e. the action induced by $\ZZ_2$ acting on $T^*X$ as $-1 \cdot (x,t) = (-x,-t)$). Take the families $\wt{\Ee}^{x_0}_{2m}$, $\mathring{\Ee}^{x_0}_{n,k,a}$ and $\ol{\Ee}^{x_0}_{n,w_2}$ when $(n,w_2)$ is $(2m,1)$, $(2m+1,0)$ or $(2m+1,1)$. These families are parametrized by $T^*Z_{m'} = T^*X \times \stackrel{m'}{\dots} \times T^*X$ for some $m'$. By Remarks \ref{rm Ee2_1 cong Ee2_2 if they are related by Gamma} and \ref{rm wtVv_z1 S-equivalent to wtVv_z2 if z1 and z2 are related by Gamma}, two points $z_1, z_2 \in T^*Z_{m'}$ parametrize isomorphic (S-equivalent) Higgs bundles on these families if and only if $z_1 = \gamma \cdot z_2$ for some $\gamma \in \Gamma_m$. By Remarks \ref{rm Ee2_1 cong Ee2_2 if they are related by Gamma} and \ref{rm ovlineVv_z1 sim ovlineVv_z2 if they are related by the action of Delta},  $(\overline{\Ee}^{x_0}_{2m,0})_{z_{1}}$ is isomorphic (S-equivalent) to $(\overline{\Ee}^{x_0}_{2m,0})_{z_{2}}$ if and only if there exists $\gamma \in \Delta_m$ such that $\gamma \cdot z_1 = z_2$.}
\erm

Recall Corollary \ref{co conn comp of MmMG are classified by the conn comp of MG}.

\btm
\label{tm MmG is the normalization of MmMG}
Let $G$ be $\Sp(2m,\CC)$, $\Ort(n,\CC)$ or $\SO(n,\CC)$. Consider the moduli spaces $\MmM(G)_d$ associated to the usual moduli problems.

\begin{enumerate}[(i)]
\item We have bijective morphisms
\begin{align*}
\wt{\eta}^{x_0}_m : \Sym^m (T^*X/\ZZ_2) \lra & \MmM(\Sp(2m,\CC)),
\\
\mathring{\eta}^{x_0}_{n,k,a} : \Sym^{(n-k)/2} (T^*X/\ZZ_2) \lra & \MmM(\Ort(n,\CC))_{k,a},
\\
\ol{\eta}^{x_0}_{2m,0} : \quotient{T^*X \times \stackrel{m}{\dots} \times T^*X}{\Delta_m} \lra & \MmM(\SO(2m,\CC))_{0},
\\
\ol{\eta}^{x_0}_{n,w_2} : \Sym^{m'} (T^*X/\ZZ_2) \lra & \MmM(\SO(n,\CC))_{w_2}, 
\end{align*}
where $m' = m$ if $(n,w_2) = (2m+1,0)$, $m' = m-1$ if $(n,w_2) = (2m+1,1)$ and $m' = m-2$ if $(n,w_2) = (2m,1)$. 

\item $\MmM(G)_d$ is irreducible.

\item The normalization of $\MmM(G)_d$ is isomorphic to the source of the bijection.

\item $\Sing(\MmM(G)_d)$ coincides with the set of points represented by polystable $G$-Higgs bundles for which at least two of the direct summands of the underlying polystable Higgs bundle are isomorphic.

\item If $\dim(\MmM(G)_d) \geq 2$, $\Sing(\MmM(G)_d)$ has codimension $2$.
\end{enumerate}
\etm

\proof
{\it (i)} The corresponding families $\wt{\Ee}^{x_0}_{2m}$, $\mathring{\Ee}^{x_0}_{n,k,a}$ or $\ol{\Ee}^{x_0}_{n,w_2}$ induce morphisms 
\[
\nu_{G,d} : T^*Z_{G,d} = T^*X \times \stackrel{m'}{\dots} \times T^*X  \lra \MmM(G)_d,
\]
where $m'$ depends on $G$ and $d$. By Remark \ref{rm wtEe_z1 S-equivalent to wtEe_z2 if z1 and z2 are related by Gamma} when $\Gamma_{G,d}$ is the corresponding finite group, this map factors through $T^*Z_{G,d}/\Gamma_{G,d}$ giving the required bijective morphism.

{\it (ii)} follows from the continuity of the morphisms and the irreducibility of the sources.

{\it (iii)} follows from Zariski's Main Theorem.

For the proof of {\it (iv)} and {\it (v)} take a $G$-Higgs bundle $(E,\Theta,\Phi)$ (resp. $(E,\Theta,\Phi, \tau)$) and denote by $P$ the principal $G$-bundle associated to $(E,\Theta)$ (resp. $(E,\Theta,\tau)$) and by $\varphi$ the section of the adjoint bundle $P(\gG)$ associated to $\Phi$. If $T$ is the infinitesimal deformation space of $(E,\Theta,\Phi)$ (resp. $(E,\Theta,\Phi,\tau)$), by \cite{biswas&ramanan} one has the exact sequence
\[
H^0(X,P(\gG))\stackrel{e_0(\varphi)}{\lra} H^0(X,P(\gG)) \lra T\lra H^1(X,P(\gG)) \stackrel{e_1(\varphi)}{\lra} H^1(X,P(\gG)),
\]
where $e_i(\varphi)(\psi) = [\psi, \varphi]$ and $e_1(\varphi)$ is the Serre dual of $e_0(\varphi)$. Recall that the canonical bundle is trivial in our case. The standard representation of $G$ gives us an isomorphism between $H^i(X,P(\gG))$ and $H^i(X,\End E)^{\Theta}$, where the latter is given by the elements $\Psi$ of $H^i(X,\End E)$ such that $\Theta(u,\Psi(v)) = -\Theta(\Psi(u),v)$ for $u,v$ in $E_x$ and every point $x \in X$. The maps $e_i(\varphi)$ correspond to 
\[
\begin{array}{cccc}
e_i(\Phi) : & H^i(X,\End E)^{\Theta} & \lra & H^i(X,\End E)^{\Theta} 
\\ 
& \Psi &\longmapsto & [\Psi,\Phi].
\end{array}
\]
From this point the proof is similar to that of Theorem \ref{tm singular set of MmMGLnC}, replacing $H^0(X,\End E)$ by $H^0(X,\End E)^{\Theta}$ and taking into account in each case the dimension of $\dim(\MmM(G)_d)$. 
\qed

\brm
\nr{The bijection $\eta^{x_0}_{G,d}$ is an isomorphism if and only if $\MmM(G)_d$ is normal, but the normality of the moduli space of Higgs bundles over a smooth projective curve of genus $g=1$ is an open question.} 
\erm

As in Section \ref{sc new moduli problem}, we modify the moduli problem in such a way that the associated moduli space will be isomorphic to the normalization of $\MmM(G)_d$. 

We say that a family of semistable $\Sp(2m,\CC)$-Higgs bundles  $\Ee \to X \times Y$ is {\it locally graded} if for every $y \in Y$ there exists an open subset $U \subset Y$ containing $y$ and a set of families of Higgs bundles $(\Vv_1, \varPhi_1), \dots (\Vv_m,\varPhi_m)$ of rank $1$ and degree $0$ parametrized by $U$, such that 
\[
\Ee|_{X \times U} \sim_S \bigoplus_{i=1}^m \left( \Vv_i \oplus \Vv_i^*, \bpm  0 & - 1 \\ 1 & 0  \epm , \bpm \varPhi_{i} & \\ & -\varPhi_{i}  \epm \right).
\]

We define {\it locally graded families} of semistable $\Ort(n,\CC)$-Higgs bundles in similar terms; the only change is that there might exist a stable $\Ort(k,\CC)$-Higgs bundle $(E^{st}_{k,a},Q^{st}_{k},0)$ such that
\[
\Ee|_{X \times U} \sim_S (E^{st}_{k,a},Q^{st}_{k},0) \oplus \bigoplus_{i=1}^{\frac{n-k}{2}} \left( \Vv_i \oplus \Vv_i^*, \bpm  0 & 1 \\ 1 & 0  \epm , \bpm \varPhi_{i} & \\ & -\varPhi_{i}  \epm \right).
\]

Analogously, we say that a family of semistable $\SO(n,\CC)$-Higgs bundles $\Ee \to X \times Y$ is {\it locally graded} if for every $y \in Y$ there exists an open subset $U \subset Y$ containing $y$, families $(\Vv_i, \varPhi_i)$ parametrized by $U$ and (if $(n,w_2) \neq (2m,0)$) a stable $\SO(k,\CC)$-Higgs bundle $(E^{st}_{k,0},Q^{st}_{k},0, 1)$ such that
\[
\Ee|_{X \times U} \sim_S (E^{st}_{k,0},Q^{st}_{k},0, 1) \oplus \bigoplus_{i=1}^{\frac{n-k}{2}} \left( \Vv_i \oplus \Vv_i^*, \bpm  0 & 1 \\ 1 & 0  \epm , \bpm \varPhi_{i} & \\ & -\varPhi_{i}  \epm, \sqrt{-1} \right).
\]

For $G = \Sp(2m,\CC)$, $\Ort(n,\CC)$ and $\SO(n,\CC)$, we define a new moduli functor associating to any scheme $Y$ the set of S-equivalence classes of locally graded families of semistable $G$-Higgs bundles parametrized by $Y$. 

\bpr
\label{pr wtEe_2m has the local universal property}
The families $\wt{\Ee}^{x_0}_{2m}$, $\mathring{\Ee}^{x_0}_{n,k,a}$ and $\ol{\Ee}^{x_0}_{n,w_2}$ have the local universal property among, respectively, locally graded families of semistable $\Sp(2m,\CC)$, $\Ort(n,\CC)$ and $\SO(n,\CC)$-Higgs bundles of the appropiate topological type.
\epr  

\proof
Take any locally graded family $\Ee \to X \times Y$ of semistable $\Sp(2m,\CC)$, $\Ort(n,\CC)$ or $\SO(n,\CC)$-Higgs bundles. By definition of locally graded for every $y \in Y$ we have an open subset $U \subset Y$ containing $y$ and families $(\Vv_i, \varPhi_i)$ of Higgs bundles of rank $1$ and degree $0$ parametrized by $U$. Take the universal family $\Ee^{x_0}_{1,0} = (\Vv^{x_0}_{1,0},\varPhi_{1,0})$; for every $(\Vv_i, \varPhi_i)$ there exists $f_i : U \to T^*X$ such that $(\Vv_i,\varPhi_i) \sim_S f^*\Ee^{x_0}_{1,0}$. With the $f_i$ we can construct $f: U \to T^*X \times \dots \times T^*X$ such that the restriction $\Ee|_{X \times U}$ is S-equivalent to $f^* \wt{\Ee}^{x_0}_{2m}$, $f^* \mathring{\Ee}^{x_0}_{n,k,a}$ or $f^* \ol{\Ee}^{x_0}_{n,w_2}$.
\qed

\btm
\label{tm MmG}
There exist coarse moduli spaces $\Mm(\Sp(2m,\CC))$, $\Mm(\Ort(n,\CC))_{(k,a)}$ and $\Mm(\SO(n,\CC))_{w_2}$. Furthermore, there are isomorphisms
\begin{align*}
\wt{\xi}^{x_0}_{m} :\Sym^m (T^*X/\ZZ_2) \stackrel{\cong}{\lra} & \Mm(\Sp(2m,\CC)),
\\
\mathring{\xi}^{x_0}_{n,k,a} : \Sym^{(n-k)/2} (T^*X/\ZZ_2) \stackrel{\cong}{\lra} & \Mm(\Ort(n,\CC))_{k,a},
\\
\overline{\xi}^{x_0}_{m,0}: \quotient{T^*X \times \stackrel{m}{\dots} \times T^*X}{\Delta_m} \stackrel{\cong}{\lra} & \Mm(\SO(2m,\CC))_{0},
\\
\overline{\xi}^{x_0}_{n,w_2} : \Sym^{m'} (T^*X/\ZZ_2) \stackrel{\cong}{\lra} & \Mm(\SO(n,\CC))_{w_2}, 
\end{align*}
where $m' = m$ if $(n,w_2) = (2m+1,0)$, $m' = m-1$ if $(n,w_2) = (2m+1,1)$ and $m' = m-2$ if $(n,w_2) = (2m,1)$.
\etm

\proof
Since $T^*Z_{G,d}/\Gamma_{G,d}$ is an orbit space, the theorem follows from Proposition \ref{pr wtEe_2m has the local universal property} and Remark \ref{rm wtEe_z1 S-equivalent to wtEe_z2 if z1 and z2 are related by Gamma}. Note that the isomorphisms are defined as follows 
\[
\begin{array}{cccc}
\xi^{x_0}_{G,d} : & \quotient{T^*Z_{G,d}}{\Gamma_{G,d}}  & \lra &  \Mm(G)_d 
\\ 
& [z]_{\Gamma_{G,d}} &\longmapsto & \left[ \Ee^{x_0}_{G,d}|_{X \times z} \right]_S,
\end{array}
\]
where we write $\xi^{x_0}_{G,d}$ for $\wt{\xi}^{x_0}_{m}$, $\mathring{\xi}^{x_0}_{n,k,a}$ or $\overline{\xi}^{x_0}_{m,0}$ for the groups considered in the statement and the corresponding topological invariant. 
\qed

\brm
\nr{If an orbifold is defined as a global quotient $Z/\Gamma$, its cotangent orbifold bundle is the orbifold given by $T^*Z/\Gamma$, where the action of $\Gamma$ on $T^*Z$ is the action induced by the action of $\Gamma$ on $Z$.}

\nr{Let $G$ be $\GL(n,\CC)$, $\SL(n,\CC)$, $\PGL(n,\CC)$, $\Sp(2m,\CC)$, $\Ort(n,\CC)$ or $\SO(n,\CC)$. Denote by $\wt{M}(G)_d$ and $\wt{\Mm}(G)_d$ the orbifolds given by the quotients $Z_{G,d}/\Gamma_{G,d}$ and $T^*Z_{G,d}/\Gamma_{G,d}$. Then $\wt{\Mm}(G)_d$ is the cotangent orbifold bundle of $\wt{M}(G)_d$.}
\erm

\section{The Hitchin map}
\label{sc hitchin map}

\subsection{Description of the Hitchin map}
\label{sc description of the hitchin map}

Let $q_{n,1}, \dots , q_{n,n}$ be the standard basis for the invariant polynomials of a rank $n$ matrix. The Hitchin map is defined in \cite{hitchin-duke} by evaluating this basis on the Higgs field:
\beq
\label{eq Hitchin map of Higgs bundles}
\begin{array}{cccc}
b_{n,d} : & \Mm(\GL(n,\CC))_d  & \lra &  B_n = \bigoplus_{i = 1}^n H^0(X,\Oo) 
\\ 
& [(E,\Phi)]_S &\longmapsto & (q_{n,1}(\Phi), \dots , q_{n,n}(\Phi)).
\end{array}
\eeq
Since $H^0(X,\Oo) \cong \CC$ we have that $B_n \cong \CC^n$. Take $h = \gcd(n,d)$ and set $n' = n/h$ and $d' = d/h$. If $(E,\Phi)$ is a polystable Higgs bundle of rank $n$ and degree $d$, we have $(E,\Phi) \cong \bigoplus_{i = 1}^{h} (E_i,\Phi_i)$ where the $E_i$ are stable and $\Phi_i = \lambda_i \id_{E_i}$; moreover $\xi^{x_0}_{n',d'}((x_i,t_i)) = [(E_i,\Phi_i)]_S$, where $t_i = n' \cdot \lambda_i$. When we apply $q_{n,i}$ to $\Phi$ we obtain polynomials in $t_1, \dots, t_h$ 
\begin{align*}
q_{n,1}(\Phi) & =  \sum_{i = 1}^{h} t_i,
\\
& \vdots
\\
q_{n,n}(\Phi) & =  \left(\frac{1}{n'}t_1 \right)^{n'} \dots \left( \frac{1}{n'}t_h\right)^{n'}.
\end{align*}
The image of the Hitchin map is always contained in a subvariety of dimension $h$ which we denote by $B_{n,d}$. If $D_{\ol{\lambda}}$ is the diagonal matrix with eigenvalues $\ol{\lambda} = (\lambda_1, \dots, \lambda_h)$ we can construct the following bijective morphism
\beq
\label{eq definition of beta_nd}
\begin{array}{cccc}
\beta_{n,d} : & \Sym^h \CC  & \lra &  B_{n,d} 
\\ 
& [\ol{t}]_{\sym_h} = [t_1, \dots, t_h]_{\sym_h} &\longmapsto & \left( q_{n,1}(D_{\frac{1}{n'}\ol{t}}), \dots, q_{n,n}(D_{\frac{1}{n'}\ol{t}}) \right).
\end{array}
\eeq

Let us define the projection
\beq
\label{eq definition of pi^h}
\begin{array}{cccc}
\pi_{h} : & \Sym^h(T^*X)  & \lra & \Sym^h(\CC) 
\\ 
& [(x_1, t_1), \dots, (x_h, t_h) ]_{\sym_h} &\longmapsto & [t_1, \dots, t_h ]_{\sym_h}.
\end{array}
\eeq

It is immediate from the definitions of (\ref{eq Hitchin map of Higgs bundles}), (\ref{eq definition of beta_nd}) and (\ref{eq definition of pi^h}) that the diagram 
\beq
\label{eq commutative diagram of the Hitchin map for GLnC}
\xymatrix{
\Sym^h(T^*X)  \ar[rr]^{\pi_{h}} \ar[d]^{\cong}_{\xi^{x_0}_{n,d}} & & \Sym^h(\CC) \ar[d]^{\beta_{n,d}}_{1:1}
\\
\Mm(\GL(n,\CC))_d  \ar[rr]^{b_{n,d}} & & B_{n,d}
}
\eeq
commutes.

The Hitchin map for a classical structure group $G$ is defined in \cite{hitchin-duke} by evaluating the invariant polynomials for the adjoint representation of $G$ on its Lie algebra on the Higgs field:
\[
b_{G,d} : \Mm(G_d) \lra B(G,d).
\]

Let $C_{G,d} \subset T^*Z_{G,d}$ be the subset $(\{ x_0 \} \times \CC) \times \dots \times (\{ x_0 \} \times \CC)$, and consider the induced action of $\Gamma_{G,d}$ on $C_{G,d}$. Consider also the projection 
\[
\begin{array}{cccc}
\pi_{G,d} : & \quotient{T^*Z_{G,d}}{\Gamma_{G,d}}  & \lra &  \quotient{C_{G,d}}{\Gamma_{G,d}} 
\\ 
& [(x_1,t_1), \dots, (x_\ell,t_\ell)]_{\Gamma_{G,d}} &\longmapsto & [t_1, \dots, t_\ell]_{\Gamma_{G,d}}.
\end{array}
\]
In each case, the invariant polynomials for the adjoint representation of $G$ on its Lie algebra allow us to construct a bijective morphism
\beq
\label{eq definition of beta_Gd}
\beta_{G,d} : \quotient{C_{G,d}}{\Gamma_{G,d}} \stackrel{1:1}{\lra} B(G)_d,
\eeq
and we can extend (\ref{eq commutative diagram of the Hitchin map for GLnC}) to the rest of the classical complex Lie groups covered in this article:
\beq
\label{eq commutative diagram for Hitchin map for G}
\xymatrix{
\quotient{T^*Z_{G,d}}{\Gamma_{G,d}} \ar[rr]^{\pi_{G,d}} \ar[d]^{\cong}_{\xi^{x_0}_{G,d}} & & \quotient{C_{G,d}}{\Gamma_{G,d}} \ar[d]^{\beta_{G,d}}_{1:1}
\\
\Mm(G)_d \ar[rr]^{b_{G,d}}  & & B(G)_d .
}
\eeq

\subsection{The Hitchin fibres for the general linear case}
\label{sc Hitchin fibres for GLnC}

The set of tuples of the form
\beq
\label{eq generic element of Cn}
\overline{t}_g = (t_1, \dots, t_h),
\eeq
where $t_i \neq t_j$ if $i \neq j$, form a dense open subset of $\CC^h$. We call a point of $B_{n,d}$
{\it generic} if it is the image under $\beta_{n,d}$ of the $\sym_h$-orbit of some $\ol{t}_g$. An arbitrary point of $B_{n,d}$ is the image under $\beta_{n,d}$ of a $h$-tuple of the form  
\beq
\label{eq arbitrary element of Cn}
\overline{t}_a = (t_1, \stackrel{m_1}{\dots}, t_1, \dots, t_\ell, \stackrel{m_\ell}{\dots}, t_\ell),
\eeq
where $h = m_1 + \dots + m_\ell$.

\bpr
\label{pr fibres of pi_parentesis}
\[
\pi_{h}^{-1}([\overline{t}_g]_{\sym_h}) \cong X \times \stackrel{h}{\dots} \times X
\]
and
\[
\pi_{h}^{-1}([\overline{t}_a]_{\sym_h}) = \Sym^{m_1} X \times \stackrel{\ell}{\dots} \times \Sym^{m_\ell} X.
\]
\epr

\proof
The centralizer $Z_{\sym_h}(\overline{t}_g)$ of $\overline{t}_g$ in $\sym_h$ is trivial. Hence the centralizer of any element of $T^*X \times \dots \times T^*X$ of the form
\[
((x_1,t_1), \dots, (x_h,t_h)).
\]
is also trivial. If two tuples $((x_1,t_1), \dots, (x_h,t_h))$ and $((x'_1,t_1), \dots, (x'_h,t_h))$ lie in the same $\sym_h$-orbit, then they are related by the action of an element of $Z_{\sym_h}(\overline{t}_g)$. Since this group is trivial, it follows that $\pi_{h}^{-1}([\overline{t}_g]_{\sym_h})$ is given by the subset of $T^*X \times \stackrel{h}{\dots} \times T^*X$ which projects to $(t_1, \dots, t_h)$, which is isomorphic to $X \times \stackrel{h}{\dots} \times X$.

On the other hand, the centralizer of $\overline{t}_a$ in $\sym_h$ is 
\[
Z_{\sym_h}(\overline{t}_a) = \sym_{m_1} \times \dots \times \sym_{m_\ell},
\]
where the factor $\sym_{m_i}$ acts only on the entries of $\overline{t}_a$ equal to $t_i$. Two tuples of $T^*X \times  \dots  \times T^*X$ that project to $\overline{t}_a$ lie in the same $\sym_{h}$-orbit if and only if some element of $Z_{\sym_h}(\overline{t}_a)$ sends one tuple to the other.
Let us write $(T^*X \times \dots \times T^*X)_{\overline{t}_a}$ for the set of tuples as above that project to $\overline{t}_a$. We have
\[
\pi_{h}^{-1}([\overline{t}_a]_{\sym_h})\cong \quotient{(T^*X \times \stackrel{h}{\dots} \times T^*X)_{\overline{t}_a}}{Z_{\sym_m}(\overline{t}_a)},
\]
where the action of $\sym_{m_i}$ permutes the entries of a tuple that are pairs of the form $(x_{ij},t_{i})$.

We can easily see that $(T^*X \times \dots \times T^*X)_{\overline{t}_a} \cong X \times \dots \times X$ and then
\begin{align*}
\pi_{h}^{-1}([\overline{t}_a]_{\sym_h}) & \cong \quotient{(X \times \stackrel{h}{\dots} \times X)}{\sym_{m_1} \times \dots \times \sym_{m_\ell}}
\\
& \cong \Sym^{m_1} X \times \stackrel{\ell}{\dots} \times \Sym^{m_\ell} X.
\end{align*}
\qed

Note from (\ref{eq fibre of aj}) and (\ref{eq definition of varsigma_x0 1 0}) that $\Sym^h X$ is a fibration over $X$ with fibre $\PP^{h-1}$. 

\bco
\label{co GLnC Hitchin fibre}
The generic Hitchin fibre of $\Mm(\GL(n,\CC))_d \to B_{n,d}$ is the self-dual abelian variety $X \times \stackrel{h}{\dots} \times X$. The fibre over an arbitrary point of the base is $\Sym^{m_1} X \times \stackrel{\ell}{\dots} \times \Sym^{m_\ell} X$ which is a fibration over the self-dual abelian variety $X \times \stackrel{\ell}{\dots} \times X$ with fibre $\PP^{m_1 - 1} \times \stackrel{\ell}{\dots}\times \PP^{m_l -1}$. 
\eco

\subsection{The Hitchin fibres for the special linear and projective cases}
\label{sc Hitchin fibres of SLnC and PGLnC}

The Hitchin maps
\[
\hat{b}_n : \Mm(\SL(n,\CC)) \lra \hat{B}_n
\]
and 
\[
\check{b}_{n,\wt{d}} : \Mm(\PGL(n,\CC))_{\wt{d}} \lra \check{B}_{n,\wt{d}}
\]
are induced by the Hitchin map for Higgs bundles (\ref{eq Hitchin map of Higgs bundles}). The Hitchin base $\hat{B}_n$ is equal to the subvariety $B_n^{\tr = 0} = \bigoplus_{i=2}^n H^0(X,\Oo)$ of $B_n$, while $\check{B}_{n,\wt{d}} = B_n^{\tr = 0} \cap B_{n,d}$, where $\wt{d} = (d \mod n)$. Note that the groups $\SL(n,\CC)$ and $\PGL(n,\CC)$ are Langlands dual groups and that $\hat{B}_{n} = \check{B}_{n,0}$.

\blm
\label{lm generic fibre for SL and PGL in terms of GL}
Let $e$ be any element of $\hat{B}_n$ and $e'$ be any element $\check{B}_{n,\wt{d}}$. Then
\[
\hat{b}_{n}^{-1}(e) \cong b^{-1}_{n,0}(e) \cap \left( (\det, \tr)^{-1} (\Oo,0) \right) 
\]
and
\[
\check{b}_{n,\wt{d}}^{-1}(e') \cong \quotient{\left( b^{-1}_{n,d}(e') \cap \left( (\det, \tr)^{-1} (\Oo(x_0)^{\otimes d},0) \right) \right)}{\Pic^0(X)[n]} 
\]
\elm

\proof
Since $\hat{b}_n$ and $\check{b}_{n,\wt{d}}$ are induced by $b_{n,0}$ and $b_{n,d}$, this follows from the fact that 
\[
\Mm(\SL(n,\CC)) \cong (\det,\tr)^{-1}(\Oo,0) \subset \Mm(\GL(n,\CC))
\]
and 
\[
\Mm(\PGL(n,\CC)) \cong \quotient{(\det,\tr)^{-1}(\Oo(x_0)^{\otimes d},0)}{\Pic^0(X)[n]}.
\]
\qed

The generic elements $e_g$ and $e'_g$ of $\hat{B}_n$ and $\check{B}_{n,\wt{d}}$ come from a tuple $\overline{t}_g$ of the form (\ref{eq generic element of Cn}) such that $\sum t_i = 0$. Analogously, the arbitrary elements $e_a$ and $e'_a$ of $\hat{B}_n$ and $\check{B}_{n,\wt{d}}$ come from a tuple $\ol{t}_a$ of the form (\ref{eq arbitrary element of Cn}) such that $\sum_{i=1}^{h} t_i = 0$. 

Recall the definition of $A_\ell$ given in (\ref{eq definition of A_n}) and the isomorphism $u_\ell : Z_{\ell-1} \stackrel{\cong}{\lra} A_\ell$. Since $\Sym^m X \to X$ is a fibration over $X$, we can consider the restriction of $\Sym^{m_1} X \times \dots \times \Sym^{m_\ell} X$ to $A_\ell \subset X\times \stackrel{\ell}{\dots} \times X$ and the pull-back $u_\ell^* (\Sym^{m_1} X \times \dots \times \Sym^{m_\ell} X)|_{A_{\ell}}$.

\bpr
\label{pr arbitrary fibre of pi 2}
Let $e_g \in \hat{B}_n$ and $e'_g \in \check{B}_{n,\wt{d}}$ be generic elements. Then
\[
\hat{b}_{n}^{-1}(e_g) \cong A_n 
\]
and
\[
\check{b}_{n,\wt{d}}^{-1}(e'_g) \cong \quotient{A_h}{X[h]}. 
\]
Let $e_a \in \hat{B}_n$ and $e'_a \in \check{B}_{n,\wt{d}}$ be arbitrary elements. Then
\[
\hat{b}_{n}^{-1}(e_a) \cong u_\ell^* (\Sym^{m_1} X \times \dots \times \Sym^{m_\ell} X)|_{A_{\ell}} 
\]
and
\[
\check{b}_{n,\wt{d}}^{-1}(e'_a) \cong \quotient{u_\ell^* (\Sym^{m_1} X \times \dots \times \Sym^{m_\ell} X)|_{A_{\ell}}}{X[h]}.
\]
\epr

\proof
This follows from Lemma \ref{lm generic fibre for SL and PGL in terms of GL} and Corollary \ref{co GLnC Hitchin fibre}.
\qed

\bco
\label{co SLnC Hitchin fibre}
The generic Hitchin fibre of $\Mm(\SL(n,\CC)) \to \hat{B}_n$ is the abelian variety $X \times \stackrel{n-1}{\dots} \times X$. The arbitrary Hitchin fibre is a fibration over the abelian variety $X \times \stackrel{\ell - 1}{\dots} \times X$ with fibre $\PP^{m_1 - 1} \times \stackrel{\ell}{\dots}\times \PP^{m_l -1}$. 
\eco

\bpr
\label{pr PGLnC Hitchin fibre}
The generic fibre of the Hitchin fibration $\Mm(\PGL(n,\CC))_{\wt{d}} \to \check{B}_{n,\wt{d}}$ is the abelian variety $X \times \stackrel{h-1}{\dots} \times X$. 

The arbitrary fibre of the Hitchin fibration $\Mm(\PGL(n,\CC))_{\wt{d}} \to \check{B}_{n,\wt{d}}$ is a fibration over $X \times \stackrel{\ell-1}{\dots} \times X$ with fibre $(\PP^{m_1-1}\times \dots \times \PP^{m_\ell-1})/X[r]$, where $r = \gcd(h, m_1, \dots, m_\ell)$. 
\epr

\proof
The action of $X[h]$ on the fibration $\Sym^{m_1}X \times \stackrel{\ell}{\dots} \times  \Sym^{m_\ell} X$ gives a weighted $(m_1, \dots, m_\ell)$-action of $X[h]$ on $X \times \stackrel{\ell}{\dots} \times X$, the base of the fibration. Let us set $r = \gcd(h,m_1,\dots,m_\ell)$. The subgroup $X[r] \subset X[h]$ acts on $\Sym^{m_1}X \times \stackrel{\ell}{\dots} \times  \Sym^{m_\ell} X$ and it acts trivially on the base $X \times \stackrel{\ell}{\dots} \times X$. The quotient of $\Sym^{m_1}X \times \stackrel{\ell}{\dots} \times  \Sym^{m_\ell} X$ by the action of $X[r]$ is a fibration over $X \times \stackrel{\ell}{\dots} \times X$ with fibre $(\PP^{m_1-1} \times \dots \PP^{m_\ell-1})/X[r]$, where the action of $X[r]$ is described in Remark \ref{rm the action of X r on PP}.

The quotient of $\Sym^{m_1}X \times \stackrel{\ell}{\dots} \times  \Sym^{m_\ell} X$ by $X[h]$ is equivalent to the quotient of $(\Sym^{m_1}X \times \stackrel{\ell}{\dots} \times  \Sym^{m_\ell} X)/X[r]$ by $X[h]/X[r]$. Since $X[h]/X[r] \cong X[h/r]$ and the weighted action of $X[r]$ is trivial, the weighted $(m_1,\dots,m_\ell)$-action of $X[h]/X[r]$ on $X\times \stackrel{\ell}{\dots} \times X$ is equivalent to the weighted $(m_1/r,\dots,m_\ell/r)$-action of $X[h/r]$ on $X \times \stackrel{\ell}{\dots} \times X$. By Lemma \ref{lm the the weighted action is free iff r=1}, this action is free since $\gcd(h/r,m_1/r,\dots, m_\ell/r) = 1$. As a consequence, $(\Sym^{m_1} X \times \dots \Sym^{m_\ell}X)/X[h]$ is a fibration with fibre $(\PP^{m_1-1}\times \dots \times \PP^{m_\ell-1})/X[r]$ over $(X \times \stackrel{\ell}{\dots} \times X)/X[h/r]$. 

The weighted $(m_1/r,\dots,m_\ell/r)$-action of $X[h/r]$ restricted to $A_{\ell}$ is equivariant under $u_{\ell}$ to the weighted $(m_1/r,\dots,m_{\ell-1}/r)$-action of $X[h/r]$ on $X \times \stackrel{\ell-1}{\dots} \times X$, so the pull-back $u_\ell^* (\Sym^{m_1} X \times \dots \times \Sym^{m_\ell} X)|_{A_{\ell}}/X[h]$ is a fibration with fibre $(\PP^{m_1-1}\times \dots \times \PP^{m_\ell-1})/X[r]$ over $(X \times \stackrel{\ell-1}{\dots} \times X)/X[h/r]$.

Finally, by Lemma \ref{lm the quotient of the weighted action is again X times dots times X}, the quotient $(X \times \stackrel{\ell-1}{\dots} \times X)/X[h/r]$ is isomorphic to $X \times \stackrel{\ell-1}{\dots} \times X$.
\qed

\brm
\label{rm SLnC and PGLnC Hitchin fibres}
\nr{The generic fibre of the Hitchin fibration $\Mm(\PGL(n,\CC))_0 \to \check{B}_{n,0} = \hat{B}_{n}$ and the corresponding fibre of the Hitchin fibration $\Mm(\SL(n,\CC)) \to \hat{B}_n$ are isomorphic to $X \times \stackrel{n-1}{\dots} \times X$, which is a self-dual abelian variety.} 

\nr{The arbitrary fibre of the Hitchin fibration $\Mm(\PGL(n,\CC))_0 \to \check{B}_{n,0}$ and the corresponding fibre of the Hitchin fibration $\Mm(\SL(n,\CC)) \to \hat{B}_n$ are fibrations over $X \times \stackrel{\ell-1}{\dots} \times X$, which is a self-dual abelian variety.}
\erm

\subsection{The Hitchin fibres for the symplectic and orthogonal cases}
\label{sc Hitchin fibres for Sp2mC OnC and SOnC}

Consider the following projection induced by the natural projection $T^*X \cong X \times \CC \to \CC$,
\[
\pi_{[m']} : \Sym^{m'}(T^*X/\ZZ_2) \lra \Sym^{m'}(\CC/\ZZ_2).
\]
We set:
\bit
\item $m' = m$ if $G = \Sp(2m,\CC)$,
\item $m' = (n-k)/2$ if $G = \Ort(n,\CC)$ and the topological invariant is $(k,a)$,
\item $m' = m$ if $G = \SO(2m+1,\CC)$ and $w_2 = 0$,
\item $m' = m-1$ if $G = \SO(2m+1,\CC)$ and $w_2 = 1$,
\item $m' = m-2$ if $G = \SO(2m,\CC)$ and $w_2 = 1$.
\eit 
By (\ref{eq commutative diagram for Hitchin map for G}), we have the following commuting diagram for the cases considered above:
\[
\xymatrix{
\quotient{Z_{G,d}}{\Gamma_{G,d}} \ar[d]^{\cong}_{\xi^{x_0}_{G,d}} \ar[rr]^{\pi_{[m']}} & & \quotient{C_{G,d}}{\Gamma_{G,d}} \ar[d]^{\beta_{G,d}}_{1:1}
\\
\Mm(G)_d \ar[rr]^{b_{G,d}}  & & B(G)_d.
}
\]

The elements of $\CC \times \stackrel{m'}{\dots} \times \CC$ of the form
\[
\overline{t}_g = (t_1, \dots, t_{m'}),
\]
where $t_i \neq \pm t_j$ if $i \neq j$ and for every $i$ we have $t_i \neq 0$, form a dense open subset. For each $(G, d)$, we call a point of $B(G)_{d}$ {\it generic} if it is the image under $\beta_{G,d}$ of $[\ol{t}_g]_{\sym_{m'}}$.
An arbitrary element of $B(G)_d$ is the image under $\beta_{G,d}$ of the $\Gamma_m$-orbit of the following tuple  
\[
\overline{t}_a = (0, \stackrel{m_0}{\dots}, 0, t_1, \stackrel{m_1}{\dots}, t_1, \dots,t_\ell, \stackrel{m_\ell}{\dots}, t_\ell),
\]
where $t_i \neq 0$, $t_i \neq \pm t_j$ if $i \neq j$ and $m' = m_0 + m_1 + \dots + m_\ell$.

The following proposition applies to all the situations listed above.

\bpr
\label{pr fibre of pi_square}
\[
\pi_{[m']}^{-1}([\overline{t}_g]_{\Gamma_{m'}}) \cong X \times \stackrel{m'}{\dots} \times X
\]
and
\[
\pi_{[m']}^{-1}([\overline{t}_a]_{\Gamma_{m'}}) \cong \PP^{m_0} \times  \Sym^{m_1} X \times \dots \times \Sym^{m_\ell} X.
\]
\epr

\proof
Since $t_i \neq -t_i$ and $t_i \neq \pm t_j$ for every $i$, $j$ such that $i \neq j$, the stabilizer in $\Gamma_{m'}$ of $\overline{t}_g$ is trivial and then the stabilizer of every tuple of the form
\[
((x_1,t_1), \dots, (x_{m'},t_{m'}))
\]
is trivial too. This implies that every such tuple is uniquely determined by the choice of $(x_1, \dots, x_{m'})$, and then $\pi_{[m']}^{-1}([\overline{t}_g]_{\Gamma_{m'}})$ is isomorphic to $X \times \stackrel{m'}{\dots} \times X$.

Since the stabilizer in $\Gamma_{m'}$ of $\overline{t}_a$ is 
\[
Z_{\Gamma_{m'}}(\overline{t}_a) = \Gamma_{m_0} \times \sym_{m_1} \times \dots \times \sym_{m_\ell},
\]
we have
\begin{align*}
\pi_{[m']}^{-1}([\overline{t}_a]_{\Gamma_{m'}}) & \cong (X\times \stackrel{m'}{\dots} \times X) / Z_{\Gamma_{m'}}(\overline{t}_a) 
\\
& \cong \Sym^{m_0} (X/\ZZ_2) \times \Sym^{m_1} X \times \dots \times \Sym^{m_\ell} X.
\end{align*}
Note that $\Sym^{m_0}(X/\ZZ_2) \cong \Sym^{m_0} \PP^1 \cong \PP^{m_0}$.
\qed

Recall that $\Sym^{m} X$ is a projective bundle over $X$.

\bco
\label{co Sp2mC Hitchin fibre}
The generic fibre of the Hitchin map for $\Mm(\Sp(2m,\CC))$ is the abelian variety $X \times \stackrel{m}{\dots} \times X$. The Hitchin fibre over an arbitrary element is a fibration over $X \times \stackrel{\ell}{\dots} \times X$ with fibre $\PP^{m_0} \times \PP^{m_1-1} \times \dots \times \PP^{m_\ell-1}$.
\eco

\bco
\label{co SOnC Hitchin fibre}
The generic fibre of the Hitchin fibration for $\Mm(\SO(2m+1,\CC))_0$ is isomorphic to $X \times \stackrel{m}{\dots} \times X$. The Hitchin fibre over an arbitrary element is a fibration over $X \times \stackrel{\ell}{\dots} \times X$ with fibre $\PP^{m_0} \times \PP^{m_1-1} \times \dots \times \PP^{m_\ell-1}$.
\eco

\brm
\label{rm SOnC and Sp2mC Hitchin fibres}
\nr{Note that the groups $\Sp(2m,\CC)$ and $\SO(2m+1,\CC)$ are Langlands dual groups and the Hitchin bases for $\Mm(\Sp(2m,\CC))$ and $\Mm(\SO(2m+1,\CC))_0$ are the same. The generic fibre of the Hitchin fibration for $\Mm(\Sp(2m,\CC))$ and the corresponding fibre of the Hitchin fibration for $\Mm(\SO(2m+1,\CC))_0$ are isomorphic to $X\times \stackrel{m}{\dots} \times X$, which is a self-dual abelian variety.} 

\nr{The fibre of the Hitchin fibration for $\Mm(\Sp(2m,\CC))$ over an arbitrary point of the Hitchin base and the fibre of the Hitchin fibration for $\Mm(\SO(2m+1,\CC))_0$ over the same point are fibrations over $X\times \stackrel{\ell}{\dots} \times X$, which is a self-dual abelian variety. The fibres are isomorphic to $\PP^{m_0} \times \PP^{m_1 - 1} \times \dots \times \PP^{m_\ell-1}$.} 
\erm


Recall the finite group $\Delta_m$ defined in (\ref{eq definition of Delta_m}). Thanks to the natural projection $T^*X \cong X \times \CC \to \CC$ we define
\[
\ol{\pi}_{m} : \quotient{T^*X \times \stackrel{m}{\dots} \times T^*X}{\Delta_m} \lra \quotient{\CC \times \stackrel{m}{\dots} \times \CC}{\Delta_m}
\] 
Let us denote $b_{G,d}$, $B_{G,d}$ and $\beta_{G,d}$ by $\ol{b}_{2m,0}$, $\ol{B}_{2m,0}$ and $\ol{\beta}_{2m,0}$ when $G = \SO(2m,\CC)$ and $w_2 = 0$. By (\ref{eq commutative diagram for Hitchin map for G}) and Theorem \ref{tm MmG} we have the following commuting diagram for the Hitchin fibration associated to $\Mm(\SO(2m,\CC))_0$
\[
\xymatrix{
\quotient{T^*X \times \stackrel{m}{\dots} \times T^*X}{\Delta_m} \ar[rr]^{\ol{\pi}_{m}} \ar[d]^{\cong}_{\ol{\xi}^{x_0}_{2m,0}} & & \quotient{\CC \times \stackrel{m}{\dots} \times \CC}{\Delta_m} \ar[d]^{\ol{\beta}_{2m,0}}_{1:1}
\\
\Mm(\SO(2m,\CC))_0 \ar[rr]^{\ol{b}_{2m,0}}  & & \ol{B}_{2m,0}.
}
\]

Using the map
\[
\begin{array}{cccc}
r : & X \times \stackrel{m}{\dots} \times X  & \lra &  X \times \stackrel{m}{\dots} \times X 
\\ 
& (x_1, x_2, \dots, x_m) &\longmapsto & (-x_1, x_2, \dots, x_m),
\end{array}
\]
we define the $r$-action of $\sym_m$ on $X \times \stackrel{m}{\dots} \times X$ as follows. Suppose that for $\sigma \in \sym_m$ we denote by $f_{\sigma}$ the permutation of $X \times \stackrel{m}{\dots} \times X$ associated to $\sigma$. We define the $r$-action of $\sigma$ to be the morphism $r \circ f_{\sigma} \circ r$.

Consider the elements of $\CC^m$ of the form
\[
\overline{t}_{g} = (t_1,\dots,t_m)
\]
where $t_i \neq 0$ and $t_i \neq \pm t_j$. We call a point of $\ol{B}_{2m,0}$ {\it generic} if it is the image under $\ol{\beta}_{2m,0}$ of some $\Delta_m$-orbit of $\ol{t}_g$.
There is a special set of points of $\ol{B}_{2m,0}$ that come from tuples of the form
\[
\overline{t}_{s_1} = (-t_1, t_1, \stackrel{m_1 - 1}{\dots}, t_1, t_2, \stackrel{m_2}{\dots}, t_2,\dots, t_{\ell}, \stackrel{m_{\ell}}{\dots}, t_{\ell}),
\]
where $t_i \neq \pm t_j$ if $i \neq j$, and for every $i$ we have $t_i \neq 0$ and $m_i$ even. If our point of $\ol{B}_{2m,0}$ is given by a tuple different from $\overline{t}_{g}$ and $\ol{t}_{s_1}$ we can always find a representative of the $\Delta_m$-orbit with the form
\[
\ol{t}_{s_2} = (0, \stackrel{m_0}{\dots}, 0, t_1, \stackrel{m_1}{\dots},t_1, \dots, t_\ell, \stackrel{m_\ell}{\dots}, t_\ell ),
\]
where $t_i \neq \pm t_j$ if $i \neq j$ and for every $i > 0$ we have $t_i \neq 0$.

\bpr
\[
\overline{\pi}_m^{-1}([\overline{t}_{g}]_{\Delta_m}) = X \times \stackrel{m}{\dots} \times X,
\]
\[
\overline{\pi}_m^{-1}([\overline{t}_{s_1}]_{\Delta_m}) = \Sym^{m_1} X \times \Sym^{m_2} X  \dots \times \Sym^{m_\ell} X
\]
and
\[
\ol{\pi}_m^{-1}([\ol{t}_{s_2}]_{\Delta_m}) = \left( \quotient{(X \times \dots \times X)}{\Delta_{m_0}} \right) \times \Sym^{m_1} X \times \dots \times \Sym^{m_\ell} X.
\]
\epr 

\proof
The first result follows from the observation that the stabilizer of $\overline{t}_{g}$ is trivial.

The stabilizer of $\overline{t}_{s_1}$ is
\[
Z_{\Delta_m}(\overline{t}_{s_1}) = Z_{\Delta_{m_1}} ((-t_1, t_1, \stackrel{m_1 - 1}{\dots} t_1)) \times \sym_{m_2} \times \dots \times \sym_{m_{\ell}}. 
\]
We can check that $Z_{\Delta_{m_1}} ((-t_1, t_1, \stackrel{m_1 - 1}{\dots} t_1))$ is given by the elements $\overline{c}\sigma$ of $\Delta_{m_i}$ such that $\sigma$ sends the first entry of $(-t_1, t_1, \stackrel{m_1 - 1}{\dots} t_1)$ to the $i$-th entry and $\overline{c}$ inverts the first and the $i$-th entry. This shows that the action of $Z_{\Delta_{m_1}} ((-t_1, t_1, \stackrel{m_1 - 1}{\dots} t_1))$ on $X \times \dots \times X$ is equivalent to the $r$-action of the symmetric group. One can check that the quotient of $X \times \dots \times X$ under this action is isomorphic to the symmetric product of the curve, so
\[
\quotient{X \times \dots \times X}{Z_{\Delta_m}(\overline{t}_{s_1})} \cong \Sym^{m_1} X \times \Sym^{m_2} X \times\dots \times \Sym^{m_\ell} X.
\]

The last statement follows from the fact that the stabilizer of $\overline{t}_{s_2}$ is 
\[
Z_{\Delta_m}(\ol{t}_{a}) = \Delta_{m_0} \times \sym_{m_1} \times \dots \times \sym_{m_\ell}.
\]
\qed


\bco
\label{co SO2mC Hitchin fibre}
The generic fibre of the Hitchin fibration $\Mm(\SO(2m,\CC))_0 \to \ol{B}_{2m,0}$ is the self-dual abelian variety $X \times \stackrel{m}{\dots} \times X$. The fibre over an arbitrary point of $\ol{B}_{2m,0}$ is a fibration over the self-dual abelian variety $X \times \dots \times X$ with fibre 
\[
\PP^{m_1-1} \times \dots \times \PP^{m_\ell-1},
\]
or
\[
(X \times \stackrel{m_0}{\dots} \times X)/\Delta_{m_0} \times \PP^{m_1-1} \times \dots \times \PP^{m_\ell-1}. 
\]
\eco

\brm
\nr{We observe that $(X \times \stackrel{m_0}{\dots} \times X)/\Delta_{m_0}$ is isomorphic to the moduli space $M(\SO(2m_0,\CC))_0$. By \cite{friedman&morgan} and \cite{looijenga}, this variety is isomorphic to the quotient of the weighted projective space $\WW\PP(1,1,1,2,\stackrel{m_0 - 2}{\dots},2)$ by a finite group.}
\erm

\appendix

\section{Some results on abelian varieties}
\label{sc Some results on abelian varieties}

In this appendix, we prove two lemmas; the first is certainly well known, the second presumably so, but we have been unable to find a suitable reference.

Let $A$ be an abelian variety and let $a_0 \in A$ be the trivial element. For any integers $m_1, \dots, m_\ell$, the formula 
\[
a' \cdot (a_1, \dots, a_\ell) = (a_1 + m_1 a' , \dots, a_\ell + m_\ell a') 
\]
defines the {\it weighted $(m_1, \dots, m_\ell)$-action} of $A$ on $A \times \stackrel{\ell}{\dots} \times A$.

\blm
\label{lm the the weighted action is free iff r=1}
Let $m_1, \dots , m_\ell$ be integers and let $h$ be a positive integer. Write $r$ for $\gcd(h,m_1,\dots,m_\ell)$. The weighted $(m_1, \dots, m_\ell)$-action of $A[h]$ on $A \times \stackrel{\ell}{\dots} \times A$ is free if and only if $r = 1$. 
\elm

\proof 
Suppose $a' \in A[h]$ and $a' \cdot (a_0, \dots, a_0) = (a_0, \dots, a_0)$. Then $m_i a' = a_0$ for every $i$. This implies that $a'$ is a $m_i$-torsion element for every $i$. On the other hand, if there exists $a' \in A[h] \cap \bigcap_{i} A[m_i]$ different from $a_0$,  the $(m_1, \dots, m_\ell)$-weighted action of $a'$ is trivial and therefore the action of $A[h]$ is not free. Thus the action is free if and only if the subgroup $A[h] \cap \bigcap_i A[m_i]$ is trivial. 

It is easy to see that $A[n_1] \cap A[n_2] = A[r']$ where $r' = \gcd(n_1,n_2)$. It follows by induction that 
$
A[h] \cap \bigcap_i A[m_i] = A[r]
$. The result follows.
\qed

For every positive integer $h$, we have an exact sequence
\[
\xymatrix{ 0 \ar[r] & A[h] \ar[r] & A \ar[r]^{\mult_h} & A \ar[r] & 0, 
}
\]
where $\mult_h(a) = ha = a + \stackrel{h}{\dots} + a$. This induces an isomorphism
\beq
\label{eq definition of wt f_h} 
\wt{\mult_h} : \quotient{A}{A[h]} \stackrel{\cong}{\lra} A.
\eeq

\blm
\label{lm the quotient of the weighted action is again X times dots times X}
Consider the weighted $(m_1,\dots,m_{\ell})$-action of $A[h]$ on $(A \times \stackrel{\ell}{\dots} \times A)$. Then
\[
\quotient{(A \times \stackrel{\ell}{\dots} \times A)}{A[h]} \cong A \times \stackrel{\ell}{\dots} \times A. 
\]
\elm 

\proof
First we treat the case where $\gcd(m_1,h) = 1$. Let $p$, $q$ be integers such that $pm_1+qh=1$. Since $\gcd(p,h)=1$, the action is equivalent to the action of $A[h]$ with weights $(1, pm_2, \dots, pm_\ell)$. For this action, consider the morphism
\[
\begin{array}{ccc}
\quotient{(A \times \stackrel{\ell}{\dots} \times A)}{A[h]}  & \lra
&  \quotient{A}{A[h]} \times A \times \stackrel{\ell - 1}{\dots} \times A
\\ \left[(a_1,\dots,a_\ell)\right]_{A[h]} &\longmapsto
&\left( [a_1]_{A[h]}, a_2 - pm_2 a_1, \dots, a_\ell - pm_\ell a_1 \right).
\end{array}\]
This is in fact an isomorphism of abelian varieties since it has an inverse
\[
\begin{array}{ccc}
\quotient{A}{A[h]} \times A \times \stackrel{\ell - 1}{\dots} \times A   & \lra
&  \quotient{(A \times \stackrel{\ell}{\dots} \times A)}{A[h]}  
\\ \left( [a'_1]_{A[h]}, a'_2, \dots, a'_\ell \right) &\longmapsto
&[(a'_1, pm_2 a'_1+a'_2, \dots,pm_\ell a'_1+ a'_\ell)]_{A[h]}.
\end{array}\]
It follows from (\ref{eq definition of wt f_h}) that the lemma is true when $\gcd(m_1,h) = 1$.

For the general case, let $r_1:=\gcd(m_1,h)$ and write
\beq\label{eq: hr1}
\quotient{(A \times A \times \stackrel{\ell - 1}{\dots} \times A)}{A[h]} = \quotient{\left( (A \times A \times \stackrel{\ell - 1}{\dots} \times A)/A[r_1]  \right)}{A[h/r_1]}.
\eeq
Since $r_1$ divides $m_1$, the action of $A[r_1]$ on the first factor is trivial, so
\[
\quotient{(A \times A \times \stackrel{\ell - 1}{\dots} \times A)}{A[r_1]} \cong A \times \quotient{(A \times \stackrel{\ell - 1}{\dots} \times A)}{A[r_1]}.
\]
When $\ell=1$ the result follows from (\ref{eq definition of wt f_h}), so we can suppose inductively that the lemma is true for $A \times \stackrel{\ell - 1}{\dots} \times A$. Then
\[
A \times \quotient{(A \times \stackrel{\ell - 1}{\dots} \times A)}{A[r_1]} \cong A \times A \times \stackrel{\ell - 1}{\dots} \times A.
\]
This implies by \eqref{eq: hr1} that the original quotient is isomorphic to the quotient by an action of $A[h/r_1]$ whose first weight is $m_1/r_1$. Since $\gcd(m_1/r_1,h/r_1)=1$, the first part of the proof completes the induction and hence the entire proof.\qed

\end{document}